\documentclass[11pt,a4paper]{article}
\usepackage[margin=1in]{geometry}
\usepackage{tabularx,enumerate,enumitem,authblk,soul}
\usepackage[dvipsnames]{xcolor}
\usepackage{amsmath,amssymb,dsfont,subcaption}
\usepackage{mathtools, graphicx,color,dsfont,amsthm,bbm, mathrsfs}
\usepackage{lipsum,bm,comment}
\usepackage{natbib}
\usepackage{longtable}
\usepackage{array}
\usepackage{caption}
\usepackage[colorlinks=true,linkcolor=blue,anchorcolor=blue,citecolor=blue,filecolor=black,menucolor=black,runcolor=black,urlcolor=black]{hyperref}
\usepackage{cleveref}

\usepackage{lscape}

\theoremstyle{plain}
\newtheorem{Theor}{Theorem}[section]
\newtheorem{Lem}[Theor]{Lemma}
\newtheorem{Asmp}[Theor]{Assumption}
\newtheorem{Prop}[Theor]{Proposition}
\newtheorem{Corol}[Theor]{Corollary}

\theoremstyle{definition}
\newtheorem{Defin}[Theor]{Definition}
\newtheorem{Example}[Theor]{Example}

\newcommand{\B}{\mathbb{B}}

\newcommand{\E}{\mathbb{E}}

\newcommand{\I}{\mathbb{I}}

\renewcommand{\P}{\mathbb{P}}

\newcommand{\R}{\mathbb{R}}
\renewcommand{\S}{\mathbb{S}}

\renewcommand{\AA}{\mathcal{A}}

\newcommand{\FF}{\mathcal{F}}
\newcommand{\GG}{\mathcal{G}}

\newcommand{\LL}{\mathcal{L}}
\newcommand{\MM}{\mathcal{M}}
\newcommand{\NN}{\mathcal{N}}
\newcommand{\OO}{\mathcal{O}}
\newcommand{\PP}{\mathcal{P}}
\newcommand{\QQ}{\mathcal{Q}}
\newcommand{\RR}{\mathcal{R}}

\newcommand{\XX}{\mathcal{X}}

\newcommand{\PPP}{\mathscr{P}}

\newcommand{\RRR}{\mathscr{R}}

\newcommand{\rr}{\mathfrak{r}}
\newcommand{\rR}{\mathfrak{R}}

\newcommand{\Fpr}{F_P^{\rr}}
\newcommand{\Qpr}{Q_P^{\rr}}

\newcommand{\Fd}{{\normalfont \textbf{\textrm{F}}}}

\def\1{{\mathbb I}}
\newcommand{\Ind}[1]{\1\left[#1\right]}
\newcommand{\ie}{\textit{i.e.} }

\newcommand{\ps}[2]{\left\langle #1,#2 \right\rangle}
\newcommand{\sm}{\setminus}

\DeclareMathOperator{\diag}{diag}

\newcommand{\ve}{\varepsilon}
\newcommand{\Om}{\Omega}

\newcommand{\spa}{\quad\quad}

\usepackage{scalerel,stackengine}
\stackMath
\newcommand\wwidehat[1]{%
	\savestack{\tmpbox}{\stretchto{%
			\scaleto{%
				\scalerel*[\widthof{\ensuremath{#1}}]{\kern-.6pt\bigwedge\kern-.6pt}%
				{\rule[-\textheight/2]{1ex}{\textheight}}%WIDTH-LIMITED BIG WEDGE
			}{\textheight}% 
		}{0.5ex}}%
	\stackon[1pt]{#1}{\tmpbox}%
}

\newcommand{\RN}[1]{%
	(\textup{\uppercase\expandafter{\romannumeral#1}})%
}
\newcolumntype{H}{>{\setbox0=\hbox\bgroup}c<{\egroup}@{}}

\newenvironment{Proof}[1]
{
	\begin{proof}[Proof of #1]
}
{
	\end{proof}
}

\allowdisplaybreaks

\title{Regularized geometric quantiles and universal\\ linear distribution functionals}
\author{Dimitri Konen$^*$}
\author{Gilles Stupfler$^\dag$}
\affil{$^*$% Department of Pure Mathematics and Mathematical Statistics,
University of Cambridge, $^\dag$University of Angers
}
\date{\today}

\begin{document}
	
	\maketitle
	
\begin{abstract}
Geometric quantiles are popular location functionals to build rank-based statistical procedures in multivariate settings. They are obtained through the minimization of a non-smooth convex objective function. As a result, the singularity of the directional derivatives leads to numerical instabilities and %bad 
poor sample properties as well as surprising `phase transitions' from empirical to population distributions. To solve these issues, we introduce a regularized version of geometric distribution functions and quantiles that are provably close to the usual geometric concepts and share their  qualitative properties, both in the empirical and continuous case, while allowing for a much broader applicability of asymptotic results without any moment condition. We also show that any linear assignment of probability measures (such as the univariate %cdf
distribution function), that is also translation- and orthogonal-equivariant, necessarily coincides with one of our regularized geometric distribution functions.
 % Lack of a canonical ordering $\leadsto$ many depth and quantiles concepts, with no particular one that is the canonical choice. Novel approach: multivariate distribution functions can be fully described by their symmetry properties.
%		While this PDE equality holds in the sense of distributions for an arbitrary probability measure, when it admits a density we provide sufficient conditions to ensure that $f_P=\LL_d(R_P)$ holds pointwise. Surprisingly, the reconstruction procedure is of a local nature when $d$ is odd, and of a non-local nature when $d$ is even. We give examples of the reconstruction in $\R^2$ and $\R^3$. We use our results to characterise the regularity of depth contours. We conclude the paper with a partial counterpart to the non-localisability in even dimensions.
%		\medskip
%		\textbf{Keywords}: Huber's contamination model; Local differential privacy; Minimax optimality.
	\end{abstract}

%	\maketitle
%	
	
%	{
%			\hypersetup{linkcolor=black}
%		\tableofcontents
%	}
%	
\setcounter{tocdepth}{2}
{
\hypersetup{linkcolor=black}
\tableofcontents
}

\section{Introduction} \label{sec:Introduction}

\subsection{Background} \label{sec:Introduction:Background}

Geometric quantiles, introduced in \citet*{Cha1996}, are location functionals that provide a center-outward ordering of points in multivariate datasets, which can then be used to build multivariate rank and sign tests (\cite*{Mottonen97}), perform quantile regression (\cite*{Cha2003, ChoCha19}), and conduct inference in various statistical applications such as supervised classification (\cite*{PaiVanB2012}), testing for symmetry (\cite*{KolLi98}), and detection of outliers (\cite*{DanSer2010}, %and  analysis of extremes (\cite*{Daoetal2019}), 
to name a few. Geometric quantiles are very popular in practice because of their conceptual simplicity and broad applicability: they are defined for \emph{any} probability measure on $\R^d$, $d\ge 1$, without moment condition, and are obtained through a convex minimization problem which makes them easily computable even in very high dimensions (\cite*{Frietal2012}). This is in sharp contrast with the measure transportation concept of distribution function and quantiles  introduced in \cite*{HallEspagne21}, that is computationally intensive even for low-dimensional datasets and is well-defined only for probability measures with a density---their concept can be extended to empirical probability measures at the cost of important technicalities including smooth interpolation between atoms of the distribution in order to preserve the cyclical monotonicity of the corresponding distribution function; see also \cite*{KonHal24} for an overview of measure transporation and geometric quantiles. The study of the geometric median started in the 1980's (\cite*{MilDuc1987, VarZha2000, MagTyl2008, Motetal2010}), whereas the main theoretical properties of geometric quantiles were investigated in \cite{Cha1996}, \cite{Kol1997}, \cite*{Ser2010}, \cite*{GirStu2015,GirStu2017}, and \cite*{KonPai1}, and include the study of existence, uniqueness and convexity of the geometric quantile map and distribution function. \smallskip

Given a probability measure $P$ on $\R^d$, a geometric quantile of order $\alpha\in [0,1)$ in direction $u$ in the unit sphere $\S^{d-1}=\{x\in\R^d : \|x\|^2=\ps{x}{x}=1\}$ is defined as a minimizer over $\R^d$ of 
\begin{equation}\label{eq:DefinMalphaU}
x\mapsto 
M_{\alpha,u}^P(x)
\equiv
\int_{\R^d} \big(\|z-x\|-\|z\|\big)\, dP(z)-\ps{\alpha u}{x} 
.
\end{equation}
When $P$ is not supported on a single line of $\R^d$, then for any given $\alpha\in [0,1)$ and $u\in\S^{d-1}$, the map $M_{\alpha,u}^P$ is strictly convex so that there is a unique geometric quantile, $\mu_{\alpha,u}(P)$ say. In particular, the \emph{geometric quantile map} $Q_P:\B^d\to\R^d,\ \alpha u\mapsto \mu_{\alpha,u}(P)$ is well-defined on the open unit ball $\B^d$ of $\R^d$. In addition, when $P$ has no atoms it can be shown---see, e.g., Theorem 6.1 in \cite{KonPai1}---that $Q_P$ is in fact a homeomorphism between $\B^d$ and $\R^d$ with inverse $F_P=Q_P^{-1}:\R^d\to\B^d$ given by 
\begin{equation}\label{eq:DefinGeometricCdf}
F_P(x)
\equiv
\int_{\R^d} \frac{x-z}{\|x-z\|}\I[z\neq x]\, dP(z)
,\spa 
\forall\ x\in\R^d 
;
\end{equation}
throughout, $\I[\cdot]$ denotes the indicator function. This fact relies on observing that $\Fpr$ and $M_{\alpha,u}^{\rr,P}$ are related through $\nabla M_{\alpha,u}^P = F_P - \alpha u$. While $P$ needs to be non-atomic and not supported on a single line for $F_P$ to be the inverse of $Q_P$, the map $F_P$ itself is well-defined for \emph{any} probability measure $P$. Geometric quantiles are indexed by $\alpha u\in \B^d$, where $\alpha$ is a measure of the outlyingness of $Q_P(\alpha u)$---small values of $\alpha$ correspond to central quantiles whereas $\alpha$ large provides extreme quantiles approximately in direction $u$. The open unit ball yields a multivariate analogue of the interval $(0,1)$ in which the usual univariate distribution function takes its values, with $1/2$ yielding the innermost quantile. As such, for any $d\ge 1$, the vector field $F_P$ from (\ref{eq:DefinGeometricCdf}) is referred to as the \emph{geometric distribution function} of $P$. Empirical geometric quantiles satisfy a Bahadur-Kiefer representation (\cite*{Kol94}), enjoy attractive robustness properties (\cite*{LopRou1991, KonPai2025a, KonPai25}), have explicitly available extreme behavior (\cite*{GirStu2015, GirStu2017, PaiVir2021}), can be used to metrize the topology of weak convergence and define a multivariate analogue of Kolmogorov's distance (\cite{Kon22}), and the geometric quantile map and distribution function fully characterize the underlying probability measure (\cite*{Kol1997, Kon2025}). The multivariate concept was then extended to infinite-dimensional vector spaces (\cite*{ChaCha2014, Romon2022, PasReid2022, KonPai25}) and to manifolds (\cite*{KonPai23, Vir26}). Geometric quantiles have been criticized for their lack of affine-equivariance and for not controlling the probability content of their quantile regions; in our opinion, these issues are minor and can be dealt with through a transformation-retransformation procedure as in \cite*{Ser2010} and a relabeling of the quantile regions as in \cite{Kon2025}, respectively. \smallskip

Despite the attractive characteristics listed above, the main drawback of the geometric notion of multivariate quantiles and distribution function is two-fold. First, it was established in \cite{GirStu2017} that the shape of extreme geometric quantile contours is both unexpected and uninformative from an extreme value perspective. Second, a thick line is drawn between the properties of geometric objects in the discrete and continuous cases. Indeed, the continuity and invertibility result alluded to above and the description of extremes require that $P$ be non-atomic; besides, uniqueness of geometric quantiles and non-existence of extreme quantiles associated with an order $\alpha=1$ are always granted when $d\ge 2$ unless $P$ is supported on a line of $\R^d$ (in this case, geometric quantiles may be non-unique and extremes may exist). This yields a surprising `phase transition' from continuous to discrete, and fully supported to degenerate distributions; in particular, many population properties of geometric quantiles do not actually hold for empirical distributions. The fundamental reason behind these pathologies lies in the fact that the geometric distribution function $F_P$ from (\ref{eq:DefinGeometricCdf}) is obtained through integration of the \emph{singular} kernel $x\mapsto x/\|x\|$ so that, in particular, $F_P$ fails to be continuous at atoms of $P$. In addition to the theoretical obstacles this creates (non-differentiability, need for strong moment conditions) to derive classical M-estimation results, this singularity yields numerical instability when computing geometric quantiles through minimization of the objective (\ref{eq:DefinMalphaU}). 

\subsection{Contribution of this paper} \label{sec:Introduction:Contribution}

Having observed these two drawbacks of geometric quantiles, we thus introduce a regularized version of $F_P$ above by letting, for some $\rr:[0,\infty)\to [-1,1]$,
\begin{equation}\label{defFpr}
\Fpr(x) 
=
\int_{\R^d} \rr(\|x-z\|) \frac{x-z}{\|x-z\|}\I[z\neq x]\, dP(z) 
,\spa 
\forall\ x\in\R^d 
,
\end{equation}
so that $\Fpr$ takes values in the closed unit ball $\overline{\B^d}$. This `regularization' will be particularly relevant if $\rr$ is continuous at $0$ and satisfies $\rr(0)=0$, but we will study in Section \ref{sec:ExistUniq} the properties of the resulting regularized geometric objects for a generic $\rr$, with the usual geometric setting obtained for $\rr\equiv 1$. Our regularized geometric quantiles will then be obtained as global minimizers of the regularized objective function 
\begin{equation}\label{defMpr}
M_{\alpha,u}^{\rr,P}(x) 
\equiv 
\int_{\R^d} \big( \rR(x-z) - \rR(z)\big)\, dP(z) - \ps{\alpha u}{x}
,\spa 
\forall\ x\in\R^d 
,
\end{equation}
where 
\begin{equation}\label{eq:rR}
\rR(x) 
=
\int_0^{\|x\|} \rr(s)\, ds 
.
\end{equation}
When $M_{\alpha,u}^{\rr,P}$ is differentiable, we thus expect $\nabla M_{\alpha,u}^{\rr,P} = \Fpr - \alpha u$ so that, as for geometric quantiles, $\Fpr$ will be the inverse of the regularized quantile map. Observe that $\Fpr$ in (\ref{defFpr}) is linear in $P$ in the sense that for all Borel probability measures $P_1$ and $P_2$, we have
$$
F_{(1-s)P_1 + sP_2}^\rr = (1-s)F_{P_1}^\rr + s F_{P_2}^\rr
,\spa 
\forall\ s\in [0,1],
$$
as is the case for the classical univariate distribution functional $P\mapsto F_P(\cdot) = P[(-\infty, \cdot]]$. This linearity property plays a key role in the robustness properties of geometric quantiles established in \cite{KonPai2025a} (see Section 4 there), so that we expect our $\rr$-quantiles to share the robustness of geometric quantiles. We could have introduced a more general regularization effect by replacing $\rr(\|x-z\|)$ in~(\ref{defFpr}) by $\rho(x,z)$ for some $\rho:\R^d\times\R^d \to [-1,1]$. The choice $\rho(x,z)\equiv \rr(\|x-z\|)$, however, ensures that the resulting $\rr$-distribution function $\Fpr$ is equivariant under translations and orthogonal transformations of $P$ (see Section~\ref{sec:symmetries} for details). Furthermore, we establish in Section \ref{sec:Univer} the following universality result: any linear assignment $P\mapsto \Fd_P$ such that $\Fd_P$ takes its values in $\overline{\B^d}$ and is equivariant under translations and orthogonal transformations is necessarily of the form~(\ref{defFpr}) for some $\rr:[0,\infty)\to [-1,1]$.
% Regularize, universality, define $\rr$-quantiles, main theor, continuity in $\rr$
% Grosse motivation pour la linéarité : key to robustness, symmetry principles, influence function?\\
To ensure uniqueness and fast computability of our $\rr$-quantiles, we will restrict to regularizers $\rr$ that make the loss $\rR$ in (\ref{defMpr}) convex. For this purpose, we will consider the class of regularizers
% Restriction to nonnegative and nondecreasing $\rr$: natural (the further, the higher the cost), uniqueness of quantiles, very fast computation of the minimization problem and convergence guarantees, use of classical results on convex M-estimators for asymptotic normality. 
% $\rr$ need not be defined at $0$ but since we restrict ourselves to those $\rr$ that are non-negative, non-decreasing, and continuous, it is not restrictive to assume that $\rr$ is well-defined at $0$.
  % It will be important that $\rr(s)>0$ for all $s>0$.
\begin{equation}\label{defrR}
\RRR 
=
\Big\{ \rr:[0,\infty) \to [0,1] : \rr\in C^1,\ \rr(s) > 0 \text{ and } \rr'(s) \ge 0 \text{ for } s > 0,\  \lim_{s\to \infty} \rr(s) = 1
\Big\} 
.
\end{equation}
Combining results from Section \ref{sec:ExistUniq}, Section \ref{sec:QuantileMap}, and Section \ref{sec:Empiric}, we summarize the main properties of our $\rr$-quantiles in what follows. In the rest of the paper, we denote by $\PPP(\R^d)$ the collection of Borel probability measures over $\R^d$.

\begin{Theor}\label{TheorSummary}
Let $\rr\in\RRR$ and fix an arbitrary $P\in\PPP(\R^d)$. Assume that $\rr(0)=0$ and $\rr'(0)=0$, and that $\rr$ is Lipschitz and (strictly) increasing over $[0,\infty)$. Fix $\alpha\in [0,1)$ and $u\in\S^{d-1}$.
\begin{enumerate}[label=(\roman*)]
    \item The objective $M_{\alpha,u}^{\rr,P}$ from (\ref{defMpr}) is %of class $C^2(\R^d)$, 
    twice continuously differentiable and strictly convex over $\R^d$, and admits a unique global minimizer $\Qpr(\alpha u)$; we call it the $\rr$-quantile of order $\alpha$ in direction $u$ for $P$.

    \item The map $M_{1,u}^{\rr,P}$ admits no global minimum, \ie $P$ admits no extreme quantile associated with order $\alpha=1$ in direction $u$.

    \item The $\rr$-quantile map $\Qpr : \B^d \to \R^d$ is a homeomorphism between the open unit ball $\B^d$ and $\R^d$, with inverse $\Fpr:\R^d\to \B^d$ defined as in (\ref{defFpr}). 

    \item For any sequences $(\alpha_k)\subset [0,1)$ and $(u_k)\subset \S^{d-1}$ such that $\alpha_k \to 1$ and $u_k\to u\in\S^{d-1}$, we have $\|\Qpr(\alpha_k u_k)\|\to \infty$, and $\Qpr(\alpha_k u_k)/\|\Qpr(\alpha_k u_k)\| \to u$ as $k\to\infty$.

    \item Let $Z_1, Z_2, \ldots$ be a sample drawn from $P$ and let $P_n$ and $\hat q_n=Q_{P_n}^\rr(\alpha u)$ be the corresponding empirical measure and $\rr$-quantile. Then, $\sqrt{n}(\hat q_n - \Qpr(\alpha u)) \to \NN(0, \Sigma)$ in distribution as $n\to\infty$, for some non-degenerate covariance matrix $\Sigma=\Sigma(\alpha, u, \rr, P)$.

    \item For any compact subset $I\subset [0,1)$, $\tilde \rr\in\RRR$, $(\alpha, u)\in I\times \S^{d-1}$, and (potentially non-unique) $\tilde\rr$-quantile $\tilde q$ of order $\alpha$ in direction $u$ for $P$, there exists a constant $C=C(\rr, I, P, \|\tilde q\|)>0$ such that
    $$
    \|Q_P^{\rr}(\alpha u) - \tilde q\|^2
    \le 
    C \int_0^\infty |\rr(s) - \tilde \rr(s)|\, ds 
    .
    $$
    
    % Fix a compact subset $I\subset [0,1)$. Then, there exists a constant $C=C(I)>0$ such that 
    % $$
    % \sup_{\tilde\rr\in\RRR}
    % \sup_{P\in\PPP(\R^d)}
    % \sup_{\substack{\alpha \in I\\ u\in\S^{d-1}}} 
    % \|Q_P^{\rr}(\alpha u) - Q_P^{\tilde\rr}(\alpha u)\|^2
    % \le 
    % C \int_0^\infty \big( \rr(s) - \tilde \rr(s)\big)\, ds 
    % .
    % $$
\end{enumerate}  
\end{Theor}

This result is quite striking as it establishes that $\rr$-quantiles share the same qualitative properties as geometric quantiles, but now holds for \emph{all} probability measures, including empirical distributions or degenerate distributions that are concentrated on a single line. We further investigate in Section \ref{sec:QuantileMap} the properties of $\rr$-distribution functions and quantiles when $\rr$ %is not necessarily smoothing 
does not necessarily have a smoothing effect at $0$, \ie $\rr(0)>0$, thus including the usual geometric case: we provide new insights on the mapping properties of the quantile map associated with discrete measures and, in particular, on what we refer to as the `black hole phenomenon' around atoms, and we show in Section~\ref{sec:extremes} that extreme $\rr$-quantiles share the same qualitative behavior as their geometric counterpart and, hence, are also uninformative about the tail characteristics of the underlying distribution. Theorem \ref{TheorSummary}(vi) with $\tilde \rr\equiv 1$ ensures (quantitatively) that $\rr$-quantiles are close to their geometric antecedents if $\rr$ increases sufficiently fast to $1$. In addition, it can be shown that the constant $C$ in Theorem~\ref{TheorSummary}(vi) is uniform in $P$ ranging in any class $\FF\subset \PPP(\R^d)$ satisfying %some 
certain properties (see Section \ref{sec:ProofsSummary} for details); moreover, $\FF$ can be chosen as the collection of
% which there exists a bounded Borel set $E\subset\R^d$ such that $\inf_{P\in\FF} P[E]>0$. 
% Consequently, taking $\FF$ consisting of 
discrete distributions having at most $\lceil n(1-\alpha)/2\rceil$ atoms that differ from those of $P_n$. This implies that the breakdown point of an empirical $\rr$-quantile $\hat q_n = Q_{P_n}^\rr(\alpha u)$ as in Theorem~\ref{TheorSummary}(v) is at least as large as that of its geometric counterpart $\lceil n(1-\alpha)/2 \rceil / n$ (see Corollary~2.2 in \cite{KonPai2025a}). Our $\rr$-quantiles thus inherit the attractive robustness properties of traditional geometrical quantiles. Section~\ref{sec:Empiric} focuses on the particular case $P=P_n$ of empirical quantiles.

\section{Universality of $\rr$-distribution functions}\label{sec:Univer}

Geometric quantiles are equivariant with respect to translations and orthogonal transformations, \ie letting $T_a(z)=z+a$, for all $z\in\R^d$ and $a\in\R^d$, and  $T_U(z)=Uz$ for all $z\in\R^d$ and $d\times d$ orthogonal matrix $U$, we have
\begin{equation}\label{eq:EquivQuantile}
Q_{T_a\#P}=T_a\circ Q_P
,\spa 
\textrm{and}
\spa 
Q_{T_U\# P}\circ T_U = T_U\circ Q_P
,
\end{equation}
for all probability measures $P$ not supported on a single of $\R^d$, where $T\#P$ denotes the push-forward probability measure $(T\# P)[A]\equiv P[T^{-1}(A)]$ for all Borel sets $A\subset\R^d$ and corresponds to the law of the random vector $T(Z)$ when $Z$ has law $P$. 
% Thus letting $Q_X$ stand for $Q_P$ when $X\sim P$, the equivariance relation (\ref{eq:EquivQuantile}) reads $Q_{T(X)}(\alpha u)=T(Q_X(\alpha u))$ for all $\alpha\in [0,1)$ and $u\in\S^{d-1}$. 
Since we focus on distribution functions rather than quantile maps in this section, we express (\ref{eq:EquivQuantile}) in terms of distribution functions by formally inverting the equalities. This leads to the following equivariance relations: 
\begin{equation}\label{eq:EquivCdf}
F_{T_a\# P}
=
F_P\circ T_a^{-1}
,\spa 
\textrm{and}
\spa 
T_U^{-1}\circ F_{T_U\# P} = F_P\circ T_U^{-1}
.
\end{equation}
% Although the transformations $T$ we consider here are invertible, we prefer the formulation (\ref{eq:EquivCdf}) to $F_{T\# P}=F_P\circ T^{-1}$ as (\ref{eq:EquivCdf}) does not require in principle $T$ to be invertible. 
It is noteworthy that these equivariance relations also hold for measure-transportation quantiles and distribution functions and are, in fact, shared by essentially all depth concepts available in the literature.
% The fact the same equivariance relations hold with respect to translations and orthogonal transformations for measure-transportation quantile maps and distributions functions legitimates the definitions we adopt in Section \ref{sec:Universality}.
% \subsection{Universality}\label{sec:Universality}

We now focus on generic ways of associating a map to any probability measure (to be thought of as the corresponding distribution function). As for geometric and measure-transportation distribution functions, we wish to preserve the directional nature of the distribution function. We will thus restrict to maps defined on $\R^d$ and taking their values in the closed unit ball $\overline{\B^d}$: to any $P\in\PPP(\R^d)$, or any $P$ in a suitable sub-collection $\PP$ of $\PPP(\R^d)$, one associates a map $\Fd_P:\R^d\mapsto\B^d$. This gives rise to a map 
$$
\Fd : \PP\subset \PPP(\R^d)\to {\rm Maps}(\R^d,\overline{\B^d})
,\ 
P\mapsto \Fd(P)\equiv  \Fd_P
.
$$
We call any such $\Fd$ a \emph{distribution function mechanism}. If $\Fd_P$ is invertible we naturally call the inverse $(\Fd_P)^{-1}$ the quantile map of $P$. Although we will mostly be interested in the case $\PP=\PPP(\R^d)$, the main result of this section also applies when $\PP$ is a strict subset of $\PPP(\R^d)$. When $\PP \subsetneq \PPP(\R^d)$, we will make the following assumptions on $\PP$. %When $A$ is an arbitrary subset of a vector space, 
For an arbitrary subset $A$ of a vector space, we let ${\rm span}(A)$ denote the linear space spanned by $A$.

% we denote by ${\rm cone}(A)$ the convex cone spanned by $A$, \ie 
% $$
% {\rm cone}(A)
% =
% \Big\{ \sum_{i=1}^N \lambda_i a_i : a_i\in A,\ \lambda_i\geq 0,\ N\in\N\Big\} 
% .
% $$

\begin{Asmp}\label{Asmp:ClassP}
% For a convex collection of probability densities $\GG$, we have $\PP=\{dP = g\, dx : g\in\GG\}$. In addition, we assume that $\GG$ is a convex set, ${\rm span}(\GG)$ is dense in $L^1(\R^d)$, and for all $h\in {\rm span}(\GG)$ we have $|h|\in {\rm span}(\GG)$.
The class $\PP$ satisfies $\PP=\{dP = g\, dx : g\in\GG\}$, where $\GG\subset L^1(\R^d)$ is a convex collection of probability densities such that ${\rm span}(\GG)$ is dense in $L^1(\R^d)$ and $|h|\in {\rm span}(\GG)$ whenever $h\in {\rm span}(\GG)$.
    % , denoted by $\PPP(\R^d)\cap L^1(\R^d)$, such that ${\rm cone}(\GG)$ is dense in $\PPP(\R^d)\cap L^1(\R^d)$ with respect to the $L^1(\R^d)$-norm. Further assume that $|f-g|\in{\rm cone}(\GG)$ for all $f,g\in {\rm cone}(\GG)$. 
    % that for any continuous and compactly supported density $f$, there exists a sequence $(f_k)\subset {\rm span}(\GG)$ such that $f\leq f_k$ for all $k$ and $(f_k)$ converges to $f$ in $L^1(\R^d)$ as $k\to\infty$.
    % such that ${\rm span}(\GG)$ is dense in $L^1(\R^d)$,
    % % (ii) there exists a non-negative sequence $(p_k)\subset {\rm span}(\GG)$ such that $(p_k)\to 0$ in $L^1(\R)$
    % and for any non-negative $f\in L^1(\R^d)$ there exists a sequence $(f_k)\subset {\rm span}(\GG)$ such that $f_k\geq f$ for all $k$ and $(f_k)\to f$ in $L^1(\R^d)$.
    % % and for any non-negative $f\in L^1(\R)$ there exists $\overline{f}\in {\rm span}(\GG)$ such that $\overline{f}\geq f$. 
    % We then consider the corresponding class $\PP=\{dP = g\, dx : g\in\GG\}$ of probability measures.
\end{Asmp}

% In addition, for any non-negative map $f\in L^1(\R^d)$, there exists a sequence $(f_k)\subset\GG$ such that $f_k\geq f$ for all $k$ and $(f_k)\to f$ in $L^1(\R^d)$.

% Assumption \ref{Asmp:ClassP} immediately implies that ${\rm span}(\GG)$ is dense in $L^1(\R^d)$, but further requires that a one-sided approximation (from above) can be achieved. 
% It will be clear from the main result of this section (see Theorem \ref{TheorUniqueCdf} below) that Assumption \ref{Asmp:ClassP}, which is only used in Step 2 of the proof, can be relaxed into the following weaker condition to obtain the conclusion of Step 2: for any $h\in {\rm span}(\GG)$ there exists a sequence $(h_k)\subset {\rm span}(\GG)$ such that $h^-\leq h_k$ almost everywhere for all $k$ and $(h_k)\to h^-$ in $L^1(\R^d)$ as $k\to\infty$. 
Canonical examples of such classes $\GG$ are: continuous (or smooth) densities, bounded densities, or any restriction of the previous classes to densities $f$ for which there exists $\lambda=\lambda(f)>0$ such that $\lambda\leq f\leq \lambda^{-1}$. Restriction of these classes to compactly supported densities also gives rise to a collection $\GG$ that satisfies Assumption \ref{Asmp:ClassP} (or the weak version mentioned above).\smallskip

While many concepts of distribution functions can be defined in principle, not any mechanism yields a valid concept %For instance, in analogy with the univariate case, one could require that $\Fd_P$ be a monotone map (in a suitable sense) for any $P\in\PP$, leading to the measure transportation notion of distribution function from \cite{HallEspagne21}. 
unless further properties are required. 
In the literature, equivariance with respect to all translations, orthogonal or affine transformations are commonly required. Other properties, such as monotonicity along rays from the center of the distribution or maximality of the corresponding depth function at the center of symmetry have also been considered in the literature. Such requirements can be debated, but the quantile map should at the very least be equivariant with respect to translations and orthogonal transformations. Nonetheless, these mild equivariance properties are already enough to enforce strong symmetries on a %cdf 
distribution function mechanism. In line with (\ref{eq:EquivCdf}), we will thus assume the following equivariance relations on $\Fd$. We start with equivariance with respect to translations.

\begin{Asmp}\label{Asmp:EquivTrans}
    The mechanism $\Fd$ is translation-equivariant: for any $a\in\R^d$ and induced map $T_a(z)=z+a$ for all $z\in\R^d$, the equality $\Fd_{T_a\#P} = \Fd_P\circ T_{-a}$ holds for all $P\in\PP$.
\end{Asmp}

When $P$ admits a density $f_P$ with respect to the Lebesgue measure, we will abusively write $\Fd(f_P)$ instead of $\Fd_P=\Fd(P)$. Consequently, when all probability measures in $\PP$ have a density, it will be convenient to consider $\Fd$ as acting on the corresponding space of densities rather than on $\PP$ itself. The action of any diffeomorphism $\varphi:\R^d\to\R^d$ on $P$, through $P\mapsto \varphi\# P$, corresponds to  transforming the density of $P$ through $f_P\mapsto f_{\varphi\# P}$ with
$$
f_{\varphi\# P}(x) 
\equiv 
\frac{f_P(\varphi^{-1}(x))}{|\det(J_x \varphi)|}
,\spa 
\forall\ x\in\R^d
,
$$
where $J_x\varphi$ stands for the Jacobian matrix of $\varphi$ at $x$. In particular, we have $f_{T_a\# P} = f_P \circ T_{-a}$ for all $x\in\R^d$ and $a\in\R^d$. Assumption \ref{Asmp:EquivTrans} %then rewrites 
can then be rewritten as $\Fd(f_P\circ T_{-a}) = (\Fd(f_P))\circ T_{-a}$ for all $a\in\R^d$ and any $P\in\PP$ with density $f_P$, or, equivalently,
\begin{equation}\label{eq:EquivCdfDensityTrans}
\Fd(f_P\circ T_a) = \Fd(f_P)\circ T_a
,\spa 
\forall\ a\in\R^d 
.
\end{equation}
We now turn to equivariance with respect to orthogonal transformations.

\begin{Asmp}\label{Asmp:EquivOrtho}
    The mechanism $\Fd$ is orthogonal-equivariant: for any $d\times d$ orthogonal matrix~$U$ and induced map $T_U(z)=Uz$ for all $z\in\R^d$, the equality $T_U^{-1}\circ \Fd_{T_U\#P} = \Fd_P\circ T_{U}^{-1}$ holds for all $P\in\PP$.
\end{Asmp}

A similar reasoning as for Assumption \ref{Asmp:EquivTrans} entails that Assumption \ref{Asmp:EquivOrtho} %rewrites, 
can be rewritten, for any $P\in\PP$ with density $f_P$, as
\begin{equation}\label{eq:EquivCdfDensityOrtho}
T_U\circ \Fd(f_P\circ T_U) = \Fd(f_P)\circ T_U
,\spa 
\forall\ U~ \text{orthogonal}
.
\end{equation}
Finally, observe that the univariate distribution functional $P\mapsto P[(-\infty, \cdot]]$ is linear in $P\in\PPP(\R^d)$. The same holds for the geometric distribution function $P\mapsto \int_{\R^d\sm\{\cdot\}} (\cdot-z)/\|\cdot-z\|\, dP(z)$. This property is at the core of the robustness properties of geometric quantiles, as it gives a meaningful and tractable interpretation of the geometric distribution function of a contaminated version $(1-s)P + sQ$ of $P$; see Section 4 in \cite{KonPai2025a}.
% In what follows, we will further restrict to linear mechanisms, \ie the assignement $P\mapsto \Fd_P$ is linear in $P\in\PPP(\R^d)$. Although it may seem restrictive, we argue that this assumption is most natural: in dimension $d=1$ the usual distribution function, that assigns any $P\in\PPP(\R)$ to the map $z\mapsto P[(-\infty,z])$, is indeed linear in $P$, and the same holds for the geometric distribution function as can be seen from (\ref{eq:DefinGeometricCdf}). 
The terminology `linear' is slightly abusive since the collection $\PPP(\R^d)$ of Borel probability measures is not a vector space. However, in the two previous examples (the usual univariate and the multivariate geometric distribution functions), the distribution function mechanism is in fact well-defined for any finite measure and, subsequently, for all (finite) signed measures, which form a vector space on which the previous mechanisms are indeed linear. Nonetheless, because we do not want to assume \emph{a priori} that a mechanism is defined for all (finite) signed measures, but only for a subcollection $\PP\subset\PPP(\R^d)$, we adopt the following notion of `linearity' that rather exploits the convex nature of $\PPP(\R^d)$.

\begin{Asmp}\label{Asmp:Linear}
    The collection $\PP\subset\PPP(\R^d)$ is convex and the distribution function mechanism $\Fd:\PP\to {\rm Maps}(\R^d,\overline{\B^d})$ is linear in the sense that, for all $P_1,P_2\in\PP$, we have
    \begin{equation}\label{eq:LinearConvex}
	\Fd\big((1-s)P_1 + sP_2\big) = (1-s)\Fd(P_1) + s \Fd(P_2)
    ,\spa 
    \forall\ s\in [0,1]
    .
	\end{equation}
\end{Asmp}

% Interpretation in terms of the influence function + limit argument ?\\

%Recalling the definition of $F_P^\rr$ from~(\ref{DefinFpr}), we 
We can now state the main result of this section.

\begin{Theor}\label{TheorUniqueCdf}
     Let $\PP\subset \PPP(\R^d)$ and $\Fd:\PP\to {\rm Maps}(\R^d, \overline{\B^d})$ satisfying Assumptions \ref{Asmp:ClassP}, \ref{Asmp:EquivTrans}, \ref{Asmp:EquivOrtho}, and \ref{Asmp:Linear}. Then $\Fd$ extends linearly to $\PPP(\R^d)$ and there exists a bounded and measurable map~$\rr:(0,\infty)\to [-1,1]$ such that
		\[
		(\Fd_P)(x)
        =
		F_P^\rr(x)
        \equiv
            \int_{\R^d} \rr(\|x-z\|) \frac{x-z}{\|x-z\|} \I[z\neq x]\, dP(z),
		\]
		for all $x\in\R^d$ and $P\in\PPP(\R^d)$.
\end{Theor}

This motivates the following definition which will feature prominently throughout.

\begin{Defin}\label{DefinFpr}
Let $P\in\PPP(\R^d)$ and $\rr\in\RRR$. The \emph{$\rr$-geometric distribution function} of $P$ is the map $F_P^\rr:\R^d\to\R^d$ defined by 
$$
\Fpr(x)
\equiv 
% \E\big[ (\nabla \rR)(x-Z)\I[Z\neq x]\big]
% =
\E\bigg[\rr(\|x-Z\|)\frac{x-Z}{\|x-Z\|}\I[Z\neq x]\bigg]
% \int_{\R^d} \rr(\|x-z\|)\frac{x-z}{\|x-z\|}\I[z\neq x]\, dP(z)
,\spa 
\forall\ x\in\R^d 
,
$$
where $Z$ is a random $d$-vector with distribution $P$.
\end{Defin}

We recover the usual geometric distribution function for $\rr\equiv 1$, which is equivariant under translations and orthogonal transformations. The more general form in Definition \ref{DefinFpr} further exhausts \emph{all} possible such equivariant distribution functions that are also linear functionals of $P$, by virtue of Theorem \ref{TheorUniqueCdf}. In the following sections, we will study how the choice of $\rr$ impacts the properties of the resulting $\rr$-concept of distribution function and quantiles.

\section{Existence and uniqueness}\label{sec:ExistUniq}

In this section, we establish existence of our regularized $\rr$-geometric quantiles as well as conditions under which they are unique. Recall the definition of the class $\RRR$ of regularizers from~(\ref{defrR}), and the loss $\rR$ from (\ref{eq:rR}).

% Following remarks in Section~\ref{sec:Introduction}, we introduce the class of regularizers that will be relevant for our purposes. Let
% $$
% \RRR 
% =
% \bigg\{ \rr:[0,\infty) \to [0,1] : \rr\in C^1,\ \rr(s) > 0 \text{ and } \rr'(s) \ge 0 \text{ for } s > 0,\  \lim_{s\to \infty} \rr(s) = 1
% \bigg\} 
% .
% $$
% For $\rr\in\RRR$ fixed, define
% \begin{equation}\label{eq:rR}
% \rR(x)
% =
% \int_{0}^{\|x\|} \rr(s)\, ds
% .
% \end{equation}

\begin{Defin}\label{DefinQuantiles}
    Let $P\in\PPP(\R^d)$ and $\rr\in\RRR$ with corresponding $\rR$ defined as in (\ref{eq:rR}). Fix $\alpha\in [0,1]$ and $u\in\S^{d-1}$. We say that $q_{\alpha,u}^\rr=q_{\alpha,u}^\rr(P)$ is a $\rr$-geometric quantile of order $\alpha$ in direction $u$ for $P$ if and only if it is a global minimizer over $\R^d$ of the objective function 
    \begin{equation}\label{defM}
    x\mapsto 
    M_{\alpha,u}^{\rr,P}(x)
    \equiv 
    \int_{\R^d} \big(\rR(z-x)-\rR(z)\big)\, dP(z)
    -
    \ps{\alpha u}{x}
    .
    \end{equation}
\end{Defin}
The choice $\rr\equiv 1$ in the definition above provides $\rR(x)=\|x\|$ so that the resulting $\rr$-quantiles coincide with the usual geometric quantiles from \cite{Cha1996}. The term $\rR(z)$ in the integral defining $M_{\alpha,u}^{\rr, P}$ may look superfluous but allows avoiding moment conditions on $P$ by a standard triangle inequality-type argument; in particular, the existence result of $\rr$-geometric quantiles in Theorem \ref{TheorExist} below holds for an arbitrary probability measure.\smallskip

Since $\rr\in C^1$, then $\rR$ is of class $C^2$ on $\R^d\sm\{0\}$ with 
\begin{equation*}\label{eq:JacR}
(\nabla \rR)(x)
=
\rr(\|x\|)\frac{x}{\|x\|}
,
\end{equation*}
and, for $\nabla^2 \rR(x)$ the Hessian matrix of $\rR$ at $x$,
\begin{equation}\label{eq:HessR}
(\nabla^2 \rR)(x)
=
\rr'(\|x\|)\frac{xx^T}{\|x\|^2} + \frac{\rr(\|x\|)}{\|x\|}\Big({\rm I}_d - \frac{xx^T}{\|x\|^2}\Big)
.
\end{equation}
Consequently, the fact that $\rr\ge 0$ and $\rr'\ge 0$ entails that $\rR$ is convex over $\R^d$. In particular, the objective function $M_{\alpha,u}^{\rr,P}$ to be minimized in (\ref{defM}) is convex. We start with the following general existence result.

\begin{Theor}\label{TheorExist}
    Fix $P\in\PPP(\R^d)$ and $\rr\in\RRR$. Let $\alpha\in [0,1]$ and $u\in\S^{d-1}$. Then, 
    \begin{enumerate}[label=(\roman*)]
    \item The map $x\mapsto M_{\alpha,u}^{\rr,P}(x)$ is continuous over $\R^d$.
    \item If $\alpha<1$, then $P$ admits at least one $\rr$-geometric quantile $q_{\alpha,u}^\rr$ of order~$\alpha$ in direction $u$.
    \end{enumerate}
\end{Theor}

Observe that, as for the classical geometric quantiles, this existence result holds for an abitrary probability measure without imposing moment conditions. The case $\alpha=1$, corresponding to extreme quantiles, is not covered by Theorem \ref{TheorExist}. Indeed, depending on the regularizer $\rr$, such quantiles might or might not exist. This is the content of the next proposition; see also Section \ref{sec:extremes} for quantitative versions of this result.
% his is further investigated in . 

% $$
% \rR_\infty 
% \equiv 
% \int_0^\infty \big(1-\rr(s)\big)\, ds
% .
% $$

\begin{Prop}
\label{PropQuantileOrderOne}
Fix $P\in \PPP(\R^d)$ and $\rr\in\RRR$. Let $u\in \S^{d-1}$. Then $P$ admits an (extreme) $\rr$-geometric quantile $q_{1,u}^\rr$ if and only if the following holds: (i) $P$ is supported on a halfline with direction $-u$, and (ii) there exists $\tau\geq 0$ such that $\rr(s)=1$ if and only if~$s\geq \tau$. In that case, letting $y_0\in\R^d$ be such that $\{y_0-\lambda u : \lambda\geq 0\}$ is the smallest halfline with direction~$-u$ (with respect to the inclusion) on which $P$ is supported, a point $x\in\R^d$ is such a quantile if and only if $x= y_0 + \lambda u$ for some $\lambda\geq \tau$.
%	, and (ii) the integral $\int_{\R^d} (\|z\| + u'z)\, dP(z)$ exists in $\R$.
\end{Prop}

\begin{landscape}

%%%%%%%%%%%%%%%%%%%%%%%%%%%%%%%%%%
\begin{figure}[h!]
\centering
\includegraphics[width=0.9\linewidth]{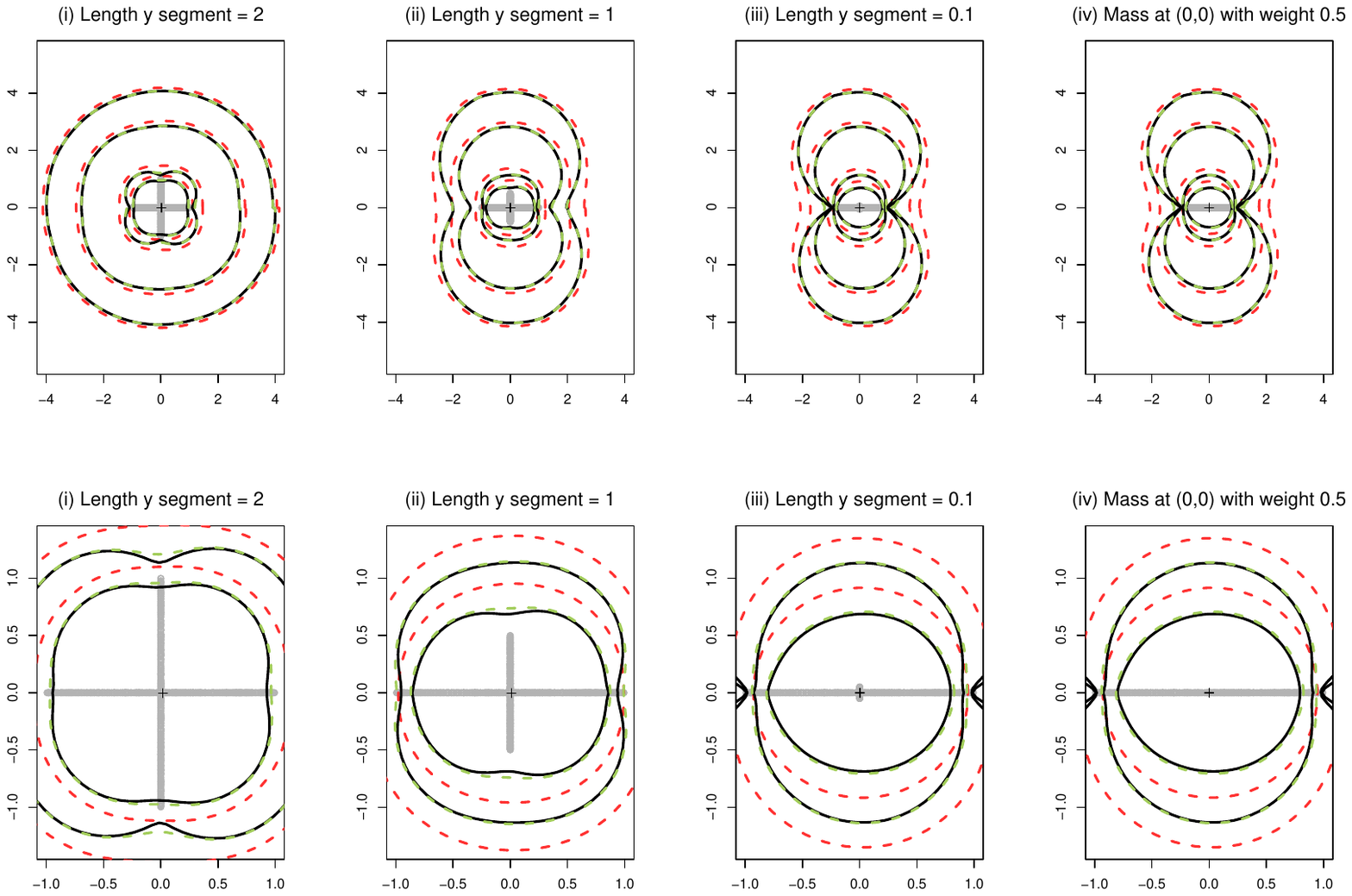}
\vspace{-1cm}
\caption{Graphical representation of the empirical quantiles $q_{\alpha,u}^\rr$, where $\rr(s)=\rr_{\beta}(s)=1-(1+s)^{-\beta}$, for $\beta=5$ (dashed red curves) and $\beta=10$ (dashed green curves), compared with the classical geometric quantiles $q_{\alpha,u}^{1}$ (full black curves), at levels $\alpha\in \{0.9,0.95,0.99,0.995\}$. The experiments use $n=1000$ data points generated from an equally weighted mixture of the uniform distribution on $[-1,1] \times \{ 0 \}$ and, from left to right, (i) the uniform distribution on $\{ 0 \} \times [-1,1]$, the uniform distribution on $\{ 0 \} \times [-1/2,1/2]$, the uniform distribution on $\{ 0 \} \times [-1/20,1/20]$, and a Dirac mass at the origin. The cross marks the geometric median of the data. The bottom row is a series of zoomed-in versions of the top panels on the square $[-1,1]\times [-1,1]$.}
\label{fig:extreme1to4}
\end{figure}
%%%%%%%%%%%%%%%%%%%%%%%%%%%%%%%%%%

\end{landscape}

\begin{landscape}

%%%%%%%%%%%%%%%%%%%%%%%%%%%%%%%%%%
\begin{figure}[h!]
\centering
\includegraphics[width=0.9\linewidth]{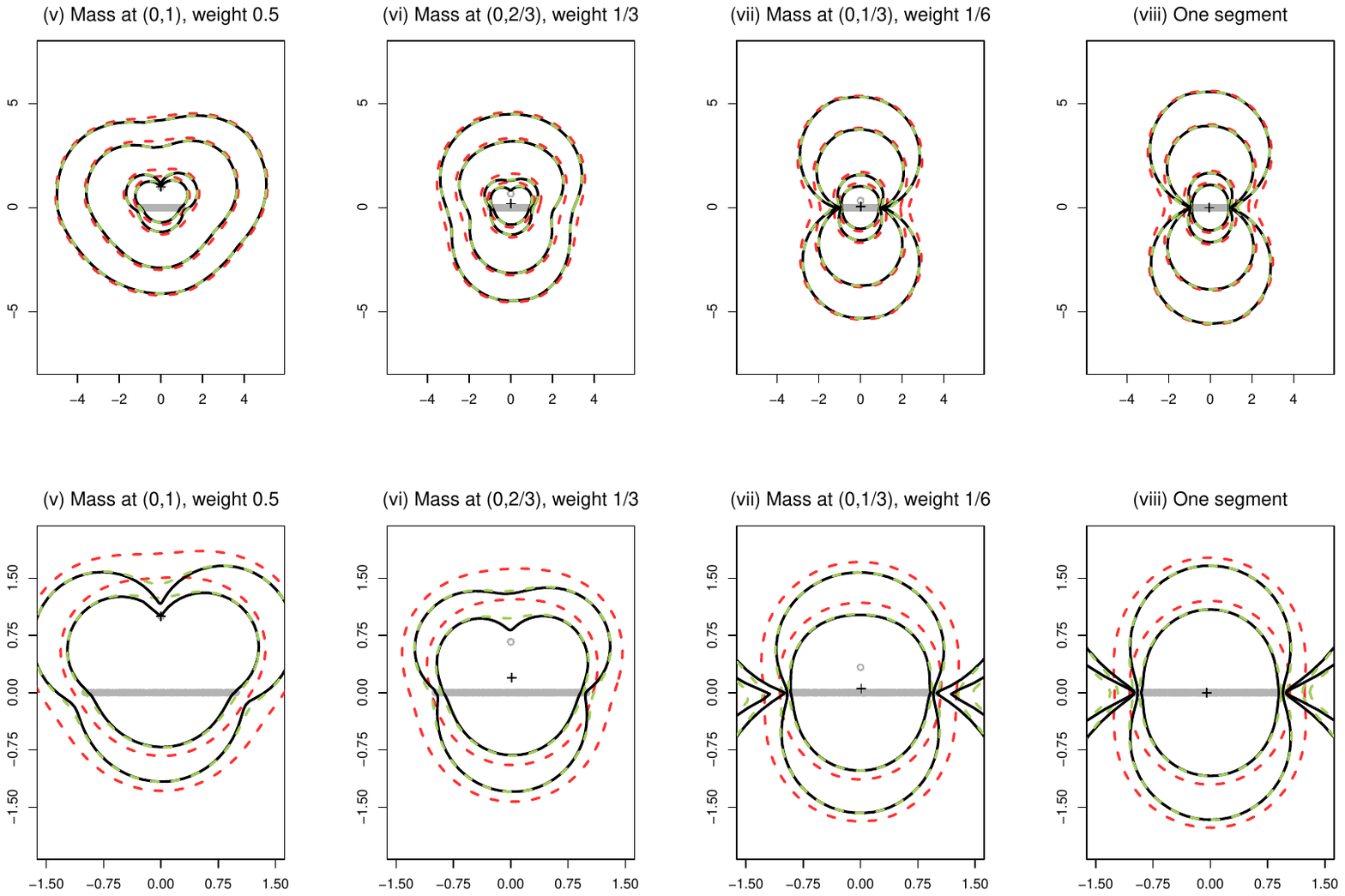}
\vspace{-1cm}
\caption{Graphical representation of the empirical quantiles $q_{\alpha,u}^\rr$, where $\rr(s)=\rr_{\beta}(s)=1-(1+s)^{-\beta}$, for $\beta=5$ (dashed red curves) and $\beta=10$ (dashed green curves), compared with the classical geometric quantiles $q_{\alpha,u}^{1}$ (full black curves), at levels $\alpha\in \{0.9,0.95,0.99,0.995\}$. The experiments use $n=1000$ data points generated from a weighted mixture of the uniform distribution on $[-1,1] \times \{ 0 \}$ and, from left to right, a Dirac mass at (i) $(0,1)$ with weight $1/2$, (ii) $(0,2/3)$ with weight $1/3$, (iii) $(0,1/3)$ with weight $1/6$. The rightmost panels use $n=1000$ data points uniformly generated on $[-1,1] \times \{ 0 \}$. The cross marks the geometric median of the data. The bottom row is a series of zoomed-in versions of the top panels on the square $[-1.5,1.5]\times [-1.5,1.5]$.}
\label{fig:extreme5to8}
\end{figure}
%%%%%%%%%%%%%%%%%%%%%%%%%%%%%%%%%%

\end{landscape}

It follows from Proposition \ref{PropQuantileOrderOne} that no probability measure over $\R^d$ admits extreme $\rr$-geometric quantiles when $\rr(s)<1$ for all $s>0$. When $\rr$ is eventually constant and equal to $1$, extreme $\rr$-quantiles in a given direction $u\in\S^{d-1}$ exist only when $P$ is supported on a halfline $\LL$ with direction $-u$, which then reduces to computing the usual extreme quantile of a univariate distribution. We illustrate this in Figures~\ref{fig:extreme1to4} and~\ref{fig:extreme5to8}. In both, we generated data from distributions that gradually converge to %the uniform 
a distribution concentrated
on a segment in $\R^2$. In each case, we compute $\rr$-quantiles for regularizers $\rr$ of the form
\[
\rr(s)=\rr_{\beta}(s)=1-\frac{1}{(1+s)^{\beta}}
\]
that will feature prominently throughout this article. Let us make a few remarks. (i) The larger~$\beta$, the faster $\rr_{\beta}(s)$ converges to 1 as $s\to\infty$, and the closer extreme $\rr$-geometric quantiles will be to the usual extreme geometric quantiles. (ii) When the distribution considered is not supported on a single half-line, then $\rr$-geometric quantiles at level $\alpha=1$ never exist, regardless of the choice of $\rr$. (iii) When the distribution is supported on a single half-line $\LL$ with direction $-u$, then $\rr_{\beta}$-geometric quantiles at level $\alpha=1$ still do not exist (since $\rr_\beta(s)<1$) although their geometric counterparts exist in the direction $u$ opposite to that of $\LL$; the same holds when $\rr$ is eventually constant equal to $1$, and the longer it takes $\rr$ to reach its eventually constant value, the further away the $\rr$-quantiles in direction $u$ will be from $\LL$.\smallskip

By definition, $\rr$-geometric quantiles are global minimizers of the convex objective function~$x\mapsto M_{\alpha, u}^{\rr,P}(x)$. When $M_{\alpha,u}^{\rr, P}$ is differentiable, these quantiles are thus characterized by the first-order condition~$\nabla M_{\alpha, u}^{\rr,P}(x)=0$; the gradient of $M_{\alpha,u}^{\rr,P}$ will further play a key role to establish the mapping properties of $\rr$-geometric quantiles in the Section \ref{sec:QuantileMap} below. By convexity, minimizers of $M_{\alpha,u}^{\rr, P}$ are in fact characterized by the weaker directional first-order conditions $(\partial M_{\alpha,u}^{\rr,P})/(\partial v)\ge 0$ for all $v\in\S^{d-1}$. While $M_{\alpha,u}^{\rr,P}$ may not always be globally differentiable---when $\rr\equiv1$, differentiability fails precisely at the atoms of $P$ (Theorem 5.1 in \cite{KonPai1} with $\rho(t)=t$)---the following result establishes existence of all directional derivatives of~$M_{\alpha,u}^{\rr,P}$, irrespective of $\rr\in\RRR$ and $P\in\PPP(\R^d)$.

\begin{Prop}\label{Propdirectionalderiv}
Let $P\in\PPP(\R^d)$ and $\rr\in\RRR$. Fix~$\alpha\in[0,1]$ and $u\in\S^{d-1}$. Let~$Z$ be a random $d$-vector with distribution~$P$. Then, for any~$x\in\R^d$ and~$v\in\R^d\setminus\{0\}$, the directional derivative  
$$
\frac{\partial M_{\alpha, u}^{\rr,P}}{\partial v}(x)
=
\lim_{t\stackrel{>}{\to} 0} \frac{M_{\alpha, u}^{\rr,P}(x+tv)-M_{\alpha, u}^{\rr,P}(x)}{t}
$$
exists and is given by
    $$
    \frac{M_{\alpha, u}^{\rr,P}}{\partial v}(x)
=	
    \rr(0) \|v\| P[\{x\}]  + \ps{v}{\E\bigg[\rr(\|x-Z\|)\frac{x-Z}{\|x-Z\|}\I[Z\neq x]\bigg]-\alpha u}
    .
    $$
\end{Prop} 

In view of the directional first-order conditions alluded to above, Proposition \ref{Propdirectionalderiv} entails that a point $x\in\R^d$ is a $\rr$-geometric quantile for $P$ of order $\alpha\in [0,1]$ in direction $u\in\S^{d-1}$  if and only if 
\begin{equation}\label{eq:DirCond}
\| 
% \E\bigg[\rr(\|x-Z\|)\frac{x-Z}{\|x-Z\|}\I[Z\neq x]\bigg]
\Fpr(x)
-
\alpha u 
\| 
\leq 
\rr(0)P[\{x\}]
.
\end{equation}
It further follows from Proposition~\ref{Propdirectionalderiv} that for $M_{\alpha,u}^{\rr,P}$ to be differentiable at~$x\in\R^d$, it is necessary that~$\rr(0)P[\{x\}]=0$. The next result shows that this condition is sufficient, and relates $\nabla M_{\alpha,u}^{\rr,P}$ to the $\rr$-distribution function $\Fpr$. 

% The following result characterizes precisely when the objective function $M_{\alpha, u}^{\rr,P}$. We start with the following result. 

\begin{Theor} 
	\label{TheorDiffM}
	Let $P\in\PPP(\R^d)$ and $\rr\in\RRR$. Fix~$\alpha\in[0,1]$ and $u\in\S^{d-1}$. Then, 
    \begin{enumerate}[label=(\roman*)]
    \item $M_{\alpha, u}^{\rr,P}$ is differentiable at $x\in\R^d$ if and only if $\rr(0)P[\{x\}]=0$, in which case the corresponding gradient is
	$$
	\nabla M_{\alpha, u}^{\rr,P}(x)
	=
	\Fpr(x)-\alpha u
	.
        $$
	\item If~$\rr(0)P[\{x\}]=0$ for all~$x$ in an open set~$U\subset\R^d$, then~$M_{\alpha, u}^{\rr,P}\in C^1(U)$.
    \end{enumerate}
\end{Theor}

It follows from Theorem \ref{TheorDiffM} that when $\rr(0)=0$, the objective function $M_{\alpha,u}^{\rr, P}$ is differentiable on all of $\R^d$. This is in sharp contrast with the case of classical geometric quantiles, obtained for $\rr\equiv 1$, for which the corresponding $M_{\alpha,u}^{\rr,P}$ fails to be differentiable at atoms of $P$, thus resulting in less desirable properties of empirical geometric quantiles such as a lack of invertibility and continuity of the quantile map in a neighborhood of atoms. 
% \com{In the following sections, we will be considering regularizers $\rr$ satisfying $\rr(0)=0$, such as 
% $$
% \rr_\beta(s)
% =
% \frac{s^\beta}{1+s^\beta}
% ,
% $$
% where $\beta>0$ controls the rate at which $\rr$ increases to $1$. [Ecrire la bonne fonction $\rr$ utilisée dans les simus.]} 
Using (\ref{eq:DirCond}), we are able to provide a rather precise characterization of uniqueness of $\rr$-geometric quantiles in the next result.

\begin{Theor}\label{TheorUniqueness}
    Let $P\in\PPP(\R^d)$ and $\rr\in\RRR$. Fix $\alpha\in [0,1)$ and $u\in\S^{d-1}$. 
    \begin{enumerate}[label=(\roman*)]
    \item If $\rr$ is strictly increasing over $[0,\infty)$, then $q_{\alpha,u}^\rr$ is unique. 
    \end{enumerate}
    Otherwise, the following holds: 
    \begin{enumerate}[label=(\roman*),resume]
    \item If $P$ is not supported on a single line, then $q_{\alpha,u}^r$ is unique.
    % $P$ admits a unique $\rr$-geometric quantile $q_{\alpha,u}^\rr$ of order $\alpha$ in direction~$u$. 
    \item If~$\alpha>0$ and~$P$ is supported on a line $\LL$ with direction $v\in \S^{d-1}\sm\{\pm u\}$, then $q_{\alpha,u}^{\rr}$ is unique. In addition, if $\rr(0)>0$ and $q_{\alpha,u}^{\rr}\in\LL$ then $q_{\alpha,u}^{\rr}$ is an atom of $P$, whereas $q_{\alpha,u}^{\rr}\in\R^d\sm\LL$ if $\rr(0)=0$. 
    \item If~$P$ is supported on a line $\LL=\{x_0+\lambda u : \lambda\in\R\}$ with direction $u$, then any $\rr$-geometric quantile~$q_{\alpha,u}^\rr$ belongs to $\LL$ and writes $q_{\alpha,u}^\rr=x_0+\lambda_\alpha u$, where $\lambda_\alpha$ is a $\rr$-geometric quantile of order $\alpha$ in direction $1$ for $P_{u,x_0}$, with $P_{u,x_0}$ the distribution of $\ps{u}{Z-x_0}$ when $Z$ has distribution $P$.
    \end{enumerate}
\end{Theor}

The proof of (i)-(ii) proceeds by showing that $M_{\alpha,u}^{\rr,P}$ is strictly convex on $\R^d$ under the stated assumptions. Consequently, a regularizer $\rr$ that is strictly increasing over $[0,\infty)$ and satisfies $\rr(0)=0$ always results in unique $\rr$-quantiles of any order $\alpha\in [0,1)$, non-existence of extreme $\rr$-quantiles of order $\alpha=1$, and a strictly convex differentiable objective function $M_{\alpha,u}^{\rr,P}$. In (iii), we show that the $\rr$-quantile in a direction different than that of the support of $P$ remains unique, whereas (iv) is a special case of a more general result stated in Proposition \ref{PropChaudInHyper} below in the setting of distributions supported on a linear subspace of $\R^d$.
% , that establishes that, when $P$ is supported on a hyperplane $H$, then for any $u$ in the unit sphere of $H$, \ie $u\in H-x_0$ with $\|u\|=1$, then $q_{\alpha,u}^\rr\in H$.

% A few comments are in order. In (ii), when $P$ is non-atomic, then all $\rr$-geometric quantiles $q_{\alpha,u}$ with $\alpha>0$ belong to $\R^d\sm\LL$.

% In (iii) when $\alpha=0$, then $u$ is irrelevant and can be taken equal to $v$ so that $q_{\alpha,u}^\rr$ satisfies the form given in (iii).

% (iv) is a special case of the following more general result, that establishes that, when $P$ is supported on a hyperplane $H$, then for any $u$ in the unit sphere of $H$, \ie $u\in H-x_0$ with $\|u\|=1$, then $q_{\alpha,u}^\rr\in H$.

\section{Symmetries}\label{sec:symmetries}

As we explained in Section \ref{sec:Introduction} and Section \ref{sec:Univer}, traditional geometric quantiles enjoy natural equivariance properties. In this section, we investigate similar features for our regularized geometric quantiles.\smallskip 

When $V$ is a $k$-dimensional linear subspace of $\R^d$ and $H$ is the affine space $H=V+x_0$, for some $x_0\in\R^d$, we define the unit sphere $\S(H)$ of $H$ as 
$$
\S(H)
=
\big\{u\in V : \|u\|=1\big\} 
.
$$
In the following result, we show that when a probability measure is supported on an affine space $H$, then $\rr$-quantiles in directions belonging to $\S(H)$ also belong to $H$. This extends %Theorem~2.7(iv) 
Theorem~\ref{TheorUniqueness}(iv) established in the setting where $H$ is a line.
% $V=(v_1\, \ldots\, v_k)$ is the $d\times k$ matrix whose image is ${\rm span}(v_1,\ldots, v_k)$, let $H=\{x_0 + Vy:y\in\R^k\}$ be a $k$-dimensional 

\begin{Prop}
	\label{PropChaudInHyper}	
	Let $V=\{Oy : y\in\R^k\}$ be a $k$-dimensional linear subspace of $\R^d$, for some $d\times k$ matrix $O$ with orthonormal columns. Fix $x_0\in\R^d$ and let $H$ be the affine space $H=V+x_0$. Let~$P\in\PPP(\R^d)$ be supported on $H$. Fix~$\alpha\in[0,1)$ and~$u\in\S(H)$. Then, any $\rr$-geometric quantile $q_{\alpha,u}^\rr$ for $P$ belongs to~$H$, and %writes 
  can be written $q_{\alpha,u}^\rr=x_0+O\xi_{\alpha,u}$, where~$\xi_{\alpha,u}$ is a $\rr$-geometric quantile of order~$\alpha$ in direction~$O^Tu$ for $P_{O,x_0}$, with $P_{O,x_0}$ the distribution of $O^T(Z-x_0)$ when $Z$ has distribution $P$.
\end{Prop}

% We conclude this section with the following orthogonal- and translation-equivariance result, that will be particularly relevant in Section~\ref{sec:Spher} when considering the particular case for which~$P$ is spherically symmetric (the proof readily follows from Definition~\ref{DefinMQuantile}).  

The next proposition is an extension of the equivariance of geometric quantiles under orthogonal transformations and translations to their regularized counterparts; and follows straight from Definition~\ref{DefinQuantiles} and the first-order condition~(\ref{eq:DirCond}). It will be particularly relevant in Proposition~\ref{TheorSpher} below when considering the particular case %for 
in which~$P$ is spherically symmetric.

\begin{Prop}
\label{PropEquivariance}
		Let $P\in\PPP(\R^d)$ and $\rr\in\RRR$. Fix~$\alpha\in[0,1)$ and~$u\in\S^{d-1}$. Let $b\in\R^d$, $O$ be a $d\times d$ orthogonal matrix, and denote as~$P_{O,b}$ the distribution of~$OZ+b$ when~$Z$ has distribution~$P$. If~$q_{\alpha,u}^\rr$ is a $\rr$-geometric quantile of order~$\alpha$ in direction~$u$ for $P$, then~$Oq_{\alpha,u}^\rr+b$ is a $\rr$-geometric quantile of order~$\alpha$ in direction~$Ou$ for~$P_{O,b}$. 
\end{Prop}

Traditional geometric quantiles fail to be equivariant under general affine transformations. Nonetheless, they can be made affine-equivariant through a transformation-retransformation approach; see, e.g., \cite{Ser2010}. We now focus on the special case for which~$P\in\PPP(\R^d)$ is spherically symmetric about some location~$\mu_0\in\R^d$, \ie $P[\mu_0+OB] = P[\mu_0 + B]$ for any Borel subset~$B\subset\R^d$ and~$d\times d$ orthogonal matrix~$O$. By virtue of the translation-equivariance established in Proposition \ref{PropEquivariance} above, we will restrict, without any loss of generality, to the case~$\mu_0=0$. It follows from Proposition~\ref{PropEquivariance} that, if~$P$ is spherically symmetric about the origin of~$\R^d$ and satisfies~$P[\{0\}]<1$, with~$d\geq 2$ say, then any $\rr$-quantile contour $\{q_{\alpha, u}^\rr : u\in\S^{d-1}\}$, with~$\alpha\in [0,1)$, is a hypersphere (uniqueness of these $\rr$-quantiles follows from Theorem~\ref{TheorUniqueness}(ii) since~$P$ is then not supported on a line). Proposition~\ref{PropEquivariance} further implies that, for any~$\alpha\in[0,1)$ and~$u\in\S^{d-1}$, the~$\rr$-quantile $q_{\alpha,u}^\rr$ is invariant under all orthogonal transformations fixing~$u$; in particular, $q_{\alpha,u}^\rr$ belongs to the line spanned by~$u$, which is most natural. This motivates the next result.

% If $P$ is spherically symmetric, $q_{\alpha,u}^\rr$ is unique there exists a non-decreasing map $h:[0,\infty)\to [0,1)$ such that $h(0)=0$ and
%         $$
%         \Fpr(x) 
%         =
%         h(\|x\|) \frac{x}{\|x\|}
%         ,\spa 
%         \forall\ x\in\R^d 
%         .
%         $$

\begin{Prop} 
	\label{TheorSpher}
	Fix $P\in\PPP(\R^d)$ and $\rr\in\RRR$. Assume that~$P$ is spherically symmetric about the origin of~$\R^d$ and that $p\equiv P[\{0\}]<1$. Then $q_{\alpha,u}^\rr$ is unique for every $\alpha\in [0,1)$ and there exists a continuous and non-decreasing map $\QQ=\QQ^\rr:[0,1)\to [0,\infty)$ such that $q_{\alpha,u}^\rr = \QQ(\alpha) u$ for all~$\alpha\in [0,1)$ and~$u\in\S^{d-1}$. In addition, $\QQ(\alpha)=0$ for all $\alpha\in [0,\rr(0)p]$ and $\QQ$ is (strictly) increasing over $(\rr(0)p,1)$.
\end{Prop}

% Notice that when $P[\{0\}]<1$ or $\rr(s)<1$ for all $s>0$, then $h$ does not take the value $1$ since $P$ is not supported on a line with 

Proposition \ref{TheorSpher} excludes the degenerate situation where $P$ consists of a single Dirac mass at the origin of $\R^d$. In this case, one can directly show that the origin is a $\rr$-geometric quantile of any order $\alpha\in [0,\rr(0)]$ in any direction $u\in\S^{d-1}$, and that any $x\in\R^d\sm\{0\}$ is a $\rr$-geometric quantile of order $\alpha=\rr(\|x\|)\in (0,1]$ in direction $x/\|x\|\in\S^{d-1}$.\smallskip
% In other words,  for $\alpha\in (0,\rr(0)]$ we have $q_{\alpha,u}^\rr=\{ su : s\geq 0,\ \rr(s)=\alpha\}\cup\{0\}$, and if $\rr(0)<1$ then for $\alpha\in (\rr(0),1]$ we have $q_{\alpha,u}^\rr=\{ su : s>0,\ \rr(s)=\alpha\}$.

% In case~(ii), the origin of~$\R^d$ may be a $\rho$-quantile of order~$\alpha>0$ in direction~$u$. Actually, it can be shown that~(a) for~$\rho(t)=t$, the origin is a $\rho$-quantile of order~$\alpha>0$ in direction~$u$ if and only~$\alpha\leq P[\{0\}]$.
% % (and that, in that case, it is the unique such quantile). 
% Moreover, (b) provided that~$\psi_+(0)P[\{0\}]+P[\|Z\|\in (0,\infty)\setminus\mathcal{D}_\rho]=0$ where~$Z$ has distribution~$P$ (a condition that always holds for~$\rho(t)=t^p$ with~$p>1$), the origin cannot be a $\rho$-quantile of order~$\alpha>0$ in direction~$u$, so that all these quantiles then belong to~$\{\lambda u:\lambda> 0\}$ (for the sake of completeness, we prove~(a)--(b) in Appendix~\ref{sec:appendixSpher}; see Proposition~\ref{PropSpherMuInOrigin}).

\begin{landscape}

%%%%%%%%%%%%%%%%%%%%%%%%%%%%%%%%%%
\begin{figure}[h!]
\centering
\includegraphics[width=0.9\linewidth]{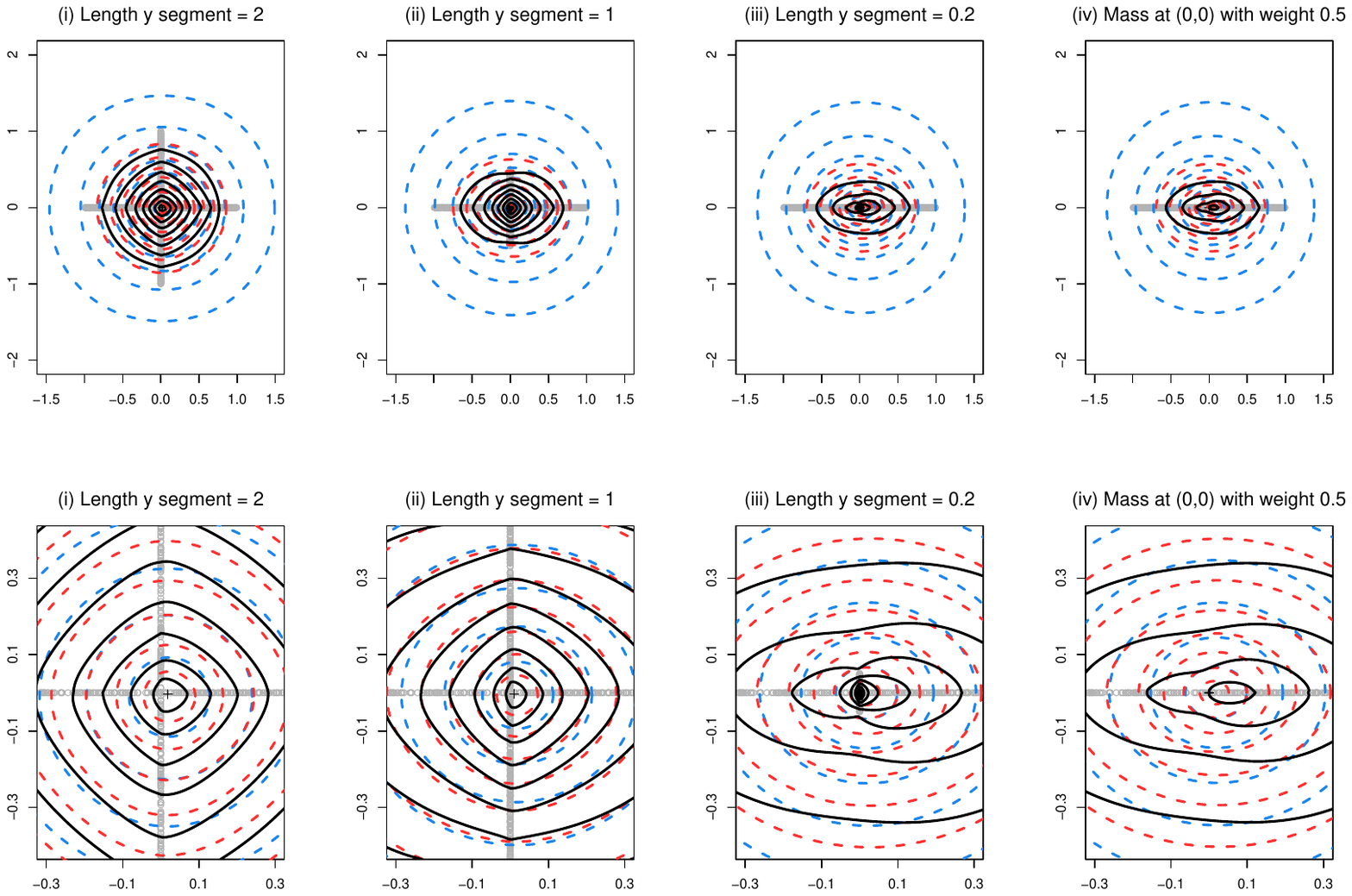}
\vspace{-1cm}
\caption{Graphical representation of the empirical quantiles $q_{\alpha,u}^\rr$, where $\rr(s)=\rr_{\beta}(s)=1-(1+s)^{-\beta}$, for $\beta=2$ (dashed blue curves) and $\beta=5$ (dashed red curves), compared with the classical geometric quantiles $q_{\alpha,u}^{1}$ (full black curves), at levels $\alpha\in \{0.1,0.2,0.3,0.4,0.5,0.6,0.7,0.8\}$. The data points are the same as in Figure~\ref{fig:extreme1to4}. The bottom row is a series of zoomed-in versions of the top panels on the square $[-0.3,0.3]\times [-0.3,0.3]$.}
\label{fig:central1to4}
\end{figure}
%%%%%%%%%%%%%%%%%%%%%%%%%%%%%%%%%%

\end{landscape}

\begin{landscape}

%%%%%%%%%%%%%%%%%%%%%%%%%%%%%%%%%%
\begin{figure}[h!]
\centering
\includegraphics[width=0.9\linewidth]{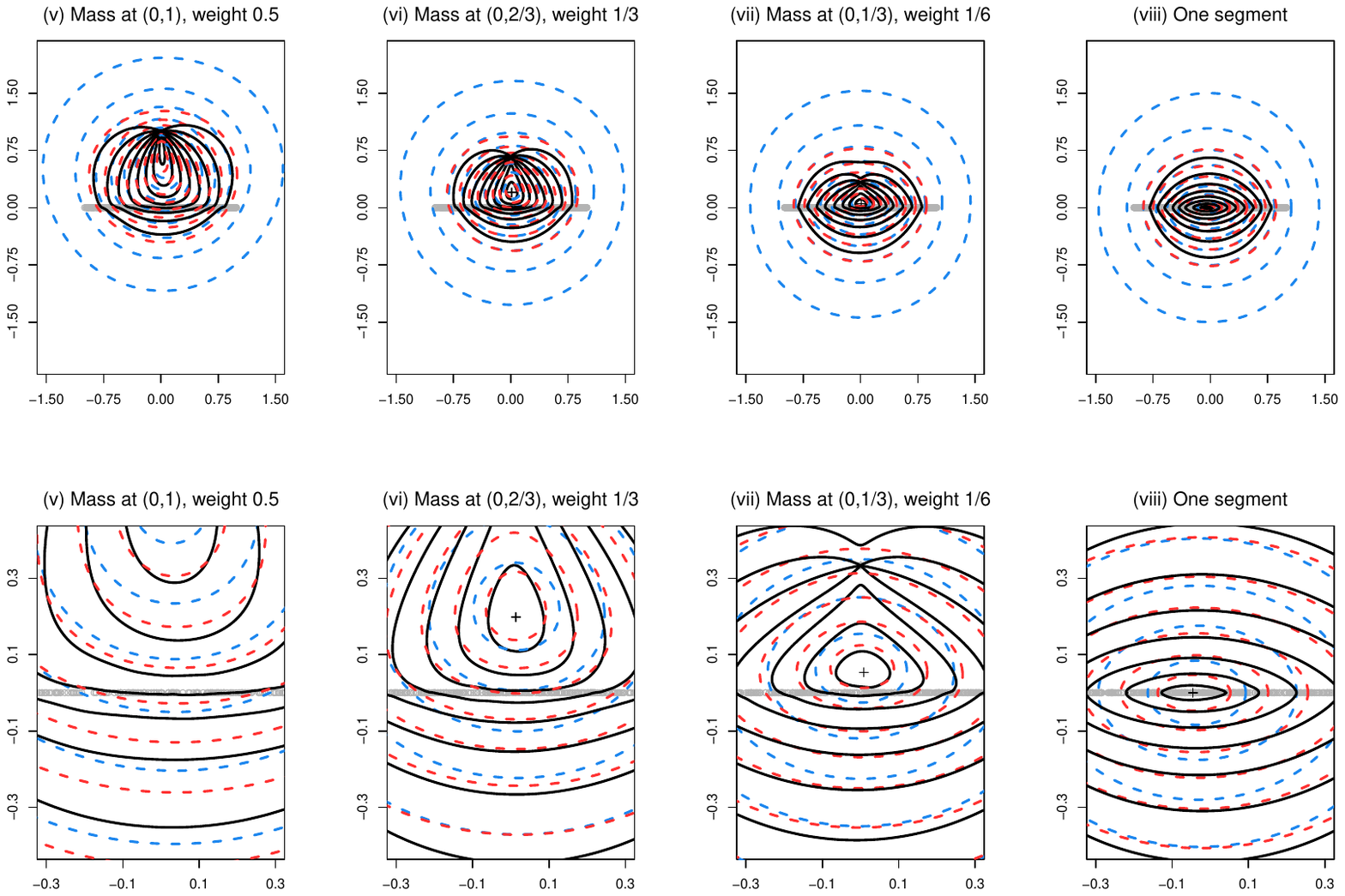}
\vspace{-1cm}
\caption{Graphical representation of the empirical quantiles $q_{\alpha,u}^\rr$, where $\rr(s)=\rr_{\beta}(s)=1-(1+s)^{-\beta}$, for $\beta=2$ (dashed blue curves) and $\beta=5$ (dashed red curves), compared with the classical geometric quantiles $q_{\alpha,u}^{1}$ (full black curves), at levels $\alpha\in \{0.1,0.2,0.3,0.4,0.5,0.6,0.7,0.8\}$. The data points are the same as in Figure~\ref{fig:extreme5to8}. The bottom row is a series of zoomed-in versions of the top panels on the square $[-0.3,0.3]\times [-0.3,0.3]$.}
\label{fig:central5to8}
\end{figure}
%%%%%%%%%%%%%%%%%%%%%%%%%%%%%%%%%%

\end{landscape}

Figures~\ref{fig:central1to4} and~\ref{fig:central5to8} illustrate the sensitivity of $\rr$-geometric quantiles to a lack of symmetry of the underlying distribution in the same setting as Figures~\ref{fig:extreme1to4} and~\ref{fig:extreme5to8}, where data are drawn from a symmetric distribution but, as empirical samples, are not strictly speaking symmetric themselves. It is apparent from these figures that classical geometric quantiles are much more sensitive to lack of symmetry than their regularized counterparts, especially at small quantile levels $\alpha$. This indicates that tests of symmetry based on classical geometric quantiles should be more powerful in practice than those based on $\rr$-geometric quantiles.

\section{Mapping properties and the black hole phenomenon}\label{sec:QuantileMap}

In this section, we investigate the mapping properties of our regularized geometric quantiles, in the spirit of \cite{Kol1997} and \cite{KonPai1}. When $\rr$-quantiles are unique, we start by defining the associated $\rr$-quantile map.

\begin{Defin}
    Let $P\in\PPP(\R^d)$ and $\rr\in\RRR$. Assume that $P$ is not supported on a line or that $\rr$ is increasing. Since $q_{\alpha,u}^\rr(P)$ is unique for all $\alpha\in [0,1)$ and $u\in\S^{d-1}$, we define the \emph{$\rr$-geometric quantile map} of $P$ as $\Qpr:\B^d\to\R^d,\ \alpha u\mapsto q_{\alpha,u}^\rr(P)$.
\end{Defin}

The next theorem is an extension of the continuity result of the classical geometric quantile map for `non-univariate' probability measures. Introducing a regularization, however, allows us to obtain such a continuity property even when the underlying probability is degenerate and supported on a single line provided $\rr$ is (strictly) increasing over $[0,\infty)$. We recall that the distribution function $\Fpr$ was introduced in Definition~\ref{DefinFpr}.

\begin{Theor}\label{TheorContQuant}
    Let $P\in\PPP(\R^d)$ and $\rr\in\RRR$. Assume that $P$ is not supported on a line or that~$\rr$ is increasing. Then, the quantile map $\Qpr:\B^d\to\R^d$ is continuous, and the distribution function $\Fpr:\R^d\to\B^d$ is injective.
\end{Theor}

\cite{GirStu2017} observed that traditional geometric quantiles associated with an increasing order $\alpha_n\to 1$ exit any bounded set as $n\to\infty$, provided the underlying probability measure is non-atomic and not supported on a single line; see Theorem 2.1 there. This (at first surprising) phenomenon is due to the non-existence of extreme geometric quantiles associated with an order $\alpha=1$, even for distributions with bounded support, and compactness of the closed unit ball in finite dimensions combined with the continuity of the geometric quantile map. Theorem 6.2 in \cite{KonPai1} shows that the same holds for $\rho$-quantiles associated with a loss function $\rho(s)=s^p$ for $p\in [1,2]$ or, more generally, such that $s\mapsto s^2/\rho(s)$ is convex over $(0,\infty)$. In our context, the choice of the regularizing effect $\rr$ can produce $\rr$-quantiles that share this property even when the probability measure has atoms and may be supported on a single line.

\begin{Theor}\label{TheorHomeoQuantiles}
Let $P\in\PPP(\R^d)$ and $\rr\in\RRR$. Assume that $P$ is not supported on a line or that~$\rr$ is increasing. Further assume that $\rr(0)=0$ or that $P$ is non-atomic. Then, the quantile map~$\Qpr:\B^d\to\R^d$ is a homeomorphism with (continuous) inverse $\Fpr:\R^d\to\B^d$. 
% In addition, if $P$ is non-atomic then~$\Fpr$ is continuous; in particular, $\Qpr$ and $\Fpr$ are homeomorphisms.
\end{Theor}

This result highlights equally important aspects. On the one hand, extreme $\rr$-quantiles will eventually exit any bounded set. On the other hand, $\rr$-quantiles $q_{\alpha,u}^\rr$ associated with any order $\alpha\in [0,1)$ and direction $u\in\S^{d-1}$ span the whole space, and the inverse of the quantile map is the $\rr$-distribution function $\Fpr$, %introduced in Definition \ref{DefinFpr}, 
which is available in closed form as a simple expectation.\smallskip

Theorem \ref{TheorHomeoQuantiles} requires that $P$ be non-atomic, hence does not seem to apply to atomic (hence, in particular, empirical) distributions. One might think that this is an artifact of the proof techniques used to establish these results, since real-data examples and simulations show that extreme quantiles do have increasingly large norms and explore the whole space even for discrete distributions. In the remainder of this section, we investigate this potential dichotomy between the atomic and non-atomic cases. We start by observing that, when $\rr(0)>0$ and $P$ has atoms, Theorem \ref{TheorHomeoQuantiles} never holds % First, it is straightforward to see that 
because $\Fpr$ is continuous at a point $z\in\R^d$ if and only if $\rr(0)P[\{z\}]=0$. %to the effect that $\Fpr$ is not globally continuous on $\R^d$ due the presence of atoms. 
Let us put this fact aside to focus, instead, on invertibility and further assume that $P$ is not supported on a line or that $\rr$ is increasing, so that all $\rr$-geometric quantiles exist and are unique (Theorem \ref{TheorUniqueness}). The characterization (\ref{eq:DirCond}) entails that if an atom $x\in\R^d$ is a $\rr$-quantile of order $\alpha\in [0,1)$ in direction $u\in\S^{d-1}$, then $x$ remains a quantile of order $\beta\in [0,1)$ in direction $v\in \S^{d-1}$ for any $\beta v$ such that $\|\Fpr(x)-\beta v\|\le \rr(0)P[\{x\}]$, \ie for all $\beta v$ in the closed ball $B_x\subset\B^d$ centered at $\Fpr(x)$ with radius $\rr(0)P[\{x\}]>0$ (observe that $B_x$ is indeed a subset of $\B^d$, for if that were not the case then there would exist $u\in B_x$ intersecting the boundary of~$\B^d$, \ie $\|u\|=1$, hence $x$ would be an extreme quantile which would contradict Proposition \ref{PropQuantileOrderOne}). This means that an atom $x\in\R^d$ `attracts' quantile orders in a closed neighborhood $B_x$ of $\Fpr(x)$, and the larger the mass of $x$ the larger the size of this neighborhood, which is why we refer to this as the `black hole' phenomenon.  This readily implies that the $\rr$-quantile map $\Qpr$ of $P$ is not injective over $\B^d$. This also shows that the $\rr$-distribution function $\Fpr$ is not surjective onto $\B^d$. Indeed, would there exist $y\in\R^d$ such that $\Fpr(y)\in B_x\sm\{\Fpr(x)\}$, then letting $v_y\equiv \Fpr(y)\in\B^d$ would yield $\|\Fpr(y)-v_y\|=0(\le \rr(0)\P[\{y\}])$ so that $\Qpr(v_y)=y$ by virtue of (\ref{eq:DirCond}). But $\Qpr(v_y)=x$ since $v_y\in B_x$, so that $x=y$ by uniqueness of $\rr$-quantiles, which contradicts the fact that $\Fpr(x)\ne \Fpr(y)$. Departure from surjectivity of $\Fpr$ is thus exclusively caused by its discontinuities at atoms.\smallskip

Recalling that $\Fpr$ is injective under our assumptions (Theorem~\ref{TheorContQuant}), then by removing atoms from the picture one might hope to achieve invertibility between appropriately punctured sets. Letting $\AA$ denote the set of atoms of $P$, the next result establishes that $\Fpr$ is a homeomorphism between $\R^d\sm \AA$ and $\B^d\sm \cup_{x\in \AA} B_x$. In addition, when atoms are not removed, invertibility is in fact preserved (while continuity unavoidably fails at atoms) provided $\B^d\sm\cup_{x\in\AA} B_x$ is augmented with the values of $\Fpr$ at the atoms of $P$.

 \begin{Theor}\label{TheorHomeoGeneral}
Let $P\in\PPP(\R^d)$ and $\rr\in\RRR$. Assume that $P$ is not supported on a line or that $\rr$ is increasing. Let $\AA$ denote the set of atoms of $P$. Further assume that $\rr(0)>0$ and that~$\AA$ is non-empty. For all $x\in \AA$, let $B_x$ be the closed ball centered at $\Fpr(x)$ with radius $\rr(0)P[\{x\}]$, and $B_x^\circ$ be the punctured ball $B_x\sm\{\Fpr(x)\}$. Then, 
\begin{enumerate}[label=(\roman*)]
\item $\Fpr$ is a homeomorphism between $\R^d\sm\AA$ and $\B^d\sm\cup_{x\in\AA} B_x$ with inverse $\Qpr$. 
\item $\Fpr$ is a bijection between $\R^d$ and~$\B^d\sm \cup_{x\in\AA} B_x^\circ$ with inverse $\Qpr$.
\end{enumerate}
\end{Theor}

\begin{figure}[h!]
\includegraphics[width=\linewidth]{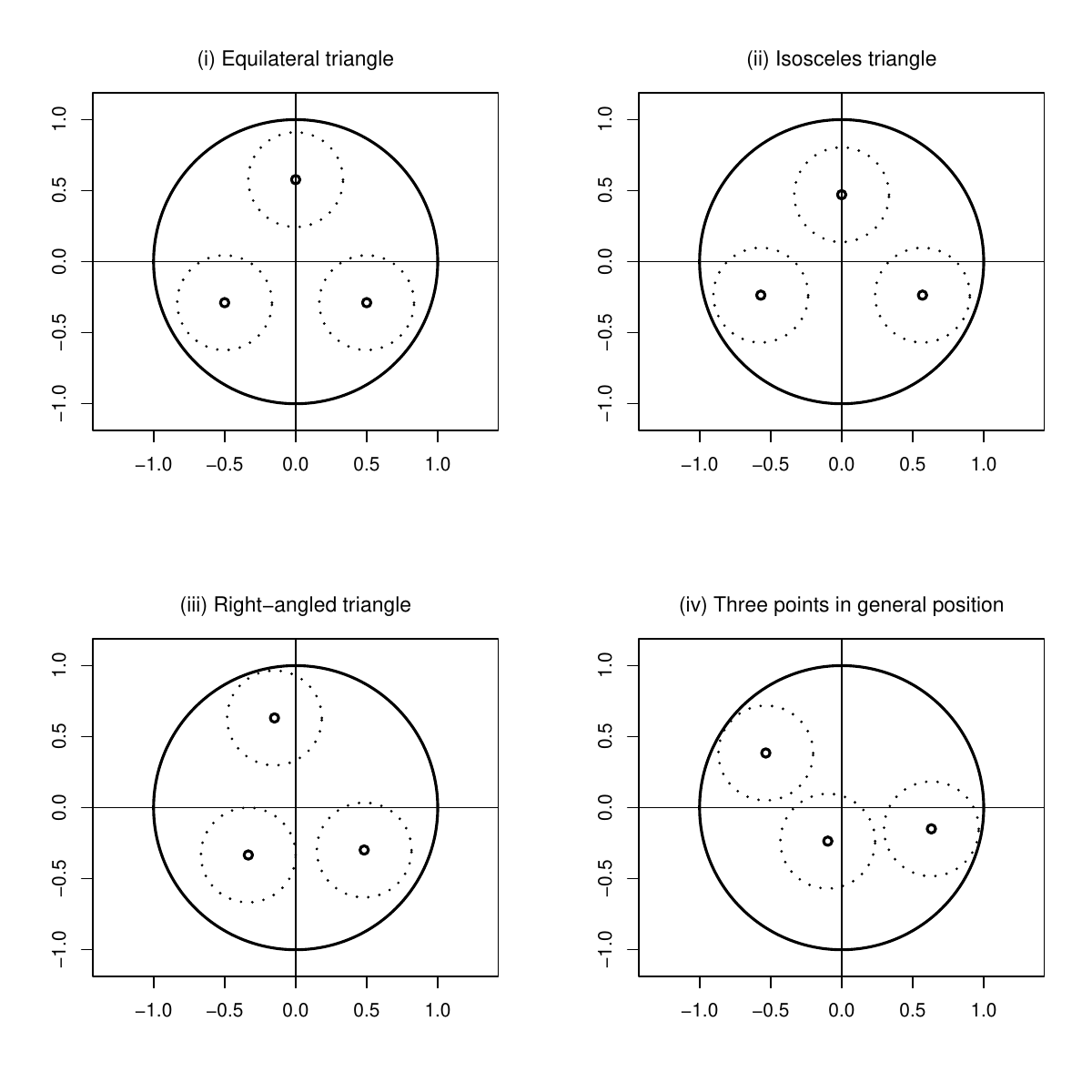}
\caption{Black holes of classical geometric quantiles corresponding to four uniform distributions on the vertices of triangles in $\mathbb{R}^2$, denoted by $ABC$, with (i) $A(-1,-1/\sqrt{3})$, $B(1,-1/\sqrt{3})$ and $C(0,\sqrt{3}-1/\sqrt{3})$ (equilateral triangle), (ii) $A(-1,-1/3)$, $B(1,-1/3)$ and $C(0,2/3)$ (isosceles triangle), (iii) $A(-1/3,-2/3)$, $B(2/3,-2/3)$ and $C(-1/3,4/3)$ (right-angled triangle) and (iv) $A(0,-1/3)$, $B(1,-1/3)$ and $C(-1,2/3)$ (points in general position). The three black holes are delimited by a dashed circle and centered at points located at the black circular marks.}
\label{fig:threedirac_Fpmap}
\end{figure}

\begin{landscape}

%%%%%%%%%%%%%%%%%%%%%%%%%%%%%%%%%%
\begin{figure}[h!]
\centering
\includegraphics[width=0.95\linewidth]{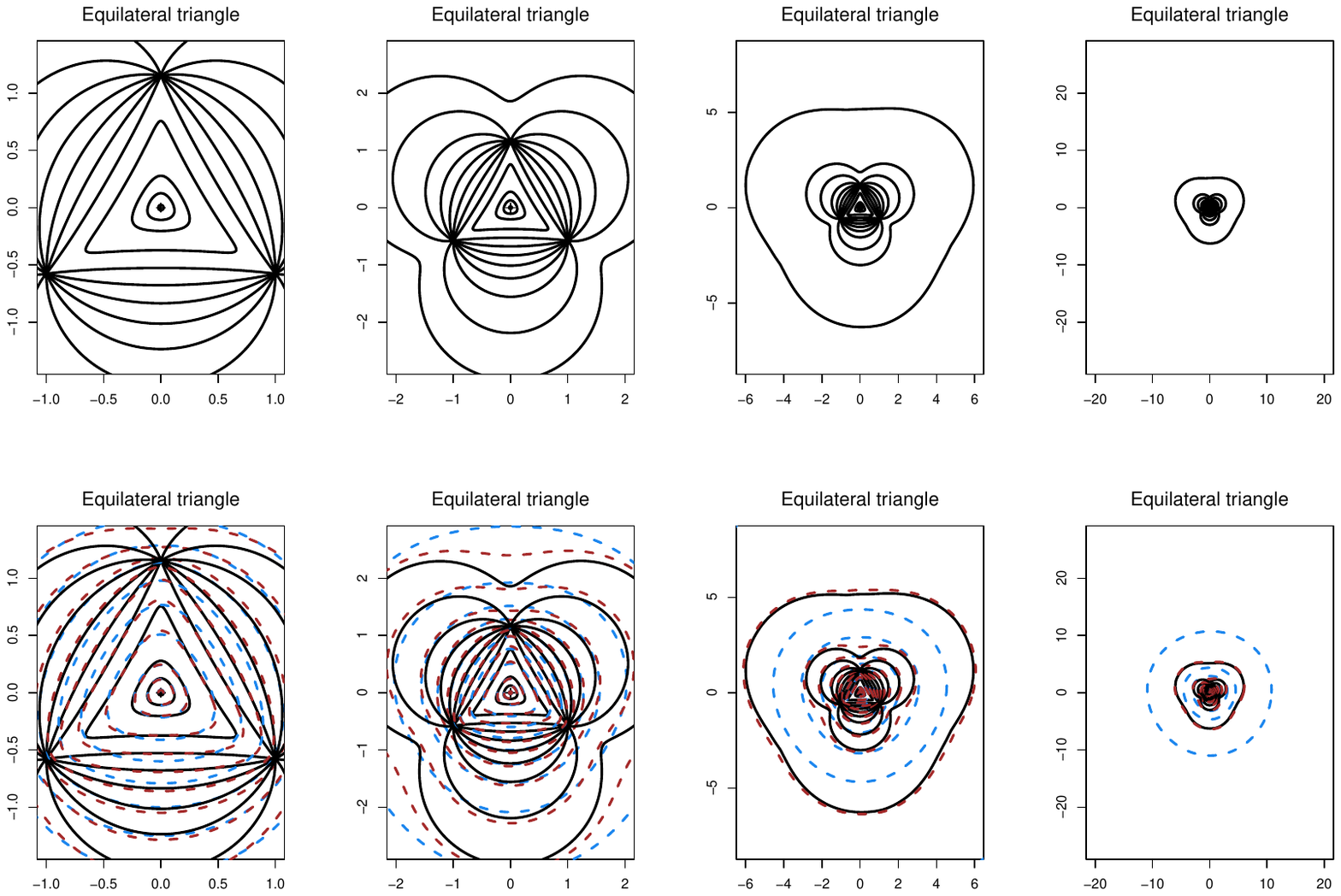}
\vspace{-1cm}
\caption{Classical geometric quantiles and $\rr_{\beta}$-geometric quantiles relative to the first uniform distribution considered in Figure~\ref{fig:threedirac_Fpmap}, {\it i.e.}~the uniform distribution on the vertices of the triangle $ABC$ in $\mathbb{R}^2$, with $A(-1,-1/\sqrt{3})$, $B(1,-1/\sqrt{3})$ and $C(0,\sqrt{3}-1/\sqrt{3})$ (equilateral triangle). Top row: Geometric contours only, bottom row: Geometric contours with superimposed $\rr_{\beta}$-geometric quantile contours for $\beta=2$ (dashed blue lines) and $\beta=4$ (dashed brown lines). The contours drawn are those corresponding to levels $\alpha\in \{0.01, 0.05, 0.1, 0.2, 0.3, 0.4, 0.5, 0.6, 0.7, 0.8, 0.9, 0.95, 0.99 \}$.}
\label{fig:threedirac_1}
\end{figure}
%%%%%%%%%%%%%%%%%%%%%%%%%%%%%%%%%%

\end{landscape}

\begin{landscape}

%%%%%%%%%%%%%%%%%%%%%%%%%%%%%%%%%%
\begin{figure}[h!]
\centering
\includegraphics[width=0.95\linewidth]{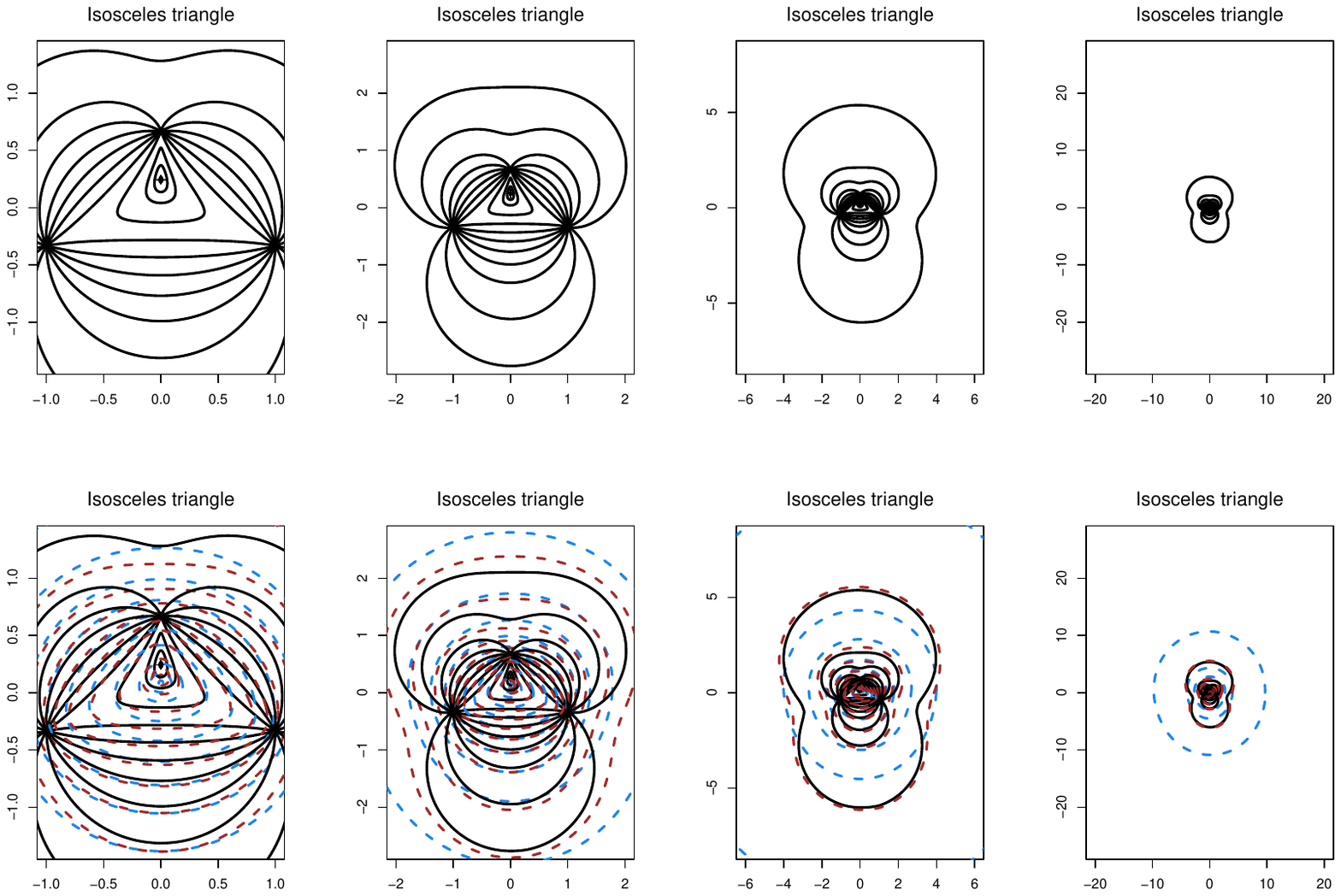}
\vspace{-1cm}
\caption{Classical geometric quantiles and $\rr_{\beta}$-geometric quantiles relative to the second uniform distribution considered in Figure~\ref{fig:threedirac_Fpmap}, {\it i.e.}~the uniform distribution on the vertices of the triangle $ABC$ in $\mathbb{R}^2$, with $A(-1,-1/3)$, $B(1,-1/3)$ and $C(0,2/3)$ (isosceles triangle). Top row: Geometric contours only, bottom row: Geometric contours with superimposed $\rr_{\beta}$-geometric quantile contours for $\beta=2$ (dashed blue lines) and $\beta=4$ (dashed brown lines). The contours drawn are those corresponding to levels $\alpha\in \{0.01, 0.05, 0.1, 0.2, 0.3, 0.4, 0.5, 0.6, 0.7, 0.8, 0.9, 0.95, 0.99 \}$.}
\label{fig:threedirac_2}
\end{figure}
%%%%%%%%%%%%%%%%%%%%%%%%%%%%%%%%%%

\end{landscape}

\begin{landscape}

%%%%%%%%%%%%%%%%%%%%%%%%%%%%%%%%%%
\begin{figure}[h!]
\centering
\includegraphics[width=0.95\linewidth]{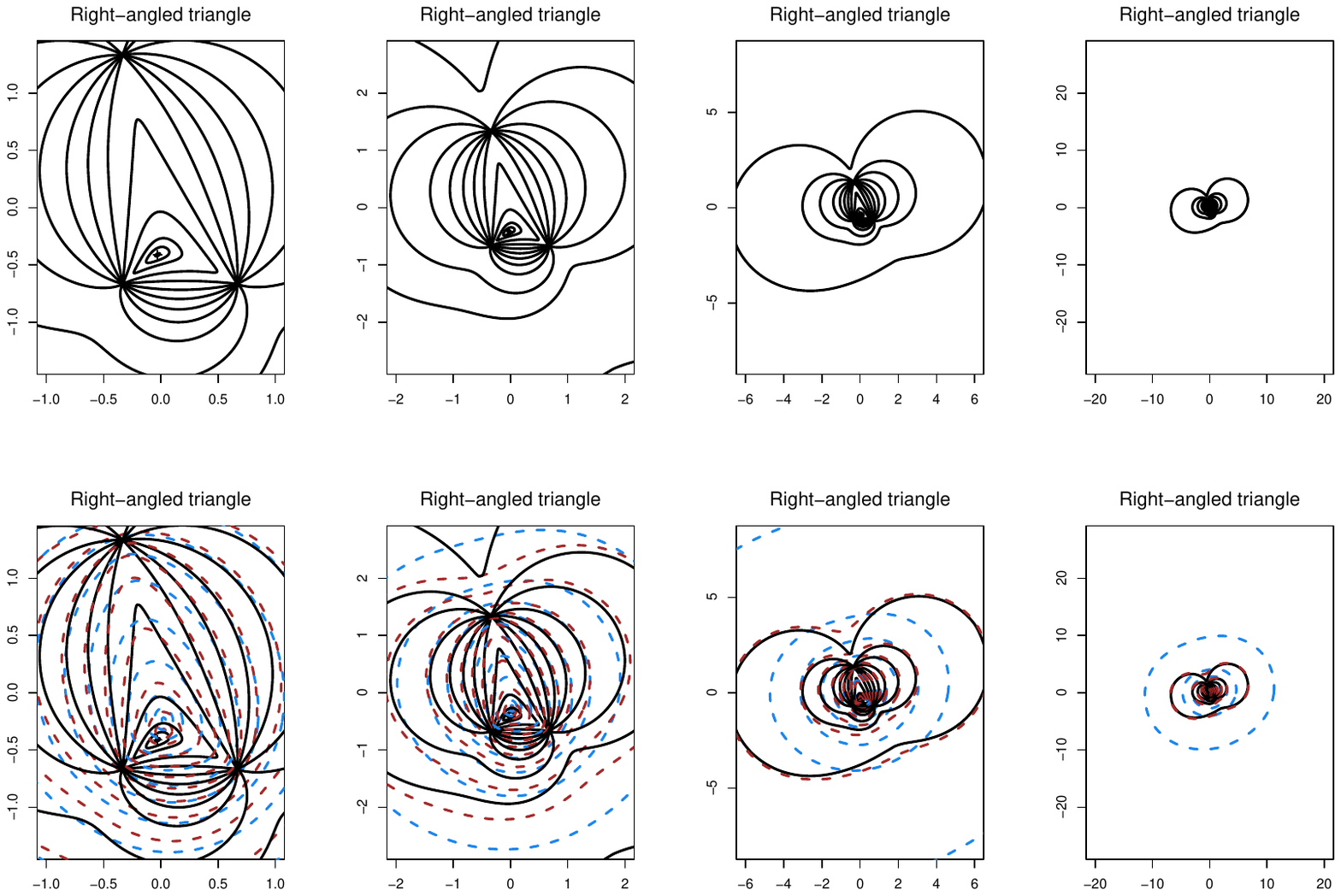}
\vspace{-1cm}
\caption{Classical geometric quantiles and $\rr_{\beta}$-geometric quantiles relative to the third uniform distribution considered in Figure~\ref{fig:threedirac_Fpmap}, {\it i.e.}~the uniform distribution on the vertices of the triangle $ABC$ in $\mathbb{R}^2$, with $A(-1/3,-2/3)$, $B(2/3,-2/3)$ and $C(-1/3,4/3)$ (right-angled triangle). Top row: Geometric contours only, bottom row: Geometric contours with superimposed $\rr_{\beta}$-geometric quantile contours for $\beta=2$ (dashed blue lines) and $\beta=4$ (dashed brown lines). The contours drawn are those corresponding to levels $\alpha\in \{0.01, 0.05, 0.1, 0.2, 0.3, 0.4, 0.5, 0.6, 0.7, 0.8, 0.9, 0.95, 0.99 \}$.}
\label{fig:threedirac_3}
\end{figure}
%%%%%%%%%%%%%%%%%%%%%%%%%%%%%%%%%%

\end{landscape}

\begin{landscape}

%%%%%%%%%%%%%%%%%%%%%%%%%%%%%%%%%%
\begin{figure}[h!]
\centering
\includegraphics[width=0.95\linewidth]{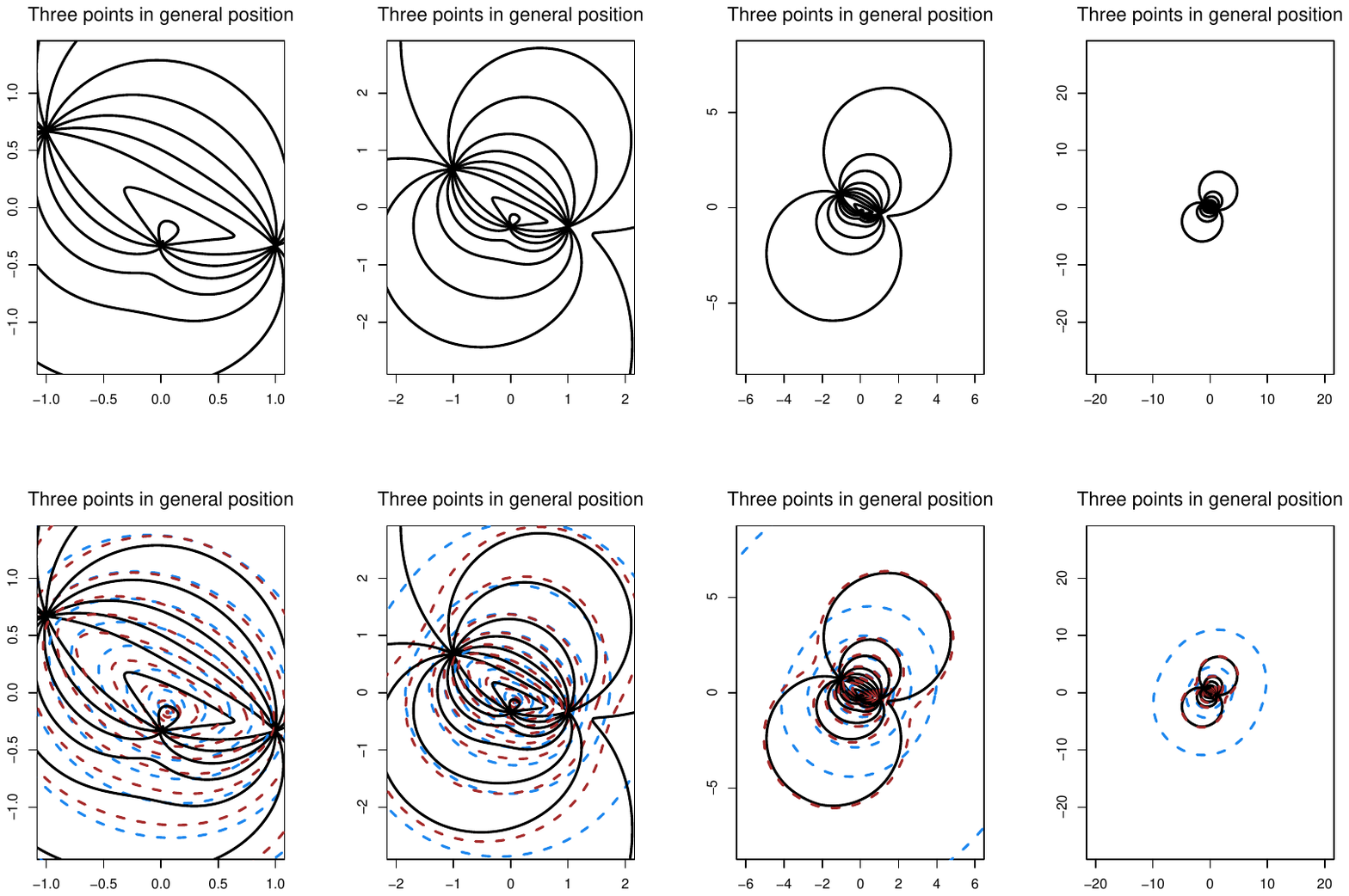}
\vspace{-1cm}
\caption{Classical geometric quantiles and $\rr_{\beta}$-geometric quantiles relative to the fourth uniform distribution considered in Figure~\ref{fig:threedirac_Fpmap}, {\it i.e.}~the uniform distribution on the vertices of the triangle $ABC$ in $\mathbb{R}^2$, with $A(0,-1/3)$, $B(1,-1/3)$ and $C(-1,2/3)$ (points in general position). Top row: Geometric contours only, bottom row: Geometric contours with superimposed $\rr_{\beta}$-geometric quantile contours for $\beta=2$ (dashed blue lines) and $\beta=4$ (dashed brown lines). The contours drawn are those corresponding to levels $\alpha\in \{0.01, 0.05, 0.1, 0.2, 0.3, 0.4, 0.5, 0.6, 0.7, 0.8, 0.9, 0.95, 0.99 \}$.}
\label{fig:threedirac_4}
\end{figure}
%%%%%%%%%%%%%%%%%%%%%%%%%%%%%%%%%%

\end{landscape}

The black hole phenomenon is illustrated in Figures~\ref{fig:threedirac_Fpmap} to~\ref{fig:threedirac_4}, in the case of the usual geometric quantiles obtained with $\rr\equiv 1$. Figure~\ref{fig:threedirac_Fpmap} shows the black holes for four uniform distributions on the vertices of triangles in $\mathbb{R}^2$, and Figures~\ref{fig:threedirac_1} to~\ref{fig:threedirac_4} represent the geometric contours of these uniform distributions for levels $\alpha$ ranging between $0.01$ and $0.99$. The attraction phenomenon at the vertices of the triangles is clearly seen. For instance, in Figure~\ref{fig:threedirac_Fpmap}(i), geometric quantiles in the directions of the three vertices are progressively trapped by the vertices until $\alpha$ is so large that they escape. In Figure~\ref{fig:threedirac_Fpmap}(iv), geometric quantiles associated with a small level $\alpha$ are all equal to a single point and progressively become nondegenerate quantile contours as $\alpha$ increases, because the black hole corresponding to one of the vertices contains the origin. By contrast, none of the $\rr_{\beta}$-quantile contours is affected by the black hole phenomenon since $\rr_{\beta}(0)=0$, as can be seen in the bottom panels of Figures~\ref{fig:threedirac_1} to~\ref{fig:threedirac_4} computed for $\beta=2$ and $\beta=4$.\smallskip

% From Theorem \ref{TheorHomeoGeneral}, we recover Theorem \ref{TheorHomeoQuantiles} as particular case since, if $P$ is non-atomic or $\rr(0)=0$, then $B_x=\{\Fpr(x)\}$ and $B_x^\circ=\emptyset$, so that $\Fpr$ is a bijection between $\R^d$ and $\B^d$, whereas if $P$ is non-atomic then $\AA=\emptyset$ so that $\Fpr$ is a homeomorphism between $\R^d$ and $\B^d$.

    % \gs{Two ways of illustrating: either draw quantile curves for a lot of values of $\alpha$ in an atomic case (triangle...) and they should stagnate, or draw a vector field (for each $v$ draw $q(v)$).}

We now use Theorem \ref{TheorHomeoGeneral} to establish a few inequalities which, to the best of our knowledge, were unknown in the literature even in the case of geometric quantiles. When $P$ is not supported on a line or $\rr$ is increasing, define the black hole $B_x$ around an arbitrary $x\in\R^d$ as above if $x$ is an atom of $P$ and as $B_x=\{\Fpr(x)\}$ otherwise.  Observe that the $B_x$'s are disjoint balls: if there exists $v\in B_x\cap B_y$, %with $x\ne y$, 
then $\|\Fpr(x)-v\|\le \rr(0)P[\{x\}]$ and $\|\Fpr(y)-v\|\le \rr(0)P[\{y\}]$ by virtue of (\ref{eq:DirCond}), to the effect that $x$ and $y$ are both $\rr$-geometric quantiles of order $\|v\|$ in direction $v/\|v\|$ (the direction may be taken arbitrarily if $v=0$), which yields $x=y$ by uniqueness of $\rr$-quantiles (Theorem \ref{TheorUniqueness}). %, hence contradicts the fact that $x\ne y$. 
Since $B_x$ and $B_y$ are closed balls centered at $\Fpr(x)$ and $\Fpr(y)$ with radius $\rr(0)P[\{x\}]$ and $\rr(0)P[\{y\}]$, respectively, this provides
\begin{equation}\label{DiffRanks}
\|\Fpr(x)-\Fpr(y)\|
> 
\rr(0)\big(P[\{x\}]+P[\{y\}]\big) 
% =:
% \kappa_{x,y}
,\spa 
\forall\ x\neq y
.
\end{equation}

\begin{Example}
Consider the empirical distribution $P_n$ of a sample 
$
\XX_n=\{X_1,\ldots, X_n\}
$ 
drawn independently from a common distribution $P\in\PPP(\R^d)$, and assume that the $X_i$'s are pairwise distinct and not collinear. Denote the corresponding black hole centered at $F_{P_n}^\rr(X_i)$ with radius $\rr(0)/n$ as $B_i$. Then $F_{P_n}^\rr$ is a homeomorphism of $\R^d\sm \XX_n$ onto $\B^d\sm \cup_{i=1}^n B_i$. This result explains why empirical distributions share the qualitative properties of continuous ones outside $\XX_n$ and $\cup_{i=1}^n B_i$, such as invertibility and continuity of the $\rr$-distribution function and $\rr$-quantile map. We can further quantify the rate of the `phase transition' from the properties of $P_n$ to those of $P$. Indeed, recall that the $B_i$'s are disjoint balls of radius $\rr(0)/n$. Therefore, denoting by $\omega_d\equiv \lambda(\B^d)$ the Lebesgue measure of $\B^d$, we find 
$$
\lambda(\B^d\sm \cup_{i=1} B_i) 
=
\omega_d\Big(1 - n \Big(\frac{\rr(0)}{n}\Big)^d\Big) 
=
\omega_d\Big( 1 -  \frac{\rr(0)^d}{n^{d-1}} \Big) 
.
$$
Consequently, if one samples from $\B^d$ according to a distribution admitting a bounded density then, with probability of the order $1-n^{-(d-1)}$, the resulting samples will belong to a region of~$\B^d$ where $F_{P_n}^\rr$ and $Q_{P_n}^\rr$ are homeomorphic to each other, thus quantifying the `phase transition' between the empirical and population properties.
% that the $\rr$-quantile map $\Qpr$ defines a homeomorphism over a subset of $\B^d$ whose Lebesgue measure increases
Maybe even more important is the fact that there is no phase transition at all when $\rr(0)=0$ since, in this case, the $\rr$-quantile map and the $\rr$-distribution function are invertible and inverse of one another, even if $P_n$ is purely discrete. Finally, when $\rr(0)>0$, then (\ref{DiffRanks}) yields
$$
 \|F_{P_n}^{\rr}(X_i)-F_{P_n}^{\rr}(X_j)\|
 \geq 
 \frac{2\rr(0)}{n}
 ,\spa 
 \forall\ i\ne j
 .
 $$
 This is the $\rr$-geometric analogue of the equi-spacing of the ranks of univariate distributions, although the $\rr$-geometric ranks $F_{P_n}^\rr(X_1),\ldots, F_{P_n}^\rr(X_n)$ are no longer uniformly distributed over~$\B^d$ when $d\ge 2$; see Section~4.1 in \cite{KonHal24} for a discussion.
    % Comparison with univariate. Sampling rate in $\B^d$
    % This yields 
    % $$
    % \Big( \|\Fpr(x)\| - \rr(0)P[\{x\}] \Big) + \Big( \|\Fpr(y)\| - \rr(0) P[\{y\}] \Big)
    % \geq 
    % 0
    % ,\spa 
    % \forall\ x,y\in\AA
    % ,
    % $$
    % so that taking $y=x$
   % Also, if the $\rr$-geometric median $m_P^\rr$ of $P$ is an atom, then
   %  $$
   %  \|\Fpr(x)\|
   %  \geq 
   %  \rr(0)\big(P[\{x\}] + P[\{m_P^\rr\}]\big)
   %  .
   %  $$
    % so that
    % $$
    % \|\Fpr(x)\|
    % \geq 
    % \rr(0) P[\{x\}]
    % ,\spa 
    % \forall\ x\in\AA 
    % .
    % $$
    
% In empirical case $X_1,\ldots, X_N$, it suffices to list all atoms $X_1,\ldots, X_N$, compute their corresponding image in $\B^d$ through $U_i=F_{P_n}^{\rr}(X_i)$ for $i=1,\ldots, N$. Letting $B_N$ be the punctured ball $\B^{(N)} = \B^d\sm \cup_{i=1}^N \overline{\B}_{\rr(0)/N}(U_i)$, $\rr$-geometric quantiles $q_{\alpha,u}^{\rr,(N)}$ will be unique as soon as $\alpha u\in B_N$. There are $N$ atoms each forming a blackhole of radius $\rr(0)/N$ so that the total blackhole volume does not exceed 
% $$
% {\rm vol}(\B^{(N)})
% = 
% N \bigg(\frac{\rr(0)}{N}\bigg)^d \omega_d
% = 
% \omega_d\ \rr(0) \frac{1}{N^{d-1}}
% ,
% $$
% where $\omega_d$ is the volume of the unit ball of $\R^d$. In particular, $Q_{P_n}^{\rr}$ is a homeomorphism over a punctured ball $\B^{(N)}$ such that that is exponentially close to $\B^d$ as $N\to\infty$: in the sense of volumes, ${\rm vol}(\B^d\sm B_N)\leq \omega_d \rr(0) / N^{d-1}$, in Hausdorff sense?\\
\end{Example}

    % The values of $B_i\sm \{U_i\}$ do not belong to the image of $\Fpr$. Indeed, if there exists $x$ such that $v\equiv \Fpr(x)\in B_i\sm\{U_i\}$, then $\|\Fpr(x)-v\|=0\leq \rr(0)P[\{x\}]$, so that $x$ is a quantile of order $v$. But the only quantile of order $v$ is $X_i$, so we must gave $x=X_i$, which contradicts the fact that $\Fpr(x)\neq U_i$. So $F_{P_n}$ will be bijective between $\R^d\sm A_{P_n}$ and $\B^{(N)}$, and between $\R^d$ and $\B^{(N)}\cup F_{P_n}^{\rr}(A_{P_n})=\B^{(N)}\cup \{U_1,\ldots, U_N\}$. And $Q_{P_n}^{\rr}$ is surjective between $\B^d$ and $\R^d$, but constant (equal to $X_i$) on $B_i$. 

Taking %$y$ in (\ref{DiffRanks}) equal to the $\rr$-geometric median of $P$, $m_P^\rr$ say, the characterization (\ref{eq:DirCond}) yields 
$x$ in (\ref{eq:DirCond}) equal to the $\rr$-geometric median of $P$, $m_P^\rr$ say, we find $\rr(0)P[\{m_P^\rr\}]-\|\Fpr(m_P^\rr)\|\geq 0$. Combining this with (\ref{DiffRanks}) for $y=m_P^\rr$, we deduce that from the triangle inequality that
\begin{equation}\label{eq:FLowerBound}
\|\Fpr(x)\|
> 
\rr(0) P[\{x\}]
,\spa 
\forall\ x\in\R^d\sm\{m_P^\rr\}
.
\end{equation}
Since $P$ is not supported on a line or $\rr$ is increasing, this ensures that $\rr(0)P[\{x\}]<1$. Because $B_x$ is centered at $\Fpr(x)$, has radius $\rr(0)P[\{x\}]\in [0,1)$, and is a closed subset of the open unit ball~$\B^d$, we deduce that 
\begin{equation}\label{eq:FUpperBound}
\|\Fpr(x) \|
<
1-\rr(0)P[\{x\}]
,\spa 
\forall\ x\in\R^d 
.
\end{equation}
Combining (\ref{eq:FUpperBound}) with (\ref{eq:FLowerBound}) yields $\rr(0)P[\{x\}] < 1/2$ for all $x\in\R^d\sm\{m_P^\rr\}$. This also entails that any atom $x\in\R^d$ satisfying $\rr(0) P[\{x\}]\ge 1/2$ is automatically a $\rr$-geometric median of $P$. For such an $x\in\R^d$, the corresponding black hole $B_x$ necessarily contains $(\alpha=)0$ and has radius $\rr(0)P[\{x\}]\ge 1/2$. In particular, we have $\Qpr(v)=x$ for all $v$ in a ball of radius larger than $1/2$ regardless of the size of $\|x\|$, which is quite striking. We illustrate this phenomenon in the following example.

\begin{Example} 
Consider the distribution $P_{X_t}$ of the random (mixture) vector in $\R^2$
$$
X_t
\equiv 
\begin{cases}
    Y & \text{ if } \sigma = 1\\
    (t,0) & \text{ if } \sigma=0
    ,
\end{cases}
$$
where $Y$ is uniformly distributed over the unit disk $\B^2$ and $\sigma\sim {\rm Bern}(1/2)$ is independent of $Y$. For simplicity, consider $\rr\equiv 1$ leading to geometric quantiles. Then, $P_{X_t}$ is not supported on a line so that all geometric quantiles exist and are unique (Theorem \ref{TheorUniqueness}). Thus, the point $(t,0)$ is the (unique) $\rr$-median of $P_{X_t}$, since $P[\{(t,0)\}]=1/2$, and the corresponding black hole $B_{(t,0)}$ has radius $1/2$ and contains $(0,0)$ since $(t,0)$ is the geometric median. When $t=0$, we have $F_{P_{X_0}}((0,0))=(0,0)$ so that the origin remains the unique geometric quantile associated with any $v\in\B^d$ such that $\|v\|\le 1/2$ (this also follows from Proposition \ref{TheorSpher} in this case by sphericity of $X_0$). Using Lemma A.1 in \cite{KonPai2025a} and invariance of $F_{P_{X_t}}((t,0))$ under axial symmetry about the horizontal axis, one shows that 
$$
F_{P_{X_t}}^{\rr}((t,0))
=
\Big(\frac12 + O\Big(\frac{1}{t}\Big), 0\Big)
$$
as $t\to +\infty$, for some non-positive remainder $O(1/t)$. The center of the black hole $B_{(t,0)}$ shifts towards $(1/2,0)$ and, in the limit, the point $(t,0)$ is thus a geometric quantile of (approximately) any order $v\in\B^d$ such that $\|v-(1/2,0)\|\le 1/2$.
% ; this can be strengthened to $O(1/t^2)$ by appealing to Lemma 5.3 in \cite{Romon2022} by symmetry of the support of $Y$. 
% So the center of the black hole $B_{(t,0)}$ shifts towards $(1/2,0)$ and we see that the limiting black hole still contains $(\alpha=)0$.
\end{Example}

    % Resultats d'inversibilité et homéo quand il y a des atomes (Koltchinski, mais jamais formulé comme tel pour des quantiles) :

% $\nu = (d-2)/2$

\section{Behaviour of extreme $\rr$-quantiles}
\label{sec:extremes}

% \textcolor{brown}{{\bf GS:} Include simulations} 

% \textcolor{red}{Define $\RRR$ to be the convex set of sub-distribution functions that are not identically 0 and, for $r\in \RRR$, set $\rr_\infty\equiv \sup_{s>0}r(s)=\lim_{s\to\infty} r(s)>0$. Since the $M$-quantiles $\Qpr(v)$ we consider can only be defined up to order $\|v\|=r_\infty$, a reasonable notion of extreme $M$-quantile is of the form $\Qpr(v_n)$ with $\| v_n \|<r_\infty$ and $\| v_n \| \to r_\infty$.} 

We start by two qualitative results establishing that, unless $P$ is supported on a half-line (in which case, the classical univariate extreme value theory applies), extreme $\rr$-quantiles of $P$ in direction $u$ escape to infinity in a direction parallel to $u$.

\begin{Prop}\label{PropExtremesQuali}
    Let $P\in\PPP(\R^d)$ and $\rr\in\RRR$. Assume that $P$ is not supported on a line or $\rr$ is increasing. For any sequences $(\alpha_k)\subset [0,1)$ and $(u_k)\subset\S^{d-1}$ such that $\alpha_k \to 1$ and $u_k\to u$,
    \[
    \|\Qpr(\alpha_k u_k)\|\to \infty,
    \spa 
    \text{and}
    \spa
    \frac{\Qpr(\alpha_k u_k)}{\|\Qpr(\alpha_k u_k)\|}
    \to 
    u
    \spa 
    \text{as}
    \ k\to\infty.
    \]
    %If $r(t)<1$ for all $t>0$, then the conclusion holds even if $P$ is supported on a halfline.\\
    % , and $r(t)=1$ for all $t>t_1$, then the conclusion also holds if $P$ is not supported on $\R^d\sm B_{t_1}$.
\end{Prop}

This result was first established in \cite{GirStu2017} (Theorem~2.1) for geometric quantiles obtained for $\rr\equiv 1$, then under weaker assumptions in \cite{PaiVir2021} (Theorem~2 and Theorem~3). 
% Even for the usual geometric quantiles, Proposition~\ref{PropExtremesQuali} improves upon previous results and holds under optimal assumptions: it holds for any probability measure that is not supported on a half-line, and fails for any probability measure supported on a half-line by virtue of Proposition~\ref{PropQuantileOrderOne}. 
Our result, however, holds for an \emph{arbitrary} probability measure as soon as $\rr(s)<1$ for all $s\ge 0$, which prevents the existence of extreme quantiles. In particular, we find that, like extreme geometric quantiles, extreme $\rr$-quantiles exit any compact set, even if the underlying distribution has a compact support. This property was shown by~\cite{KonPai1} to be also shared by extreme multivariate $L^p$-quantiles for $p\in [1,2]$, which includes the class of extreme multivariate expectiles studied by~\cite*{Heretal2018}. From an extreme value perspective, this is of course a drawback of geometric quantiles; in the univariate case, the asymptotic behavior of extreme quantiles of a given distribution depends on its tail characteristics in a non-trivial way (see, e.g., Theorem~1.1.6 in~\cite*{HaaFer2006}). This feature, however, turns out to be useful in some inferential applications, such as supervised classification based on the max-depth approach, where one leverages the fact that the quantile-based depth of any location in $\R^d$ (even lying outside the bulk of the data) is non-zero; see, e.g., \cite*{GhoCha2005A}, \cite*{Lietal2012},  \cite*{Franetal2019}, and \cite{KonPai23}.\smallskip

 Theorem 2.2 in~\cite{GirStu2017} implies that the growth rate of extreme geometric quantiles of probability distributions with finite second moment is determined by their covariance matrix, and the size of an extreme geometric quantile is asymptotically the largest in directions where the variance is the smallest. For $\rr$-quantiles obtained with $\rr(s)<1$, the loss function featuring in the minimization problem defining the corresponding $\rr$-quantiles puts less weight on the center of the distribution (see Definition \ref{DefinQuantiles}). In what follows, we thus investigate the asymptotic behavior of extreme $\rr$-quantiles and whether it might reveal additional extreme value characteristics on the extreme value behavior of multivariate distributions. We start by introducing a quantity that will play a key role. For all $\beta>0$, define (when it exists)
\begin{equation}\label{defEll}
\ell_{\rr}(\beta) 
\equiv 
\lim_{s\to \infty} s^\beta \big(1-\rr(s)\big) 
.
\end{equation}
% The quantity $\ell_{\rr}(\beta)$ is intuitively important for the calculation of extreme $M$-quantiles, since it gives information about how quickly $r(s)$ tends to 1 as $s\to\infty$: 
Geometric quantiles, obtained for $\rr\equiv 1$, satisfy $\ell_{\rr}(\beta)=0$ for any $\beta>0$. More generally, $\ell_\rr(\beta)<\infty$ implies that $\rr(s)$ increases to $1$ at least as fast as $s^{\beta}$ when $s\to\infty$. In addition, $\ell_{\rr}(\beta)=\ell\in (0,\infty)$ means that $1-\rr$ is regularly varying at infinity with index $-\beta$ and expands as $\rr(s)=1-\ell s^{-\beta}+{\rm o}(s^{-\beta})$ as $s\to\infty$. Note that depending on the choice of $\rr\in\RRR$, we may have $\ell_{\rr}(\beta)=\infty$ for some values of $\beta$. If, however, $\ell_{\rr}(\beta_0)<\infty$ for some $\beta_0>0$, then $\ell_{\rr}(\beta)=0$ for all $\beta\in (0,\beta_0)$. In particular, there is at most one value of $\beta>0$ such that $\ell_{\rr}(\beta)\in (0,\infty)$.

% \textcolor{red}{We then have the following results characterizing the asymptotic behavior of the direction and norm of an extreme $M$-quantile.}

\begin{Theor}\label{TheoExtremesDirNorm}
    Let $P\in\PPP(\R^d)$ and $\rr\in\RRR$. Assume that $P$ is not supported on a line or that $\rr$ is increasing. Further assume that $P$ has at most finitely many atoms or that $\rr(0)=0$. Fix $\beta>0$ and assume that $\ell_{\rr}(\beta)<\infty$. Let $(\alpha_k)\subset [0,1)$ and $(u_k)\subset\S^{d-1}$ be sequences such that $\alpha_k \to 1$ and $u_k\to u\in \S^{d-1}$ as $k\to\infty$. 
    \begin{enumerate}[label=(\roman*)]
    \item Assume that $P$ has a finite first moment. If $\beta\in (0,1)$, then
    \[
    \| \Qpr(\alpha_k u_k) \|^{\beta} \left( \frac{\Qpr(\alpha_k u_k)}{\| \Qpr(\alpha_k u_k) \|} - \alpha_k u_k \right) \to \ell_{\rr}(\beta) u
    .
    \]
    If $\beta\ge 1$ then, letting $Z$ denote a random $d$-vector with law $P$,
    \[
    \| \Qpr(\alpha_k u_k) \| \left( \frac{\Qpr(\alpha_k u_k)}{\| \Qpr(\alpha_k u_k) \|} - \alpha_k u_k \right) \to \ell_{\rr}(1) u + \E\big[Z-\ps{Z}{u} u\big]
    .
    \]
    Moreover, if $\beta \in (0,1]$, we have
    \[
    \| \Qpr(\alpha_k u_k) \|^{\beta} (1-\alpha_k) \to \ell_{\rr}(\beta)
    . 
    \]
    \item Assume that $P$ has a finite second moment. If $\beta\in (1,2)$, then
    \[
    \| \Qpr(\alpha_k u_k) \|^{\beta} (1-\alpha_k) \to \ell_{\rr}(\beta).
    \]
    If $\beta\ge 2$, then, letting $\Sigma$ denote the covariance matrix of $P$, we have 
    \[
    \| \Qpr(\alpha_k u_k) \|^2 (1-\alpha_k) \to \ell_{\rr}(2) + \frac{1}{2} \big(\operatorname{Tr}(\Sigma) - \ps{u}{\Sigma u}\big) >0
    .
    \]
    \end{enumerate}
\end{Theor}

This readily provides the following directional asymptotics for extreme $\rr$-quantiles.

\begin{Corol}\label{CoroExtremesDir}
    Let $P\in\PPP(\R^d)$ and $\rr\in\RRR$. Assume that $P$ is supported on a line or that $\rr$ is increasing. Further assume that $P$ has at most finitely many atoms or that $\rr(0)=0$. Let $(\alpha_k)\subset [0,1)$ and $(u_k)\subset\S^{d-1}$ be sequences such that $\alpha_k \to 1$ and $u_k\to u\in \S^{d-1}$ as $k\to\infty$. Let $Z$ denote a random $d$-vector with law $P$.
    \begin{enumerate}[label=(\roman*)]
    \item If $\ell_{\rr}(1)<\infty$ and $\E[\|Z\|]<\infty$, then
    \[
    \| \Qpr(\alpha_k u_k) \| \left( \frac{\Qpr(\alpha_k u_k)}{\| \Qpr(\alpha_k u_k) \|} - u_k \right) \to \E\big[Z-\ps{Z}{u} u\big]
    .
    \] 
    \item If $\ell_{\rr}(2)<\infty$ and $\E[\|Z\|^2]<\infty$, then 
    \[
    \frac{1}{\sqrt{1-\alpha_k}} \left( \frac{\Qpr(\alpha_k u_k)}{\| \Qpr(\alpha_k u_k) \|} - u_k \right) \to  \frac{\sqrt{2}\ \E\big[Z-\ps{Z}{u} u\big]}{2\ell_{\rr}(2) + \operatorname{Tr}(\Sigma) - \ps{u}{\Sigma u}}
    .
    \]
    \end{enumerate}
\end{Corol}

A few comments are in order for distributions with finite second moment, in light of Theorem~\ref{TheoExtremesDirNorm} and Corollary~\ref{CoroExtremesDir} above. (i) When $\rr(s)$ converges to 1 faster than $s^2$ as $s\to\infty$, then $\ell_{\rr}(2)=0$ and we recover the exact same asymptotics as geometric quantiles observed in Theorem~2.2 in~\cite{GirStu2017}. In particular, extreme $\rr$-quantile contours of elliptical distributions are orthogonal to their density level plots. (ii) If $\rr(s)$ converges to 1 precisely at rate $s^2$ as $s\to\infty$, then the asymptotic behavior of extreme $\rr$-quantiles does not reveal any further information about the extreme value behavior of the underlying distribution, also taking note that the extra term compared to Theorem~2.2 in~\cite{GirStu2017} depends exclusively on the choice of regularizer $\rr$. As the value of $\ell_{\rr}(2)$ increases, \ie as the regularizer $\rr$ puts more weight on extreme observations, the resulting $\rr$-quantile contours become increasingly spherical. (iii) If $\rr(s)$ converges to 1 at a rate $s^\beta$ slower than $s^2$, then the size of an extreme $\rr$-quantile $\Qpr(\alpha_k u_k)$ is asymptotically equivalent to $(\ell_{\rr}(\beta))^{-1/\beta} (1-\alpha_k)^{-1/\beta}$ and extreme $\rr$-quantile contours are spherical. These are universal asymptotics in the sense that they do not even depend on any characteristic of the underlying distribution. The asymptotic behavior of the direction, meanwhile, only depends on $\ell_{\rr}(\beta)$ and the first moment of the distribution.\smallskip

\begin{landscape}
%
%%%%%%%%%%%%%%%%%%%%%%%%%%%%%%%%%%
\begin{figure}[h!]
\centering
\includegraphics[width=0.95\linewidth]{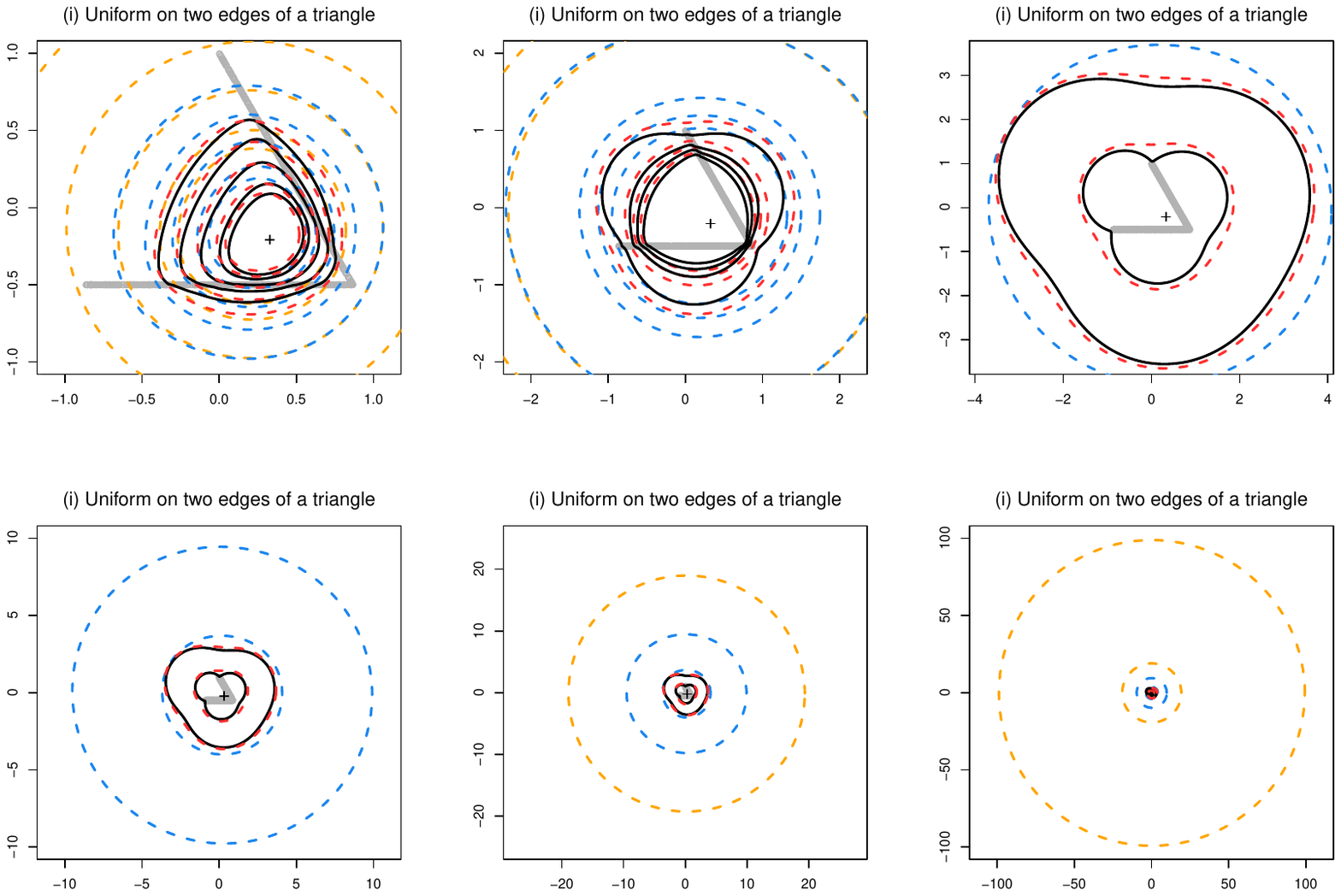}
\vspace{-1.5cm}
\caption{Graphical representation of the empirical quantiles $q_{\alpha,u}^\rr$ at six different levels of zoom, where $\rr=\rr_{\beta}$, for $\beta=1$ (dashed orange curves), $\beta=2$ (dashed blue curves) and $\beta=5$ (dashed red curves), compared with the classical geometric quantiles $q_{\alpha,u}^{1}$ (full black curves), at levels $\alpha\in \{0.25,0.3,0.4,0.5,0.6\}$ (top left panel), $\alpha\in \{0.7,0.75,0.8,0.9\}$ (top middle panel) and $\alpha\in \{0.95,0.99\}$ (remaining panels). The experiments use $n=1000$ data points generated from the uniform distribution on the union of the edges $AC$ and $BC$ of the triangle $ABC$ with vertices $A(0,1)$, $B(\sqrt{3}/2,-1/2)$ and $C(-\sqrt{3}/2,-1/2)$.}
\label{fig:continuous_1}
\end{figure}
%%%%%%%%%%%%%%%%%%%%%%%%%%%%%%%%%%
%
\end{landscape}

\begin{landscape}
%
%%%%%%%%%%%%%%%%%%%%%%%%%%%%%%%%%%
\begin{figure}[h!]
\centering
\includegraphics[width=0.95\linewidth]{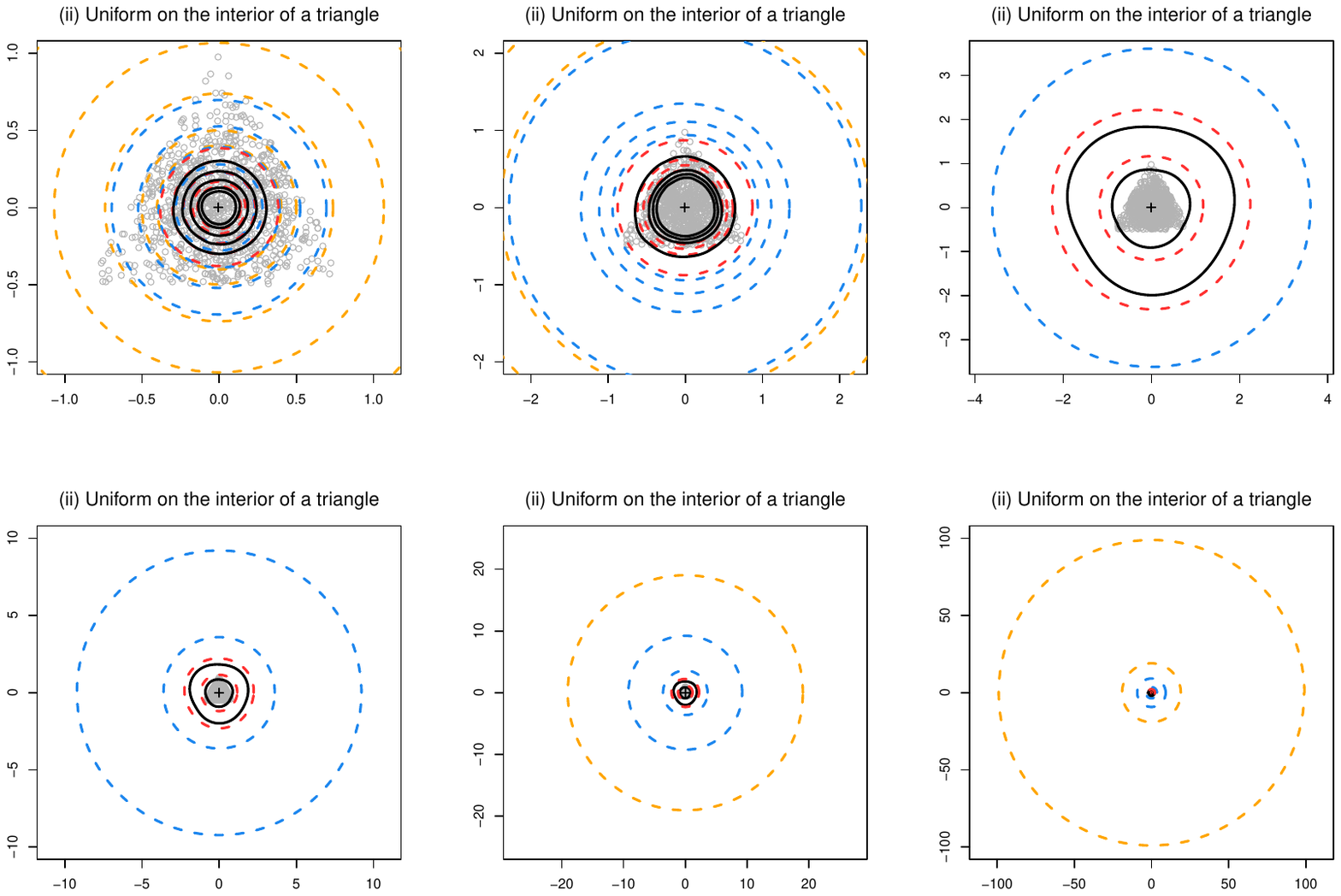}
\vspace{-1.5cm}
\caption{Graphical representation of the empirical quantiles $q_{\alpha,u}^\rr$ at six different levels of zoom, where $\rr=\rr_{\beta}$, for $\beta=1$ (dashed orange curves), $\beta=2$ (dashed blue curves) and $\beta=5$ (dashed red curves), compared with the classical geometric quantiles $q_{\alpha,u}^{1}$ (full black curves), at levels $\alpha\in \{0.25,0.3,0.4,0.5,0.6\}$ (top left panel), $\alpha\in \{0.7,0.75,0.8,0.9\}$ (top middle panel) and $\alpha\in \{0.95,0.99\}$ (remaining panels). The experiments use $n=1000$ data points generated from the uniform distribution on the interior of the triangle $ABC$ with vertices $A(0,1)$, $B(\sqrt{3}/2,-1/2)$ and $C(-\sqrt{3}/2,-1/2)$.}
\label{fig:continuous_2}
\end{figure}
%%%%%%%%%%%%%%%%%%%%%%%%%%%%%%%%%%
%
\end{landscape}

\begin{landscape}
%
%%%%%%%%%%%%%%%%%%%%%%%%%%%%%%%%%%
\begin{figure}[h!]
\centering
\includegraphics[width=0.95\linewidth]{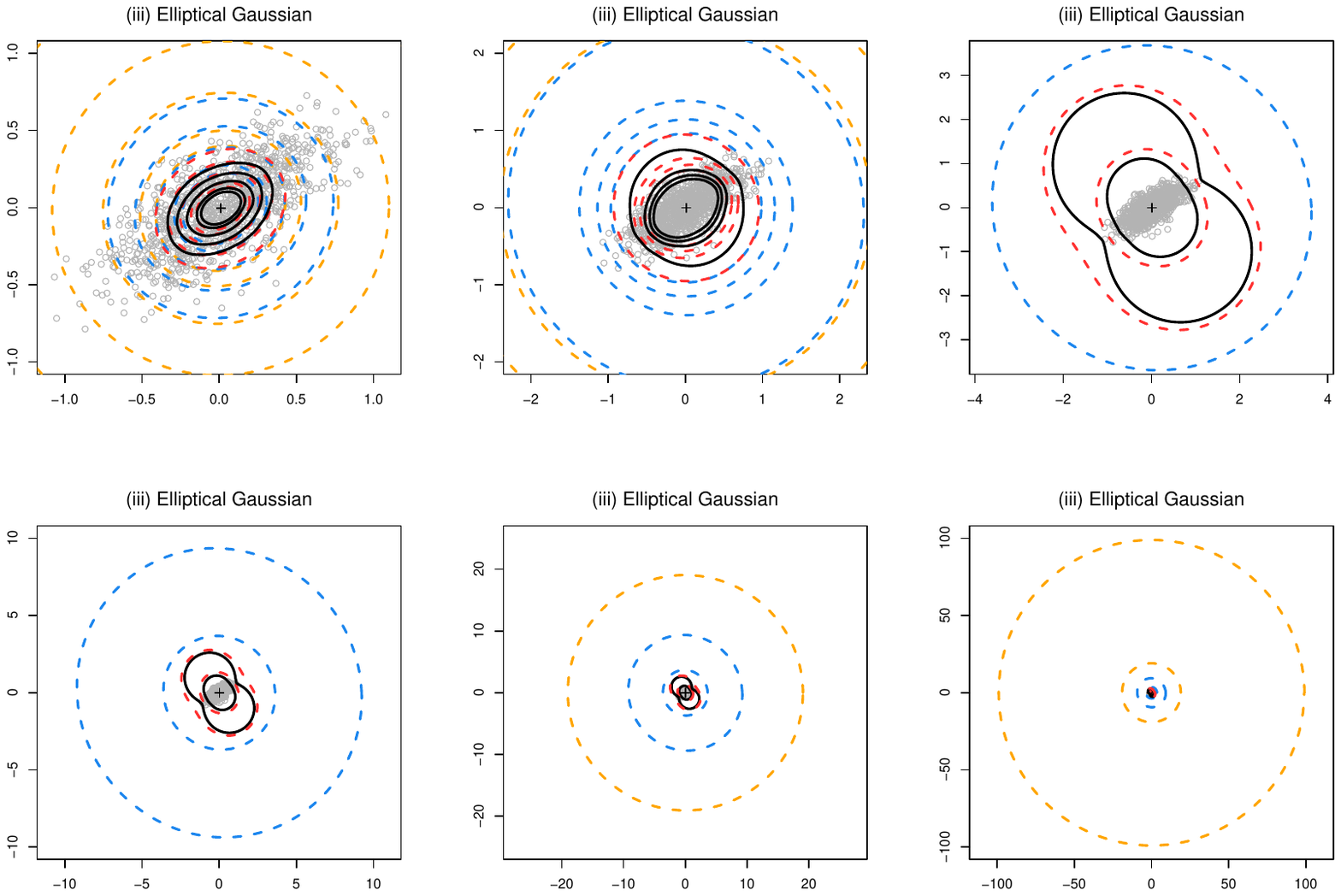}
\vspace{-1.5cm}
\caption{Graphical representation of the empirical quantiles $q_{\alpha,u}^\rr$ at six different levels of zoom, where $\rr=\rr_{\beta}$, for $\beta=1$ (dashed orange curves), $\beta=2$ (dashed blue curves) and $\beta=5$ (dashed red curves), compared with the classical geometric quantiles $q_{\alpha,u}^{1}$ (full black curves), at levels $\alpha\in \{0.25,0.3,0.4,0.5,0.6\}$ (top left panel), $\alpha\in \{0.7,0.75,0.8,0.9\}$ (top middle panel) and $\alpha\in \{0.95,0.99\}$ (remaining panels). The experiments use $n=1000$ data points generated from the Gaussian pair $(X,Y)$ with $\operatorname{Var}(X)=2/16$, $\operatorname{Var}(Y)=1/16$, $\operatorname{Cov}(X,Y)=1/16$.}
\label{fig:continuous_3}
\end{figure}
%%%%%%%%%%%%%%%%%%%%%%%%%%%%%%%%%%
%
\end{landscape}

\begin{landscape}
%
%%%%%%%%%%%%%%%%%%%%%%%%%%%%%%%%%%
\begin{figure}[h!]
\centering
\includegraphics[width=0.95\linewidth]{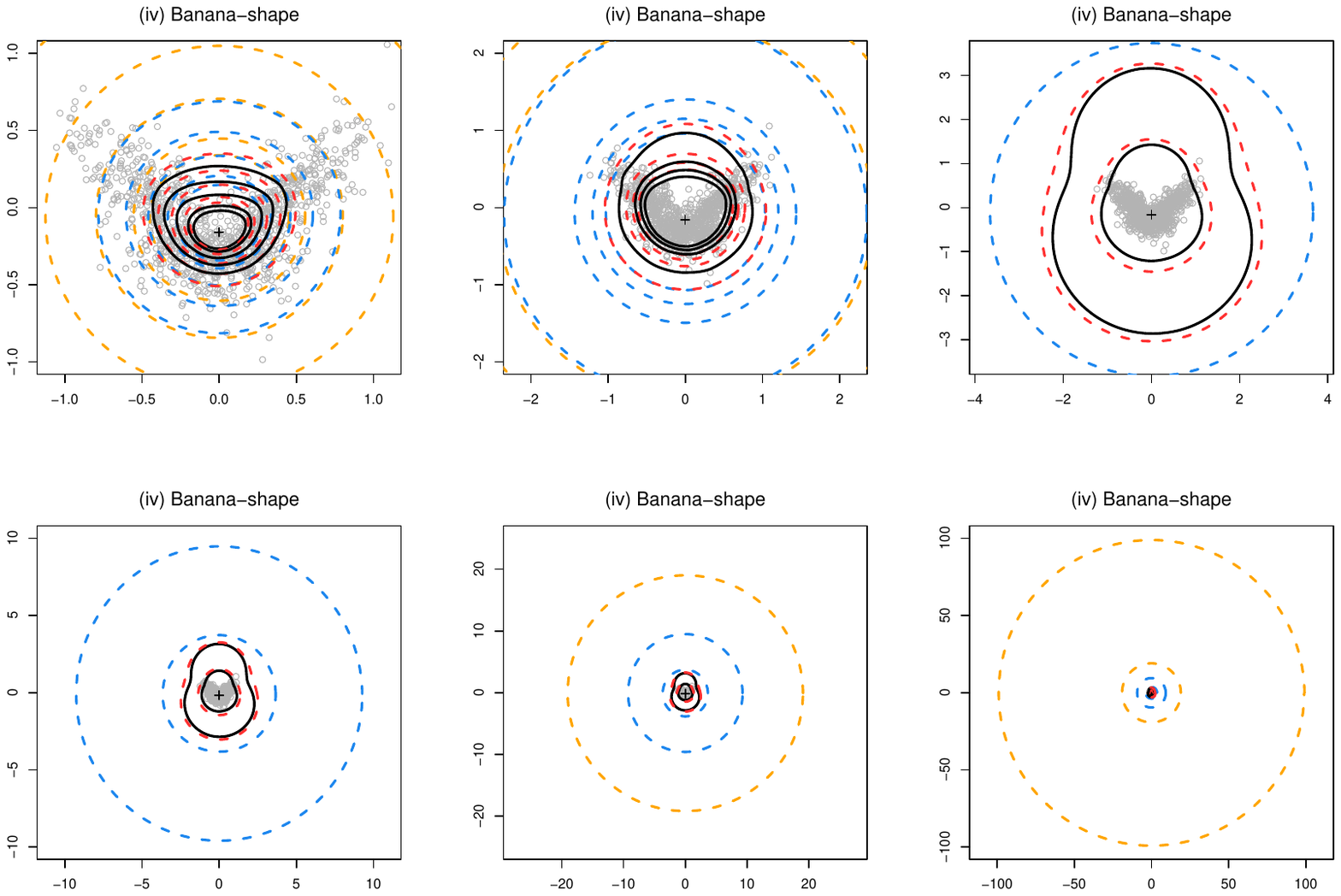}
\vspace{-1.5cm}
\caption{Graphical representation of the empirical quantiles $q_{\alpha,u}^\rr$ at six different levels of zoom, where $\rr=\rr_{\beta}$, for $\beta=1$ (dashed orange curves), $\beta=2$ (dashed blue curves) and $\beta=5$ (dashed red curves), compared with the classical geometric quantiles $q_{\alpha,u}^{1}$ (full black curves), at levels $\alpha\in \{0.25,0.3,0.4,0.5,0.6\}$ (top left panel), $\alpha\in \{0.7,0.75,0.8,0.9\}$ (top middle panel) and $\alpha\in \{0.95,0.99\}$ (remaining panels). The experiments use $n=1000$ data points generated from a mixture of three Gaussian random pairs with respective means $(-3/2,0)$, $(3/2,0)$ and $(0,-5/4)$, covariance matrices $\begin{pmatrix} 5/4 & -1 \\ -1 & 5/4 \end{pmatrix}$, $\begin{pmatrix} 5/4 & 1 \\ 1 & 5/4 \end{pmatrix}$ and $\begin{pmatrix} 1 & 0 \\ 0 & 1/4 \end{pmatrix}$, weighted with probabilities $3/8$, $3/8$ and $1/4$.}
\label{fig:continuous_4}
\end{figure}
%%%%%%%%%%%%%%%%%%%%%%%%%%%%%%%%%%
%
\end{landscape}

Figures~\ref{fig:continuous_1} to~\ref{fig:continuous_4} illustrate these phenomena, using $n=1000$ data points sampled from four different bivariate distributions, and comparing geometric quantiles with $\rr_{\beta}$-geometric quantiles for $\beta=1,2,5$. For $\beta=5$, while there is a visible difference for low values of $\alpha$, extreme $\rr_{\beta}$-geometric quantile contours are almost identical to the corresponding extreme classical geometric quantile contours. These contours are, as expected, in some sense orthogonal to the shape of the data cloud. For $\beta=2$, extreme $\rr_{\beta}-$geometric quantile contours appear to be spherical; this is due to the fact that here $\ell_{\rr}(2)=1$, which is large compared to the quantity $\frac{1}{2} (\operatorname{Tr}(\Sigma) - \ps{u}{\Sigma u})$ appearing in the asymptotics of Theorem~\ref{TheoExtremesDirNorm}(ii) regardless of the direction $u$, so that the influence of $u$ on the contours is negligible. For $\beta=1$, extreme $\rr_{\beta}$-geometric quantile contours are spherical, as predicted by Theorem~\ref{TheoExtremesDirNorm}(i). Note, moreover, that for $\beta=2$ (resp.~$\beta=1$), extreme $\rr_{\beta}$-geometric quantile contours at level $\alpha=0.99$ have norm roughly equal to 10 (resp.~100) in every direction, as suggested by Theorem~\ref{TheoExtremesDirNorm}.\smallskip

To summarize, the first-order asymptotics of extreme $\rr$-quantiles of distributions with a finite second moment do not shed any light on tail characteristics, even when the regularizer $\rr$ converges to $1$ slowly (\ie~when the $\rr$-distribution function $\Fpr$ puts substantially larger weight towards the extremes). This further implies that outliers detection procedures based on extreme $\rr$-quantiles (including geometric quantiles), such as~\cite{ChaGog2010}, should be conducted with great care.

\section{Empirical $\rr$-quantiles}\label{sec:Empiric}

% \subsection{Pointwise behaviour of empirical quantiles}\label{sec:EmpiricPoint}

We now consider estimation of the $\rr$-quantiles $\Qpr(\alpha u)$ based on a random sample~$Z_1,\ldots,Z_n$ drawn from~$P$. The canonical estimator is obtained by replacing~$P$ with the corresponding empirical distribution. In this section, we study the asymptotic properties of the resulting sample $\rr$-quantiles. To deal with cases where $\rr$-quantiles may be non-unique, we introduce the set of $\rr$-quantiles of a given order $\alpha\in [0,1)$ in direction $u\in \S^{d-1}$
$$
\MM_{\alpha,u}^P
=
\big\{
x \in\R^d : M_{\alpha,u}^{\rr, P}(x) = \inf_{y\in\R^d} M_{\alpha,u}^{\rr, P}(y)
\big\}
.
$$
Then, for any $x\in\R^d$ define the distance between $x$ and $\MM_{\alpha,u}^P$ as the smallest Euclidean distance between $x$ and any $m\in\MM_{\alpha,u}^P$, \ie 
$$
d(x,\MM_{\alpha,u}^P) 
\equiv 
\inf_{m\in \MM_{\alpha,u}^P} \|x-m\| 
.
$$
We start with the following consistency result. 

% \begin{Prop}\label{PropContMEmpiri}
%     Let $P\in\PPP(\R^d)$ and $\rr\in\RRR$. Fix $\alpha\in [0,1)$ and $u\in\S^{d-1}$. Then, the map $M_{\alpha,u}^{\rr,P}$, introduced in Definition \ref{DefinQuantiles}, is $1$-Lipschitz continuous over $\R^d$, \ie we have for all $x,x'\in\R^d$
%     $$
%     |M_{\alpha,u}^{\rr,P}(x)-M_{\alpha,u}^{\rr,P}(x')|
%     \le
%     \|x-x'\|
%     .
%     $$
% \end{Prop}

% \vspace{1mm}

% \begin{Proof}{Proposition \ref{PropContMEmpiri}}
%     A direct computation provides 
%     $$
%     M_{\alpha,u}^{\rr,P}(x)-M_{\alpha,u}^{\rr,P}(x')
%     =
%     \int_{\R^d} \big( \rR(z-x) - \rR(z-x') \big)\, dP(z) 
%     .
%     $$
%     The conclusion follows from the fact that $\rR$ is Lipschitz continuous; \com{see (...)}.
% \end{Proof}

\begin{Theor}\label{TheorConsistency}
    Let $P\in\PPP(\R^d)$ and $\rr\in\RRR$. Fix $\alpha\in [0,1)$ and $u\in\S^{d-1}$. Let $Z_1, Z_2\ldots$ be a random sample drawn from $P$. Denote by $P_n\equiv  n^{-1}\sum_{i=1}^n \delta_{Z_i}$ the induced empirical distribution, and $M_{\alpha,u}^{\rr,P_n}$ the corresponding objective function. For $\star={\rm a.s.}$ or $\star={\rm P}$, denote by $o_\star(1)$ any random variable that converges to $0$ almost surely or in $P$-probability, respectively, as $n\to\infty$. Let $(\hat q_n)$ be a sequence of estimators satisfying, as $n\to\infty$,
    $$
    M_{\alpha,u}^{\rr,P_n}(\hat q_n) 
    \le 
    \inf_{x\in\R^d} M_{\alpha,u}^{\rr,P_n}(x) + o_\star(1)
    .
    $$
    Then, 
    % an arbitrary $\rr$-quantile of order $\alpha$ in direction $u$ for the empirical distribution $P_n\equiv  n^{-1}\sum_{i=1}^n \delta_{Z_i}$. Then, 
    $
    d(\hat q_n, \MM_{\alpha,u}^P) 
    =
    o_\star(1)
    $. In particular, if $\MM_{\alpha,u}^P=\{q_{\alpha,u}^\rr\}$ is unique, then $\hat q_n - q_{\alpha,u}^\rr = o_\star(1)$.
\end{Theor}
\vspace{1mm}

% Comment on the case where $\hat q_n$ is exactly the empirical quantile.\\

Observe that if $\rr$ is strictly increasing over $(0,\infty)$, then $\hat q_n$ is automatically unique by virtue of Theorem \ref{TheorUniqueness}. When this is not necessarily the case (for geometric quantiles obtained with $\rr\equiv 1$, for instance) and when $P$ is not supported on a single line of $\R^d$, then there exists $p\in [0,1)$ such that, for all $n\ge 3$, the observations $Z_1,\ldots, Z_n$ are collinear with $P$-probability at most $p^{n-2}$. In particular, $\hat q_n$ is unique with probability at least $1-p^{n-2}$. Indeed, observe that the mutual independence of $Z_1,\ldots, Z_n$ entails that
\begin{eqnarray*}
\lefteqn{
\P[Z_1,\ldots, Z_n \text{ are collinear} ]
% }
% \\[2mm]
% && 
=
\P\big[Z_3,\ldots, Z_n \in \{Z_1 + \lambda (Z_2-Z_1) : \lambda\in\R\}\big]
}
\\[2mm]
&&=
\int_{\R^d\times \R^d} \P\big(Z_3 \in \{z_1 + \lambda (z_2-z_1) : \lambda\in\R\}\big)^{n-2}\, dP(z_1)\, dP(z_2)
% \\[2mm]
% &&
\le 
p^{n-2}
,
\end{eqnarray*}
where, letting the supremum range over all lines of $\R^d$,
$
p
\equiv 
\sup_{\LL} P[\LL]
<
1
$
%The fact that $p\leq 1$ is obvious, and we show in 
(see Lemma \ref{LemPLine}). %that $p<1$ when $P$ is not supported on a single line of $\R^d$. 
Instead of a high-probability bound valid for all $n\ge 1$, %we can alternatively make use of Borel-Cantelli's lemma to establish that, with $P$-probability $1$, the 
the upper bound in the last display combined with Borel-Cantelli's lemma implies that $Z_1,\ldots, Z_n$ are not collinear for all $n$ large enough (sample-path dependent). In particular, $\hat q_n$ is eventually unique with $P$-probability $1$. A similar result was obtained in \cite{Romon2022} in the weaker sense of $\P^*$-outer probability and was established by means of different proof techniques (see Section 3.4 there).\smallskip 

Assuming further moment conditions on $P$ or additional %regularization 
regularity properties of $\rr$, we establish in the following result that all $\rr$-quantiles enjoy a standard asymptotic normality result, relying on a neat Bahadur representation. For this purpose, we require that $P$ be not supported on a line, or that $\rr$ is strictly increasing, so that the $\rr$-quantile is unique by virtue of Theorem \ref{TheorUniqueness}; in fact, we need the slightly stronger assumption $\rr'(s)>0$ for all $s>0$, to further ensure that the Hessian matrix featuring in the distributional limit is invertible.

\begin{Theor}\label{TheorNormal}
    Let $P\in\PPP(\R^d)$. Fix $\rr\in\RRR$ (right-)differentiable at $0$ and Lipschitz continuous over $[0,\infty)$. Fix $\alpha\in [0,1)$ and $u\in\S^{d-1}$. Assume that (i) $\rr'(s)>0$ for all $s>0$ and $P$ is not supported on a single point of $\R^d$, or (ii) $P$ is not supported on a line. %; in particular, $P$ admits a unique $\rr$-quantile $q_{\alpha,u}^\rr$ (Theorem \ref{TheorUniqueness}). 
    Further assume that (iii) $\rr(0)=\rr'(0)=0$, %$\rr(0)=0$ and $\rr'(0)=0$, 
    or (iv) $q_{\alpha,u}^\rr$ is not an atom of $P$ and $P$ satisfies the moment condition
    $$
    \int_{\R^d\sm\{q_{\alpha,u}^{\rr}\}} \frac{1}{\|z-q_{\alpha,u}^{\rr}\|}\, dP(z)
    <
    \infty 
    .
    $$
    % Further assume that one of the following holds: (iv) $\rr$ is Lipschitz over $(0,\infty)$, or (v) there exists $\delta>0$ such that $t\mapsto \rr(t)/t$ is non-increasing over $(0,\delta)$, and $\rr$ is Lipschitz over $[\delta, \infty)$. 
    Let $Z_1, Z_2\ldots$ be a sample drawn from $P$. For any $n\geq 1$, let $\hat q_n$ be an arbitrary $\rr$-quantile of order $\alpha$ in direction $u$ for the empirical distribution $P_n\equiv n^{-1}\sum_{i=1}^n \delta_{Z_i}$. Then, 
% \begin{eqnarray*}
% \lefteqn{
$$
\sqrt{n}(\hat q_n-q_{\alpha,u}^\rr)
% }
% \\[2mm]
% &&\hspace{-5mm}
=
A^{-1}
\sqrt{n}\ F_{P_n}^{\rr}(q_{\alpha,u}^\rr)
% \big( \nabla \rR(q_{\alpha,u}^\rr-Z_i)  - \alpha u \big)\I[Z_i\neq q_{\alpha,u}^\rr]
% \big( \nabla \rR(q_{\alpha,u}^\rr-Z_i)  - \alpha u \big)\I[Z_i\neq q_{\alpha,u}^\rr]
% \Big( \rr(\|q_{\alpha,u}^\rr-Z_i\|)\frac{q_{\alpha,u}^\rr-Z_i}{\|q_{\alpha,u}^\rr-Z_i\|}\I[Z_i\neq q_{\alpha,u}^\rr] - \alpha u \Big)
+
o_{\rm P}(1)
\stackrel{d}{\to}
\mathcal{N}(0,\Sigma)
,
\spa 
n\to \infty 
,
% \end{eqnarray*}
$$
where~$\Sigma\equiv  A^{-1}BA^{-1}$ involves the positive definite matrices
$$
A\equiv \E\Big[ \nabla^2 \rR(Z_1-q_{\alpha,u}^\rr) \I[Z_1\neq q_{\alpha,u}^\rr] \Big]
,
$$
% $$
% A
% \equiv 
% \E\bigg[ \bigg( \rr'(\|q_{\alpha,u}^\rr-Z_1\|) \frac{(q_{\alpha,u}^\rr-Z_1) (q_{\alpha,u}^\rr-Z_1)^T}{\|q_{\alpha,u}^\rr-Z_1\|^2} + \frac{\rr(\|q_{\alpha,u}^\rr-Z_1\|)}{\|q_{\alpha,u}^\rr-Z_1\|}\Big({\rm I}_d - \frac{(q_{\alpha,u}^\rr-Z_1)(q_{\alpha,u}^\rr-Z_1)^T}{\|q_{\alpha,u}^\rr-Z_1\|^2}\Big) \bigg) \I[Z_1\neq 0]\bigg]
% % \nabla^2 M_{\alpha, u}^{\rho}(\mu^{\rho}_{\alpha,u}) $ 
% $$
and
$$
B
\equiv 
\E\Big[
\big( \nabla \rR(q_{\alpha,u}^\rr-Z_1)  - \alpha u \big)
\big( \nabla \rR(q_{\alpha,u}^\rr-Z_1)  - \alpha u \big)^T \I[Z_1\neq q_{\alpha,u}^\rr]
% (\nabla H^{\rho}_{\alpha,u}(Z_1-\mu^{\rho}_{\alpha,u}))
% (\nabla H^{\rho}_{\alpha,u}(Z_1-\mu^{\rho}_{\alpha,u}))^T
% \I[\|Z_1-\mu^{\rho}_{\alpha,u}\|\in \mathcal{D}_\rho]
\Big]
.
$$ 
\end{Theor}
\vspace{1mm}

Assumptions (ii) and (iv) in Theorem \ref{TheorNormal} are classical in the literature on asymptotic results for geometric quantiles. The novelty of Theorem \ref{TheorNormal} is that assumptions (ii) and (iv) can be replaced by (i) and (iii) which, instead, hold for \emph{any} (non-Dirac) probability measure $P$, provided~$\rr$ is regularizing enough for the singular kernel $x\mapsto x/\|x\|$ in the sense that $\rr(0)=\rr'(0)=0$; in this case, no non-atomicity or moment assumption is required, and asymptotic normality automatically holds for our regularized geometric quantiles.

\section*{Acknowledgements}
Dimitri Konen acknowledges funding from an ERC Advanced Grant (UKRI G116786). Gilles Stupfler acknowledges financial support from the French National Research Agency under the grants ANR-23-CE40-0009 (EXSTA project) and ANR-11-LABX-0020-01 (Centre Henri Lebesgue), as well as from the TSE-HEC ACPR Chair `Regulation and systemic risks' and the Chair Stress Test, RISK Management and Financial Steering of the Foundation Ecole Polytechnique.

\bibliographystyle{agsm}
\bibliography{mybib}

%%%%%%%%%%%%%%%%%%%%%%%%%%%%%%%%%%
%%%%%%%%%%%%%%%%%%%%%%%%%%%%%%%%%%

\appendix

\section{Supplementary material -- Proofs}

\subsection{Proof of Theorem \ref{TheorUniqueCdf}}\label{sec:ProofsUniver}

   The proof proceeds in four steps. In step 1, we show that $\Fd$ extends to all of ${\rm span}(\GG)$. In step 2 we show $\Fd$ gives rise a bounded linear form on $(L^1(\R^d))^d$, which allows, in Step 3, to deduce that $\Fd$ is given as a convolution with respect to a kernel $K$. Finally, we conclude in Step 4 that $K$ has a specific form imposed by the the orthogonal and translation equivariance constraints.\\
   
{\it Step 1.} Let us show that $\Fd$ extends to a linear map $\Fd:{\rm span}(\PP)\to{\rm Maps}(\R^d, \R^d)$ defined on all of ${\rm span}(\PP)$ that (necessarily) takes its values in ${\rm Maps}(\R^d,\R^d)$ instead of ${\rm Maps}(\R^d,\overline{\B}^d)$. We first extend $\Fd$ to $\PP^+\equiv \{ d\mu = g\,dx : g\in\GG^+\}$ where 
$$
\GG^+\equiv \{sf : s\geq 0,\ f\in\GG\}
.
$$
In what follows we will abusively write $\Fd(g)$ instead of $\Fd(\mu)$ when $d\mu=g\,dx$. Since any $g\in \GG^+$ can be written \emph{in a unique way} as $g=s f$ with $s\geq 0$ and $f\in\GG$, let us define for all $s \geq 0$ and $f\in\GG$
\begin{equation}\label{eq:ConvexLinear}
\Fd(s f)\equiv  s\, \Fd(f)
.
\end{equation}
% It is easy to see that $\Fd$ is well-defined on the set $C^\infty_c(\R^d, [0,\infty))$ of nonnegative functions in $C^\infty_c(\R^d)$---in the sense that if $s_1,s_2\geq 0$ and $g_1,g_2\in\EEE$ are such that $s_1g_1 = s_2g_2$ then in fact $s_1=s_2$ and $g_1=g_2$ so that $\Fd(s_1g_1)=\Fd(s_2 g_2)$---and that 
% In particular, $\Fd$ coincides with $\Fd$ on $\PP$ so that $\Fd$ genuinely extends $\Fd$. 
Fix $g_1,g_2\in\GG^+$ and let us show that 
\begin{equation}\label{eq:LinearExt}
    \Fd(g_1 + g_2)
    =
    \Fd(g_1)+\Fd(g_2)
    .
\end{equation} 
If $g_1\equiv  0$ or $g_2\equiv  0$, then (\ref{eq:LinearExt}) trivially holds. Then assume that $g_1$ and $g_2$ do not identically vanish over $\R^d$. For $i\in\{1,2\}$, define $\lambda_i \equiv  \int_{\R^d} g_i(x)\,dx$ and notice that $\lambda_i > 0$. Then, write
$$
g_1 + g_2 
=
(\lambda_1 + \lambda_2) 
\bigg\{ 
\frac{\lambda_1}{\lambda_1 + \lambda_2} \frac{g_1}{\lambda_1}
+
\frac{\lambda_2}{\lambda_1 + \lambda_2} \frac{g_2}{\lambda_2}
\bigg\}
    =:
    (\lambda_1+\lambda_2) g_{1,2}
.
$$
    Since $g_{1,2}=(g_1+g_2)/(\lambda_1+\lambda_2)\in\GG$, it follows from (\ref{eq:ConvexLinear}) that
    $$
    \Fd(g_1 + g_2)
    =
    (\lambda_1 + \lambda_2)\
    \Fd \bigg( \frac{\lambda_1}{\lambda_1 + \lambda_2} \frac{g_1}{\lambda_1}
    +
    \frac{\lambda_2}{\lambda_1 + \lambda_2} \frac{g_2}{\lambda_2} \bigg)
    .
    $$
    Now observe that $g_{1,2}$ is a convex combination of $g_1/\lambda_1\in\GG$ and $g_2/\lambda_2\in\GG$. In particular, the 'linearity' assumption of $\Fd$ (see Assumption \ref{Asmp:Linear}) yields
\begin{align*}
    \Fd(g_1 + g_2)
    =
    \lambda_1\, \Fd\Big(\frac{g_1}{\lambda_1}\Big)
    +
    \lambda_2\, \Fd\Big(\frac{g_2}{\lambda_2}\Big)
    % \\[2mm]
    =
    \Fd(g_1) + \Fd(g_2)
    ,
\end{align*}
    where the last equality follows from (\ref{eq:ConvexLinear}). Let us now extend $\Fd$ to ${\rm span}(\PP)$. Since any $h\in{\rm span}(\GG)$ writes as $h=g_2-g_1$ for some $g_1,g_2\in \GG^+$, then define for any $g_1,g_2\in \GG^+$
$$
\Fd(g_2-g_1) \equiv  \Fd(g_2) - \Fd(g_1)
.
$$
Let us show that $\Fd$ is well-defined on ${\rm span}(\PP)$, \ie if $g_1,g_2, g_1^*, g_2^*\in\GG^+$ satisfy $g_2-g_1=g_2^*-g_1^*$ then $\Fd(g_2^*)-\Fd(g_1^*)=\Fd(g_2)-\Fd(g_1)$. This follows by observing that, since $g_2 + g_1^* = g_1 + g_2^*$, then (\ref{eq:LinearExt}) yields 
$
\Fd(g_2) + \Fd(g_1^*)
=
\Fd(g_1)+\Fd(g_2^*)
.
$
% In particular, we have 
% $$
% \Fd(g_2) - \Fd(g_1)
% =
% \Fd(\psi_2)- \Fd(\psi_1)
% .
% $$
It follows that $\Fd$ is well-defined on ${\rm span}(\PP)$.

Let us now show that $\Fd$ is linear on ${\rm span}(\PP)$. Fix $\lambda\in\R$ and $h\in{\rm span}(\GG)$ with $h=g_2-g_1$ for some $g_1,g_2\in \GG^+$. 
  % If $\lambda\geq 0$, then $\lambda\, h = \lambda g_2-\lambda g_1$ with $\lambda g_i\in \PP^+$, whereas if $\lambda\leq 0$ then $\lambda\, h = |\lambda|g_1 - |\lambda|g_2$. It follows that $\lambda\, h = {\rm sign(\lambda)}(|\lambda|g_2-|\lambda|g_1)$. 
  If $\lambda=0$, then it is clear that $\Fd(\lambda\, h)=\lambda\, \Fd(h)$. If $\lambda >0$, then
$$
\Fd(\lambda\, h)
=
\Fd(\lambda g_2)-\Fd(\lambda g_1)
=
\lambda\, \Fd(g_2) - \lambda\, \Fd(g_1)
=
\lambda\, \Fd(h)
.
$$
If $\lambda < 0$, then $\lambda\, h = (-\lambda) (g_1-g_2)$, with $-\lambda >0$ so that 
$$
\Fd(\lambda\, h)
=
(-\lambda) \big(\Fd (g_1) - \Fd(g_2) \big)
=
\lambda\, \Fd(h)
.
$$
Now take $h^*\in {\rm span}(\GG)$ and $g_1^*,g_2^*\in\GG^+$ such that $h^*=g_2^*-g_1^*$. Then,
$$
h+h^*
=
(g_2 + g_2^*)-(g_1+g_1^*)
.
$$ 
Consequently, (\ref{eq:LinearExt}) entails that 
$$
\Fd(h+h^*)
=
\Fd(g_2+g_2^*) - \Fd(g_1 + g_1^*)
=
\big(\Fd(g_2)-\Fd(g_1)\big)
+
\big(\Fd(g_2^*)-\Fd(g_1^*)\big)
,
$$
which yields $\Fd(h + h^*)=\Fd(h)+\Fd(h^*)$. We conclude that $\Fd$ is linear on ${\rm span}(\PP)$.\\

{\it Step 2.} Let us now show that
\begin{equation}\label{eq:BoundedF}
\|(\Fd h)(0)\| 
\leq 
\int_{\R^d} |h(x)|\, dx 
\end{equation}
for all $h\in {\rm span}(\GG)$. First observe that, since this holds for any $h\in \GG$ by definition, the linearity of $\Fd$ entails that this also holds for any $h\in\GG^+$. Now fix $h\in {\rm span}(\GG)\subset \FF$, and let $h^+\equiv |h|\, \I[h\geq 0]$ and $h^- = |h|\, \I[h\leq 0]$. In particular, we have $h^+ = (|h|+h)/2$ and $h^-=(|h|-h)/2$. By Assumption \ref{Asmp:ClassP} we have $|h|\in\FF$, so that $h^+,h^-\in\FF$. Since $h^+\geq 0$ and $h^-\geq 0$ we have $h^+,h^-\in\GG^+$ by definition of $\GG$, so that 
$$
\|(\Fd h)(0)\| 
\leq 
\|(\Fd h^+)(0)\| 
+
\|(\Fd h^-)(0)\| 
\leq 
\int_{\R^d} h^+(x)\, dx 
+
\int_{\R^d} h^-(x)\, dx 
=
\int_{\R^d} |h(x)|\, dx 
.
$$

{\it Step 3.} From the observations made in Step 2 of this proof, we easily deduce that ${\rm span}(\GG)=\FF$. In particular, ${\rm span}(\GG)$ is dense in $L^1(\R^d)$, so that Steps 1 and 2 above entail that the linear map $\Lambda:{\rm span}(\GG)\to \R^d$, defined by $\Lambda(h)\equiv (\Fd h)(0)$ for all $h\in{\rm span}(\GG)$, extends to a bounded linear operator on $L^1(\R^d)$, which we still denote $\Lambda$. Reasoning componentwise, standard results of functional analysis---see, e.g., Theorem 6.5.7(ii) in \cite{BogRealAndFunctional}---entail that there exists $J\in (L^\infty(\R^d))^d$ such that 
$$
\Lambda(f)
=
\int_{\R^d} J(z) f(z)\, dz 
$$
for all $f\in L^1(\R^d)$. For any $h\in {\rm span}(\GG)$, the translation equivariance of $\Fd$ (see Assumption \ref{Asmp:EquivTrans} and, in particular, the version given in (\ref{eq:EquivCdfDensityTrans})) yields 
\begin{align*}
(\Fd h)(x)
&=
\big( (\Fd h) \circ T_{-x}\big)(0)
=
\big(\Fd (h \circ T_{-x})\big)(0)
    \\[2mm]
&\hspace{10mm} 
    =
    \Lambda( h\circ T_{-x} )
    =\int_{\R^d} J(z)\, (h\circ T_{-x})(z)\, dz
    \\[2mm]
    &\hspace{10mm}
=\int_{\R^d} J(z)\, h(z+x)\, dz
    =\int_{\R^d} J(-z)\, h(x-z)\, dz
    \\[2mm]
    &\hspace{10mm}
    =(K*h)(x)
,
\end{align*}
where we let $K(z)\equiv J(-z)$ for all $z\in\R^d$. Because $K\in (L^\infty(\R^d))^d$, we can extend $\Fd$ to all of $\PPP(\R^d)$ by letting 
$$
(\Fd P)(x)
\equiv 
\int_{\R^d} K(x-z)\, dP(z)
$$
for any $P\in\PPP(\R^d)$.\\

{\it Step 4.} Denote by $\OO(d)$ the group of orthogonal $d\times d$ matrices, and fix $U\in\OO(d)$. Then, the orthogonal equivariance of $\Fd$ (see Assumption \ref{Asmp:EquivOrtho} and, in particular, the version given in (\ref{eq:EquivCdfDensityOrtho})) entails that 
\begin{eqnarray*}
\lefteqn{
    \int_{\R^d} K(Ux-z)\, h(z)\, dz
    =
    (\Fd h)(Ux)
    =
    \big( (\Fd h)\circ T_U \big)(x)
}
\\[2mm]
&&\hspace{10mm}
=
\big( T_U\circ \Fd(h\circ T_U) \big)(x)
=
\int_{\R^d} U K(x-z)\, (h\circ T_U)(z)\, dz
\\[2mm]
&&\hspace{10mm}
=
\int_{\R^d} U K(x-z)\, h(Uz)\, dz
=
\int_{\R^d} U K(x-U^T z)\, h(z)\, dz
\end{eqnarray*}
for all $h\in {\rm span}(\GG)$. Since $K\in (L^\infty(\R^d))^d$ and ${\rm span}(\GG)$ is dense in $L^1(\R^d)$, we have 
$$
K(Ux-z)
=
U K(x-U'z)
$$
for all $x,z\in\R^d$ and any $U\in \OO(d)$. Since $K(Ux-z)=K(U(x-U^T z))$, we deduce that 
\begin{equation}\label{eq:OrthoK}
K(Vy)=V K(y),
\end{equation} 
for all $y\in\R^d$ and $V\in\OO(d)$. 
%         This entails that there exists a bounded map $R:\R^d $
%         % $$
%         % \int_{\R^d} K(Ux-z)\, f(z)\, dz
%         % =
%         % \int_{\R^d} K(x-U'z)\, f(z)\, dz
%         % $$
%         % for all $f\in L^1(\R^d)$ so that, by a density argument, we have
%         % $$
%         % \int_{\R^d} K(Ux-z)\, dP(z)
%         % =
%         % \int_{\R^d} K(x-U'z)\, dP(z)
%         % $$
%         % for any $P\in \PPP(\R^d)$.
% 		$$
% 		O(\Fd \varphi)(x)
% %		=
% %		(\RR_O \Fd \varphi)(x)
% 		=
% 		(\Fd \RR_O \varphi)(O x)
% 		=
% 		\int_{\R^d} K(Ox-z) \varphi(O'z)\, dz
% 		=
% 		\int_{\R^d} K(O(x-z))\varphi(z)\, dz
% 		$$
% 		for any $\varphi\in C^\infty_c(\R^d)$. A standard density argument yields 
% 		$$
% 		O(\Fd P)(x)
% 		=
% 		\int_{\R^d} K(O(x-z))\, dP(z)
% 		$$
% 		for any $P\in\PPP(\R^d)$. Taking $P$ as the dirac measure $\delta$ at $0$ then yieds 
% 		\begin{equation}\label{eq:OrthoK}
% 		O K(x)
% 		=
% 		K (Ox)
% 		\end{equation}
% 		for all $x\in\R^d$ and $d\times d$ orthogonal matrix $O$. 
    This readily yields $K(0)=0$. Now fix $\lambda>0$ and denote by $e_1$ the first vector of the canonical basis of $\R^d$. Letting $U=\diag(1,-1,\ldots, -1)$ be the diagonal (orthogonal) matrix with eigenvalues $1,-1,\ldots, -1$ and writing $K(\lambda e_1)=:(z_1,\ldots, z_d)$, then (\ref{eq:OrthoK}) yields $z_i = 0$ for $i=2,3,\ldots, d$. It follows that for all $\lambda > 0$ we have
    $$
    K(\lambda e_1) 
    =
    s_\lambda \|K(\lambda e_1)\|e_1
    $$
    for some $s_\lambda\in\{-1,1\}$. Since any $u\in \S^{d-1}$ writes $U e_1$ for some $U\in\OO(d)$, we have for all $\lambda > 0$
    $$
    K(\lambda u)
    =
    U K(\lambda e_1)
    =
    s_\lambda \|K(\lambda e_1)\| U e_1
    =
    s_\lambda \| K(\lambda e_1)\| u 
    .
    $$
    Because any $x\in\R^d\setminus\{0\}$ writes $x=\lambda u$ for some $\lambda>0$ and $u\in \S^{d-1}$, with $\lambda = \|x\|$ and $u=x/\|x\|$, there exists a bounded map $\rr:(0,\infty) \to \R$ such that
    $$
    K(x)
    =
    \rr(\|x\|) \frac{x}{\|x\|} 
    $$
    for all $x\in\R^d\setminus\{0\}$. Since $K(0)=0$, we finally have 
    $$
    K(x)
    =
    \rr(x)\frac{x}{\|x\|} \1[x\neq 0]
    .
    $$
    Recalling that $K(x)=(\Fd \delta)(x)$ for all $x\in\R^d$, with $\delta$ the Dirac probability measure at $0$, and because the map $(\Fd \delta)$ is bounded by $1$, we have $\|K(x)\|\leq 1$ for all $x\in\R^d$. Since $\rr(\lambda e_1) = s_\lambda |K(\lambda e_1)\|$ for all $\lambda>0$, then $\rr$ takes its values in $[-1,1]$, which concludes the proof.

\subsection{Proofs for Section \ref{sec:ExistUniq}}

\begin{Proof}{Theorem \ref{TheorExist}}
Part (i) follows by a routine application of the dominated convergence theorem, after observing that $|\rR(x-z)-\rR(z)|\leq \|x\|$ for all $x,z\in\R^d$ by the mean value theorem and because $|\rr|\leq 1$. (ii) Using the monotonicity of $\rr$ it is straightforward to see that $\rR(y)/\|y\|\to 1$ as $\|y\|\to\infty$. In particular, the dominated convergence theorem (which applies since $|\rR(x-z)-\rR(z)|\leq \|x\|$ for all $x,z\in\R^d$) entails that 
\begin{equation}\label{eq:AsymMalphau}
    \frac{M_{\alpha,u}^{\rr,P}(x)}{\|x\|} 
    + 
    \alpha \frac{\ps{u}{x}}{\|x\|}
    \to 1
\end{equation}
as $\|x\|\to\infty$ since $\rr(t)\uparrow 1$ as $t\uparrow \infty$. In particular, $\liminf_{\|x\|\to\infty} M_{\alpha,u}^{\rr,P}(x)/\|x\|\geq 1-\alpha > 0$. Together with the continuity of $M_{\alpha,u}^{\rr,P}$, this entails that $M_{\alpha,u}^{\rr,P}$ admits a global minimum over $\R^d$, which concludes the proof.
\end{Proof}
\vspace{3mm}

\begin{Proof}{Proposition \ref{PropQuantileOrderOne}}
We first show that the infimum value of $M_{1,u}^{\rr,P}$ over~$\R^d$ is given by 
$$
m^* \equiv  - \int_{\R^d} \big(\ps{u}{z}+\rR(z)\big)\, dP(z) - \rR_\infty
,
$$
where we let
$$
\rR_\infty 
\equiv 
\int_0^\infty \big(1-\rr(s)\big)\, ds
\in [0,\infty]
.
$$
To this end, observe that, for any $t\geq 0$, 
$$
%		F(t)=
M_{1,u}^{\rr,P}(t u)
=
\int_{\R^d} 
\big(\rR(z-tu)-\rR(z)\big)\, dP(z) - t
.
$$
Fix $z\in\R^d$. Since $\|z-tu\|-t\to -\ps{u}{z}$ as $t\to \infty$, then the fact that $\rr(s)\to 1$ as $s\to \infty$ and the mean value theorem entail that $\rR(z-tu)-\rR(tu)\to -\ps{u}{z}$ as $t\to\infty$.
% since $\rr$ is non-negative, non-decreasing, and converges to $1$ at infinity. 
We deduce that
\begin{align*}
\rR(z-tu)-t
&=
\big(\rR(z-tu)-\rR(tu)\big) + \big(\rR(tu)-t\big)
\\[2mm]
&
=\big(\rR(z-tu)-\rR(tu)\big) + \int_0^t \big(\rr(s)-1)\, ds
\\[2mm]
% &\hspace{15mm}\to 
% -\ps{u}{z}-\int_0^\infty \big(1-\rr(s)\big)\, ds
% \\[2mm]
&\to
-\ps{u}{z}-\rR_\infty\in \big[-\infty,\|z\|\big)
,
\end{align*}
as $t\to \infty$. Because $0\leq \rr\leq 1$, we have 
$$
\rR(z-tu) 
\leq 
\int_0^{\|z\|+t} \rr(s)\, ds 
= 
\int_0^{\|z\|} \rr(s)\, ds 
+
\int_{\|z\|}^{\|z\|+t} \rr(s)\, ds 
\leq 
\int_0^{\|z\|} \rr(s)\, ds 
+
t
=
\rR(z) + t
.
$$
We deduce that $\rR(z-tu)-\RR(z)-t\le 0$ for all $z\in\R^d$ and $t\ge 0$, so that Fatou's lemma yields
\begin{equation}\label{eq:MinValM}
    \limsup_{t\to\infty} M_{1,u}^{\rr,P}(tu)
    \leq 
    -\int_{\R^d} 
    \big(\ps{u}{z}+\rR(z)\big)\, dP(z)
        -\rR_\infty
    =
    m^*
    .
\end{equation}
If $m^*=-\infty$, then the infimum value of $M_{1,u}^{\rr,P}$ over $\R^d$ is~$-\infty(=m^*)$.
%		, so that $M_{1,\tilde{u}}^P$ admits no minimum over $\R^d$. 
%		Because $m^*\leq 0$, let us then assume that $m^*\in\R$. 
If $m^* > - \infty$, then straightforward computations provide for all~$x\in\R^d$
\begin{equation}\label{eq:MinValM2}
    M^{\rr, P}_{1,u}(x)
    - 
    m^*
    =
    \int_{\R^d} \bigg( \|z-x\| + \ps{u}{z-x}
        +
        \int_{\|z-x\|}^\infty \big(1-\rr(s)\big)\, ds
        \bigg)\, dP(z)
    \geq 
    0
    .
\end{equation}
It follows that, in all cases, $m^*$ is the infimum value of $M_{1,u}^{\rr,P}$ over $\R^d$. \smallskip

Let us now prove the equivalence in the statement. Assume that assumptions (i) and (ii) of the statement hold, and let $P$ be supported on the halfline~$\LL\equiv \{y-\lambda u : \lambda\geq 0\}$ with direction~$-u$.
% and that there exists $\tau\geq 0$ such that $\rr(s)=1$ for all $s\geq \tau$. Since $\rr$ is continuous and non-decreasing, we may assume without loss of generality that $\rr(s)=1$ if and only if $s\geq \tau$, \ie $\rr(s)<1$ for all $s<\tau$ if $\tau>0$ (if $\tau=0$, then $\rr(s)=1$ for all $s\geq 0$). 
In particular, $\rR_\infty$ is finite since we thus have
$$
\rR_\infty 
=
\int_0^\tau \big(1-\rr(s)\big)\, ds
.
$$
% Observe that the collection of halflines with direction $-u$ on which $P$ is supported is totally ordered and lower-bounded (with respect to the inclusion).
Then define $\lambda_0$ as
$$
\lambda_0 
\equiv 
\sup \Big\{ \lambda \geq 0 : P[\{y-s u : 0\leq s < \lambda\}] = 0\Big\} 
.
$$
Observe that $\lambda_0<\infty$ since $P$ is supported on $\LL$. Letting $y_0\equiv y-\lambda_0 u$ we deduce that $P$ is supported on the halfine $\LL_0\equiv \{y_0-\lambda u : \lambda \geq 0\}\subset \LL$ with direction $-u$ and that $\LL_0$ is the smallest halfline with direction $-u$ (with respect to the inclusion) on which $P$ is supported. For now, assume that the integral $\int_{\R^d} (\ps{u}{z} + \rR(z))\, dP(z)$ is finite. In particular, $m^*$ is finite as well since $\rR_\infty<\infty$, and (\ref{eq:MinValM2}) thus entails that $M_{1,u}^{\rr, P}$ attains its infimum value $m^*$ at any point $x$ such that  
\begin{equation}\label{eq:CondMinM1}
\|z-x\| + \ps{u}{z-x}
+
\int_{\|z-x\|}^\infty \big(1-\rr(s)\big)
=
0
\end{equation}
for $P$-almost all $z$. By assumption (ii), this happens if and only if~$\|z-x\| + \ps{u}{z-x}=0$ for $P$-almost all $z$ and~$\|z-x\|\geq \tau$ for $P$-almost all $z$, \ie if and only if $P$ is supported on the halfline~$\{x-\lambda u : \lambda\geq 0\}$ with direction $-u$ and, by minimality of $\LL_0$, $x$ is at distance at least $\tau$ from $\LL_0$. 
% The minimality of $\LL_0^-$ thus entails that $x$ belongs to the halfline $\LL_0^+\equiv \{y_0+\lambda u : \lambda\geq 0\}$. Since $\|z-x_0\|\leq \|z-x\|$ for all $z\in\LL_0$, we have 
% $$
% \|z-x_0\| + \ps{u}{z-x_0}
% +
% \int_{\|z-x_0\|}^\infty \big(1-\rr(s)\big)
% =
% 0 
% $$
% for $P$-almost all $z$, so that $x_0$ minimizes $M_{1,u}^{\rr,P}$.
We deduce that $x\in\R^d$ minimizes $M_{1,u}^{\rr,P}$ over $\R^d$ if and only if $x=y_0+\lambda u$ for some $\lambda\geq \tau$. In particular, $M_{1,u}^P$ admits a minimum. Therefore, it remains to show that $\int_{\R^d} (\ps{u}{z}+\rR(z))\, dP(z)$ is finite. For this purpose, it is enough to show that
\begin{equation}\label{eq:BoundedIntegrand}
\sup_{z\in\LL}
 \big(\ps{u}{z}+\rR(z)\big) 
<
\infty 
.
\end{equation}
For any $z\in\mathcal{L}$, say, $z\equiv y-\lambda u$ for some $\lambda\geq 0$, write 
$$
\ps{u}{z}+\rR(z)
=
\Big(\ps{u}{y} - \lambda + \|y - \lambda u\| \Big) - \int_0^{\|y-\lambda u\|} \big(1-\rr(s)\big)\, ds
=:
\phi(\lambda)-\psi(\lambda)
.
$$
Since $\|y - \lambda u\| - \lambda \to	- \ps{u}{y}$ as $\lambda\to \infty$, we have $\lim_{\lambda\to + \infty} \phi(\lambda)=0$. Because the map $\lambda \mapsto \phi(\lambda)$ is continuous over $[0,\infty)$, we deduce that $\phi$ is bounded over $[0,\infty)$. Since $\rR_\infty<\infty$ and $\psi$ is continuous over $[0,\infty)$ and eventually monotone non-decreasing with $\psi(\lambda)\to \rR_\infty$ as $\lambda\to \infty$, we also have that $\psi$ is bounded over $[0,\infty)$. This establishes (\ref{eq:BoundedIntegrand}) and concludes this part of the proof.
\smallskip

Assume now that $M_{1,u}^{\rr,P}$ admits a minimum over $\R^d$. In particular, $m^*$ is finite.
% so that $\rR_\infty$ is finite as well, which establishes condition (i) of the statement.
%		, \ie $\int_{\R^d} (\|z\| + u'z)\, dP(z)$ exists in $\R$. 
Let $x\in\R^d$ be such that $M_{1,u}^{\rr,P}(x)=m^*$. Then (\ref{eq:MinValM2}) entails that $P$ is supported on the halfline $\LL\equiv \{x-\lambda u : \lambda\geq 0\}$ with direction $-u$, which establishes condition (i) of the statement. The same observation implies that there exists $z\in\R^d$ such that $\int_{\|z-x\|}^\infty (1-\rr(s))\, ds = 0$ which, in turn, implies that there exists $T\geq 0$ such that $\rr(s)=1$ for all $t\geq T$. Letting $\tau$ be the minimum of such $T$'s establishes condition (ii) of the statement and concludes the proof.		
\end{Proof}
\vspace{3mm}

\begin{Proof}{Proposition \ref{Propdirectionalderiv}}
    Since $\rr$ is bounded, then $\rR$ is Lipschitz over $\R^d$. In particular, the dominated convergence theorem yields 
    $$
    \lim_{t\stackrel{>}{\to} 0} \frac{M_{\alpha, u}^{\rr,P}(x+tv)-M_{\alpha, u}^{\rr,P}(x)}{t}
    =
    \int_{\R^d} \lim_{t\stackrel{>}{\to} 0} \frac{\rR(x+tv-z)-\rR(x-z)}{t}\, dP(z)
    .
    $$
    Since $\rR$ is differentiable on $\R^d\sm\{0\}$ with $\nabla \rR(x)=\rr(\|x\|)x/\|x\|$, a straightforward computation provides
    $$
    \lim_{t\stackrel{>}{\to} 0} \frac{\rR(x+tv-z)-\rR(x-z)}{t}
    =
    \begin{cases}
        \ps{\nabla \rR(x-z)}{v} & \textrm{if } z\neq x,\\[2mm]
        \|v\|\rr(0) & \textrm{if } z = x.
    \end{cases}
    ,
    $$
    which establishes the result.
\end{Proof}
\vspace{3mm}

\begin{Proof}{Theorem \ref{TheorDiffM}}
    (i) Let us start with necessity. Since $M_{\alpha,u}^{\rr,P}$ is differentiable at $x_0$, the map 
    $$
    v
    \mapsto
    \Phi(v)\equiv 
    \frac{\partial M_{\alpha,u}^{\rr,P}}{\partial v}(x_0)
    $$
    is linear. In particular we have $\Phi(v)+\Phi(-v)=0$ for all $v\in\R^d\sm\{0\}$, which entails that $\rr(0)P[\{x_0\}]=0$. Therefore, Proposition \ref{Propdirectionalderiv} entails that  $\Phi(v)=\ps{\Fpr(x_0)-\alpha u}{v}$, whence $\nabla M_{\alpha,u}^{\rr,P}(x_0)=\Fpr(x_0)-\alpha u$. Let us now turn to sufficiency, and assume that $\rr(0)P[\{x\}]=0$ for all $x\in\R^d$. Fix $x_0\in\R^d$ and define $\Phi$ as previously. In particular, we have $\Phi(v)=\ps{\Fpr(x_0)-\alpha u}{v}$ for all $v\in\R^d$. Thus, define $\nabla M_{\alpha,u}^{\rr,P}(x)=(\Phi(e_1),\Phi(e_2),\ldots, \Phi(e_d))$, where $e_i$ stands for the $i$th canonical direction of $\R^d$.
    % the vector of $\R^d$ whose $i$th component is the directional derivative of $M_{\alpha,u}^{\rr,P}$ in the $i$th canonical direction of $\R^d$ at $x$ provided in 
    % Proposition \ref{Propdirectionalderiv}, and 
    Let us show that
    \begin{equation}\label{eq:DiffM}
    \lim_{x\mapsto x_0} \frac{M_{\alpha,u}^{\rr,P}(x)-M_{\alpha,u}^{\rr,P}(x_0) - \ps{\nabla M_{\alpha,u}^{\rr,P}(x_0)}{x-x_0}}{\|x-x_0\|}
    =
    0
    .
    \end{equation}
    As in the proof of Proposition \ref{Propdirectionalderiv}, the fact that $\rr$ is bounded entails that the limit in the last display can be taken under the integrals defining $M_{\alpha,u}^{\rr,P}(x_0)$ and $\nabla M_{\alpha,u}^{\rr,P}(x_0)$. The result now follows from the differentiability of the map $x\mapsto \rR(x-z)$ over $\R^d\sm\{z\}$ combined with the fact that $\rr(0)P[\{x_0\}]=0$.\smallskip

    (ii) Part (i) of this proof entails that $\nabla M_{\alpha,u}^{\rr,P}(x)=\Fpr(x)-\alpha u$ for all $x\in U$. The conclusion follows by observing that $\Fpr$ is continuous at $x$ if and only if $P[\{x\}]=0$.
\end{Proof}
\vspace{3mm}

\begin{Proof}{Theorem \ref{TheorUniqueness}}
    (i)-(ii) Because $\rR$ is convex, $M_{1,u}^{\rr,P}$ is convex over $\R^d$: for all $x,y\in\R^d$ and $t\in [0,1]$ we have
    \begin{eqnarray*}
        \lefteqn{
            \hspace{5mm}
            M_{1,u}^{\rr,P}\big((1-t)x+ty\big)   
        }
        \\[2mm]
        &&
        =
        \int_{\R^d} \Big( \rR\big((1-t)(x-z) + t(y-z)\big) - \rR(z)\Big)\, dP(z) 
        - \ps{\alpha u}{(1-t)x+ty}
        \\[2mm]
        &&
        \leq 
         (1-t)\int_{\R^d}\Big( \rR(x-z) - \rR(z)\Big)\, dP(z) 
        +
        t \int_{\R^d} \Big( \rR(y-z) - \rR(z)\Big)\, dP(z) 
            - \ps{\alpha u}{(1-t)x+ty}
        \\[2mm]
        &&
        =
        (1-t)M_{1,u}^{\rr, P}(x)+tM_{1,u}^{\rr,P}(y)
        .
    \end{eqnarray*}
    Thus assume, \emph{ad absurdum}, that there exist $x,y\in\R^d$ and $t\in [0,1]$ such that 
    $$
    M_{1,u}^{\rr,P}\big((1-t)x+ty\big)
    =
    (1-t)M_{1,u}^{\rr, P}(x)+tM_{1,u}^{\rr,P}(y)
    .
    $$
    The convexity of $\rR$ thus entails that
    $$
    \rR\big((1-t)(x-z) + t(y-z)\big)
    =
    (1-t)\rR(x-z)+t \rR(y-z) 
    $$
    for $P$-almost all $z$. Since $\rr$ is continuous and non-decreasing, the map $g(t)\equiv  \int_0^t \rr(s)\, dt$ is differentiable over $(0,\infty)$ with non-decreasing derivative equal to $\rr$. In particular, $g$ is convex over $(0,\infty)$. Since $\rr$ is non-negative, we thus have 
    \begin{eqnarray*}
        \lefteqn{
            \hspace{-5mm}
            \rR\big((1-t)(x-z) + t(y-z)\big)
            =
            \int_0^{\|(1-t)(x-z) + t(y-z)\|} \rr(s)\, ds
            \leq 
        \int_0^{(1-t)\|x-z\|+t\|y-z\|}\rr(s)\, ds
        }
        \\[2mm]
        &&
        \leq 
        (1-t) \int_0^{\|x-z\|}\rr(s)\, ds 
        + t \int_0^{\|y-z\|} \rr(s)\, ds 
        % \\[2mm]
        % &&
        % \hspace{20mm}
        =
        (1-t)\rR(x-z) +t\rR(y-z)
        ,
    \end{eqnarray*} 
    for $P$-almost all $z$. Because all the inequalities in the previous display must be equalities, the fact that $\rr(s)>0$ for all $s>0$ entails that
    $$
    \|(1-t)(x-z)+t(y-z)\|
    =
    (1-t)\|x-z\|+t\|y-z\| 
    ,
    $$
    for $P$-almost all $z$, and $g$ is affine on the segment with endpoints $\|x-z\|$ and $\|y-z\|$. The latter implies that $g'=\rr$ is constant on this segment, which proves (i). The former implies that $x-z$ and $y-z$ are colinear for $P$-almost all $z$, which, in turn, implies that $P$ is supported on the line containing $x$ and $y$, which contradicts (ii). This establishes that for (i) and (ii), the map $M_{1,u}^{\rr,P}$ is strictly convex over $\R^d$. Since $M_{1,u}^{\rr,P}$ admits a global minimizer (see Theorem \ref{TheorExist} (ii)), this minimizer is unique.\smallskip

    (iii) Let $\LL$ be the line with direction $v$ on which $P$ is supported. We established in the first part of this proof that, for all $x,y\in\R^d$, the inequality
    \begin{equation}\label{eq:StrictConvexRemark}
    M_{1,u}^{\rr,P}\big((1-t)x+ty\big)
    <
    (1-t)M_{1,u}^{\rr, P}(x)+tM_{1,u}^{\rr,P}(y)
    \end{equation}
    is strict as soon as $P$ is not concentrated on the line containing $x$ and $y$, \ie if $x\notin\LL$ or $y\notin\LL$. By Theorem \ref{TheorExist} (ii), there exists a global minimizer $x_0$ of $M_{\alpha,u}^{\rr,P}$ over $\R^d$. If $x_0\in\R^d\sm\LL$, then (\ref{eq:StrictConvexRemark}) holds for all $y\in\R^d$ with $x=x_0$ fixed, which entails that $x_0$ is the unique minimizer of $M_{\alpha,u}^{\rr,P}$, hence proves the claim. Thus, assume that $x_0\in\LL$. First consider the case $\rr(0)>0$. Any other global minimizer of $M_{\alpha,u}^{\rr,P}$ (if it exists) must belong to $\LL$ as well. Assume, \emph{ad absurdum}, that such a minimizer $x_1\in\LL\sm\{x_0\}$ exists. For all $\lambda\in [0,1]$, let $x_\lambda\equiv (1-\lambda)x_0+\lambda x_1$. Fix $\lambda\in [0,1]$. The convexity of $M_{\alpha,u}^{\rr,P}$ entails that $x_\lambda$ is a global minimizer of $M_{\alpha,u}^{\rr,P}$, so that $x_\lambda$ satisfies the first-order condition~(\ref{eq:DirCond}). Let $Z$ be a random $d$-vector with distribution $P$ and write $Z=x_\lambda+ S_\lambda v$ for some random variable $S_\lambda$. Then, (\ref{eq:DirCond}) writes 
    \begin{equation}\label{eq:LineGradCond}
    \bigg\| \E\Big[\rr(S_\lambda){\rm sign}(S_\lambda)\I[S_\lambda\neq 0]\Big]v-\alpha u\bigg\| 
    \leq 
    \rr(0)P[\{x_\lambda\}]
    .
    \end{equation}
    Since $\alpha>0$ and $v\in \S^{d-1}\sm\{\pm u\}$, the l.h.s. of the last display cannot vanish. In particular, we have $\rr(0)P[\{x_\lambda\}]>0$. Since $\rr(0)>0$, we deduce that $x_\lambda$ is an atom of $P$ for all $\lambda\in [0,1]$, which is a contradiction since the set of atoms of $P$ is at most countable. We deduce that $x_0$ is the unique minimizer of $M_{\alpha,u}^{\rr,P}$, and satisfies (\ref{eq:LineGradCond}) so that $x_0$ is an atom of $P$. Assume now that $\rr(0)=0$. Then $x_0$ satisfies (\ref{eq:LineGradCond}) with $\lambda=0$, so that the l.h.s. of (\ref{eq:LineGradCond}) must vanish, which is not possible since $\alpha>0$ and $v\in\S^{d-1}\sm\{\pm u\}$; in particular, we cannot have $x_0\in\LL$ so that $x_0\in\R^d\sm\LL$ as was to be proved.\smallskip

    (iv) This is a special case of Proposition \ref{PropChaudInHyper}.
\end{Proof}

\subsection{Proofs for Section \ref{sec:symmetries}}

\begin{Proof}{Proposition~\ref{PropChaudInHyper}}
Fix~$x\in \R^d\sm H$ and let $x_H$ be the orthogonal projection of $x$ onto $H$. Let~$c\equiv \|x-x_H\|>0$ and~$w\equiv (x_H-x)/c$. Since $P$ is supported on $H$ we have $P[\{x\}]=0$. Proposition \ref{Propdirectionalderiv} then entails that the directional derivative of $M_{\alpha,u}^{\rr,P}$ at $x$ in direction $w$ is given by
$$
\frac{\partial M_{\alpha,u}^{\rr,P}}{\partial w}(x) 
=
\int_{\R^d} \rr(\|x-z\|)\frac{\ps{w}{x-z}}{\|x-z\|}\, dP(z) - \ps{\alpha u}{w}
.
$$
Let~$Z$ be a random $d$-vector with distribution~$P$ and define the random $k$-vector~$Y\equiv O^T(Z-x_H)$ so that~$Z=x_H+OY$. Denote by $P^Y$ the distribution of $Y$. Since $w$ is orthogonal to $u$ and $Y$, we have 
\begin{eqnarray*}
\lefteqn{
    \hspace{-25mm}
    \frac{\partial M_{\alpha,u}^{\rr,P}}{\partial w}(x) 
    =
    \int_{\R^d} \rr(\|(x_H-cw)-(x_H+Oy)\|)\frac{\ps{w}{(x_H-cw)-(x_H+Oy)}}{\|(x_H-cw)-(x_H+Oy)\|}\, dP^Y(z)
}
\\[2mm]
&&=
\int_{\R^k} \rr(\|cw + Oy\|)\frac{c}{\|cw+Oy\|}\, dP^Y(y)
<
0
.
\end{eqnarray*}
In particular, $x$ is not a minimizer of $M_{\alpha,u}^{\rr,P}$. We deduce that any minimizer $q_{\alpha,u}^\rr$ of $M_{\alpha,u}^{\rr,P}$ belongs to $H$, hence writes $q_{\alpha,u}^\rr = x_0 + O\xi_{\alpha,u}$ for some $\xi_{\alpha,u}\in\R^k$.\smallskip

(ii) We just showed that all spatial quantiles of order~$\alpha$ in direction~$u$ for~$P$ belong to~$H$, hence are of the form~$z_0+V\xi_{\alpha,u}$. Define the random $k$-vector~$X\equiv O^T(Z-x_0)$ so that~$Z=:x_0+OX$, and let $P^X$ denote the distribution of $X$. By definition, $\xi_{\alpha,u}$ thus minimizes
$$
\xi\mapsto M_{\alpha,u}^{\rr,P}(x_0+O \xi)
=
\int_{\R^d}
\Big(\| (x_0+O\xi)-z \| - \| z \| \Big)
\,
dP(z)
- \ps{\alpha u}{x_0+O\xi}
$$
$$
=
\int_{\R^k}
\Big(\|\xi-x\| - \|x_0+Ox\| \Big)\,dP^X(x)
-\ps{\alpha O' u}{\xi}
-\ps{\alpha u}{x_0}
$$
(where we used the fact that~$O'O$ is the $k\times k$ identity matrix), or, equivalently, minimizes
$$
\xi\mapsto 
\int_{\R^k}
\Big(\|\xi-x\| - \|x\| \Big)\,dP^X(x)
-\ps{\alpha O' u}{\xi}
.
$$
It follows that $\xi_{\alpha,u}$ is a $\rr$-geometric quantile of order~$\alpha$ in direction~$O^T u$ for~$P^X$, which concludes the proof.  
\end{Proof}
\vspace{3mm}

\begin{Proof}{Proposition \ref{TheorSpher}.} 
    Since $p<1$, $P$ gives positive probability to $\R^d\sm\{0\}$, so that the spherical symmetry of $P$ entails that $P$ is not supported on a single line. Theorem \ref{TheorUniqueness} (ii) thus entails that $\rr$-geometric quantiles $q_{\alpha,u}^\rr$ are unique for all $\alpha\in [0,1)$ and $u\in\S^{d-1}$. Using the uniqueness result established previously, Proposition~\ref{PropEquivariance} entails that $O q_{\alpha,u}^\rr = q_{\alpha,Ou}^\rr$ for any $d\times d$ orthogonal matrix $O$, and all $\alpha\in [0,1)$ and $u\in\S^{d-1}$. In particular, for any fixed $\alpha\in [0,1)$, the map $\S^{d-1}\to \R^d,\ u\mapsto \|q_{\alpha,u}^\rr\|$ is constant. Then, fix an arbitrary $u\in\S^{d-1}$ and define the map $\QQ:[0,1)\to [0,\infty),\ \alpha \mapsto \QQ(\alpha)\equiv \|q_{\alpha,u}^\rr\|$. The spherical symmetry of $P$ implies that $\Fpr(0)=0$. It then follows from the first-order condition (\ref{eq:DirCond}) that $q_{\alpha,u}^\rr=0$ if and only if $\alpha\in [0,\rr(0)p]$, which establishes that $\QQ(\alpha)=0$ if and only if $\alpha\in [0,\rr(0)p]$. Now fix $\alpha\in (\rr(0)p,1)$ (recall that $p<1$) and $u\in \S^{d-1}$. In particular, we have $\QQ(\alpha)>0$. Proposition~\ref{PropEquivariance} entails that for any~$d\times d$ orthogonal matrix fixing $u$, $O q_{\alpha,u}^\rr$ is a $\rr$-geometric quantile of order $\alpha$ in direction~$Ou$~($=u$), \ie $Oq_{\alpha,u}^\rr=q_{\alpha,u}^\rr$. We deduce that $q_{\alpha,u}^\rr$ belongs to the line with direction $u$, so that $q_{\alpha,u}^\rr = s\QQ(\alpha) u$ for some $s\in\{\pm 1\}$. Let us now assume that $u=(1,0,\ldots, 0)=:e_1$, without loss of generality since $q_{\alpha,Oe_1}^\rr = Oq_{\alpha,e_1}^\rr$ for any $d\times d$ orthogonal matrix. Since $\QQ(\alpha)>0$, the spherical symmetry of $P$ implies that $P[\{q_{\alpha,e_1}^\rr\}]=0$. In particular, the first-order condition (\ref{eq:DirCond}) yields $\Fpr(q_{\alpha,e_1}^\rr)=\alpha e_1$, \ie letting $Z=(Z_1,\ldots, Z_d)$ be a random $d$-vector with distribution $P$, we have
    \begin{equation}\label{eq:SphericalFirstOrder}
    \E\bigg[
    \rr\big(\|s\QQ(\alpha)e_1-Z\|\big)
    \frac{s\QQ(\alpha)e_1 - Z}{\|s\QQ(\alpha)e_1-Z\|}\I[Z\neq s\QQ(\alpha)e_1]
    \bigg]
    =
    \alpha e_1
    .
    \end{equation}
    Define $\phi:(0,\infty)\to \R$ by letting 
    $$
    \phi(\lambda)
    \equiv 
    \E\bigg[
    \rr\big(\|\lambda e_1-Z\|\big)
    \frac{\lambda - Z_1}{\|\lambda e_1-Z\|}\I[Z_1\neq \lambda]\bigg]
    ,\spa 
    \forall\ \lambda > 0
    .
    $$
    Notice that, since $P$ is spherically symmeric, we have 
    $$
    \E\bigg[
    \rr\big(\|\lambda e_1-Z\|\big)
    \frac{Z_j}{\|\lambda e_1-Z\|}\I[Z\neq \lambda e_1]
    \bigg]
    =
    0
    ,
    $$
    for all $\lambda\geq 0$. In particular, (\ref{eq:SphericalFirstOrder}) is equivalent to $\phi(s\QQ(\alpha))=\alpha$. Since $P$ has no atom, except possibily at the origin of $\R^d$, a routine application of the dominated convergence theorem entails that $\phi$ is continuous over $(0,\infty)$. Therefore, since $\lim_{\lambda\downarrow 0}\phi(\lambda) = \rr(0) p$ and $\lim_{\lambda\uparrow} \phi(\lambda) = 1$, and $\alpha\in (\rr(0)p,1)$, the Intermediate Value Theorem implies that there exists $\lambda^*\in (0,\infty)$ such that $\phi(\lambda^*)=\alpha$. In particular, as we already observed, this entails that $\Fpr(\lambda^* e_1)=\alpha u$, so that $q_{\alpha,u}^{\rr}=\lambda^* e_1$. We deduce that $s=1$ and that $q_{\alpha,u}^\rr=\QQ(\alpha)u$ (recall that these quantiles are unique). Now let us show that $\QQ$ is increasing over $(\rr(0)p,1)$. Let us first show that $\QQ$ is injective on $(\rr(0)p,1)$. Let $\alpha,\alpha'\in (\rr(0)p,1)$ and assume that $\QQ(\alpha)=\QQ(\alpha')$. For an arbitrary fixed $u\in\S^{d-1}$, this implies that $q_{\alpha,u}^\rr = q_{\alpha',u}^\rr$. Since $\QQ(\alpha)>0$ and $\QQ(\alpha')>0$, these quantiles are not atoms of $P$, so that the first-order condition (\ref{eq:DirCond}) yields $\Fpr(q_{\alpha,u}^\rr)=\alpha u$ and $\Fpr(q_{\alpha',u}^\rr)=\alpha' u$. Since $q_{\alpha,u}^\rr = q_{\alpha',u}^\rr$, we thus have $\alpha=\alpha'$, which yields injectivity of $\QQ$. Since $P$ is not supported on a line, Theorem \ref{TheorContQuant}, which is established independtly of the present proof, entails that the map $(\alpha,u)\mapsto q_{\alpha,u}^\rr$ is continuous over $[0,1)\times\S^{d-1}$. In particular, $\QQ:[0,1)\to [0,\infty)$ is continuous as well. We deduce from this, and the injectivity of $\QQ$, that $\QQ$ is either increasing over $(\rr(0)p,1)$ or decreasing over $(\rr(0)p,1)$. Since $\QQ(0)=0$ and $\QQ(\alpha)\geq 0$ for all $\alpha\in [0,1)$, we deduce that $\QQ$ is increasing over $(\rr(0)p,1)$. In particular, $\QQ$ is non-decreasing over $[0,1)$.
    % $$
    % \E\bigg[\frac{s\QQ(\alpha) - Z_1}{\|s\QQ(\alpha)e_1-Z\|}\I[Z_1\neq s\QQ(\alpha)]\bigg]
    % =
    % \alpha
    % .
    % $$
\end{Proof}

\subsection{Proofs for Section \ref{sec:QuantileMap}}

\begin{Proof}{Theorem \ref{TheorContQuant}}
    Our assumptions imply that all $\rr$-geometric quantiles of order $\alpha\in [0,1)$ exist and are unique. Therefore, using the continuity result of Lemma \ref{LemExtremesQuali1} and the coercivity obtained in (\ref{eq:AsymMalphau}), the continuity of $\Qpr$ can obtained by following the same lines as the proof of Proposition~6.1 in \cite{KonPai1}. Let us turn to injectivity of $\Fpr$. Let $x,y\in\R^d$ and assume that $v_x\equiv \Fpr(x)=\Fpr(y)\equiv v_y$. Then the equality $\|\Fpr(x)-v_y\|=0$, together with the uniqueness of quantiles, implies that $\Qpr(v_y)=x$, whereas $\|\Fpr(y)-v_x\|=0$ implies that $\Qpr(v_x)=y$. Since $v_x=v_y$, we deduce that $x=y$, which concludes the proof.
\end{Proof}
\vspace{3mm}

\begin{Proof}{Theorem \ref{TheorHomeoQuantiles}}
    Since~$\rr(0)P[\{x\}]=0$ for all $x\in\R^d$, Theorem \ref{TheorDiffM} entails, for all $x\in\R^d$ and $\alpha u\in\B^d$, that $\alpha u=\Fpr(x)$ if and only if $x$ is a $\rr$-geometric quantile of order $\alpha$ in direction $u$ for $P$, which, given the uniqueness of such quantiles (Theorem \ref{TheorUniqueness}(i)--(ii)), is equivalent to $\Qpr(\alpha u)=x$. We deduce that $\Qpr\circ \Fpr$ is the identity map on $\R^d$ and $\Fpr\circ \Qpr$ is the identity map on $\B^d$. It is easy to show that the distribution function map $\Fpr$ is continuous at a point $x\in\R^d$ if and only if $\rr(0) P[\{x\}]=0$, which is guaranteed by our assumptions.  
    % so that when $P$ is non-atomic then $\Fpr$ is continuous over $\R^d$. 
    We already proved in Theorem \ref{TheorContQuant} that the quantile map is continuous over $\B^d$, so that the proof is complete.
\end{Proof}
\vspace{3mm}

\begin{Proof}{Theorem \ref{TheorHomeoGeneral}}
    (i) For all $x\in\R^d\sm\AA$, Theorem \ref{TheorDiffM} entails that $\alpha u=\Fpr(x)$ if and only if $x$ is a $\rr$-geometric quantile of order $\alpha$ in direction $u$ for~$P$, which, given the uniqueness of such quantiles (see Theorem \ref{TheorUniqueness} (i) and (ii)), is equivalent to~$\Qpr(\alpha u)=x$. We deduce that~$\Qpr\circ \Fpr$ is the identity map on $\R^d\sm\AA$. For all $\alpha u\in\B^d$, there exists a unique $x\in\R^d$ such that $\Qpr(\alpha u)=x$, characterized by the first-order condition $\|\Fpr(x)-\alpha u\|\leq \rr(0)P[\{x\}]$ (see (\ref{eq:DirCond})). If $\alpha u\in\B^d\sm \cup_{z\in\AA} B_z$, then $P[\{x\}]=0$ (would $x$ be an atom of $P$, the first-order inequality would imply that $\alpha u\in B_x$), so that the first-order condition yields $\Fpr(x)=\alpha u$. We deduce that $\Fpr\circ \Qpr$ is the identity map on $\B^d\sm \cup_{z\in\AA} B_z$. Since $\Fpr$ is continuous precisely at those points $x\in\R^d$ such that $P[\{x\}]$, then $\Fpr$ is continuous over $\R^d\sm\AA$, so that the restriction of $\Fpr$ to $\R^d\sm\AA$ is continuous as well (with respect to the induced topology). Theorem \ref{TheorContQuant} entails that $\Qpr$ is continuous over $\R^d$, so that the restriction of $\Qpr$ to $\B^d\sm \cup_{x\in\AA} B_x$ is continuous.
    \smallskip
    
    (ii) We established in the first part of this proof that $\Fpr$ is a  bijection between $E_1\equiv \R^d\sm\AA$ and $E_2\equiv \B^d\sm \cup_{x\in \AA} B_x$ with inverse $\Qpr$. Since $\Fpr$ is injective over $\R^d$ (Theorem \ref{TheorContQuant}), we deduce that $\Fpr$ is a bijection between $E_1\cup \AA$ $(=\R^d$) and $E_2\cup \Fpr(\AA)$ $(=\B^d\sm\cup_{x\in\AA} B_x^\circ)$. It remains to observe that $\Qpr(\Fpr(x))=x$ for all $x\in\AA$. Since $\|\Fpr(x)\|<1$ (would $\|\Fpr(x)\|=1$, then the first-order condition (\ref{eq:DirCond}) would imply the existence of quantile of order~$1$, which would contradict our assumptions), there exists a unique $q\in\R^d$ such that $q$ is a $\rr$-geometric quantile of order $\|\Fpr(x)\|$ in direction $\Fpr(x)/\|\Fpr(x)\|$ (the direction can be chosen arbitrarily if $\|\Fpr(x)\|=0$), \ie we have $\Qpr(\Fpr(x))=q$. In addition, $q$ is characterized by the first-order condition $\|\Fpr(q)-\Fpr(x)\|\leq \rr(0)P[\{q\}]$. Since $\rr(0)P[\{z\}]=0$ for all $z\in\R^d$ by assumption, we deduce that $\Fpr(q)=\Fpr(x)$, which, given the injectivity result from Theorem~\ref{TheorContQuant}, yields $q=x$. In other words, we have $\Qpr(\Fpr(x))=x$ for all $x\in\AA$, which concludes the proof.
\end{Proof}

\subsection{Proofs for Section~\ref{sec:extremes}}

Throughout this section, $X$ denotes a random $d$-vector with law $P$. All limits are to be understood as $k\to\infty$ unless otherwise specified.

\begin{Lem}\label{LemExtremesQuali1}
    Let $P\in\PPP(\R^d)$ and $\rr\in\RRR$. Then, for any compact subset $K\subset \R^d$, the map 
    \[
    (\alpha,u,x)\mapsto M_{\alpha,u}^{\rr,P}(x)
    =
    \int_{\R^d} \big(\rR(x-z)-\rR(z)\big)\, dP(z)
    -
    \ps{\alpha u}{x}
    \]
    is Lipschitz over $[0,1]\times \S^{d-1}\times K$.
\end{Lem}
\vspace{1mm}
\begin{Proof}{Lemma \ref{LemExtremesQuali1}}
Fix $C>0$, and let $(\alpha_1,u_1,x_1)$ and $(\alpha_2,u_2,x_2)$ in $[0,1]\times \S^{d-1}\times [-C,C]^d$. For any $z\in \R^d$, clearly 
\[
| \rR(x_1-z)-\rR(x_2-z) | 
\le
\left| 
\int_{\|x_1-z\|}^{\|x_2-z\|} \rr(s)\, ds 
\right| 
\le
\left| \|x_1-z\| - \|x_2-z\| \right| 
\le
\|x_1-x_2\|. 
\]
Then 
\begin{align*}
| M_{\alpha_1,u_1}^{\rr,P}(x_1) - M_{\alpha_2,u_2}^{\rr,P}(x_2) | 
&\le
\|x_1-x_2\| 
+ 
\lvert \ps{\alpha_1 u_1}{x_1} - \ps{\alpha_2 u_2}{x_2} \rvert \\
&\le
C |\alpha_1-\alpha_2| + C \|u_1-u_2\| + 2\|x_1-x_2\|
,
\end{align*}
which yields the required Lipschitz continuity property.
\end{Proof}
\vspace{3mm}

% \begin{Lem}\label{LemExtremesQuali2}
%     \textcolor{blue}{Fix $r\in\RRR$ such that $r_\infty=1$ and set, for any $x\in \R^d$, 
%     \[
%     R(x)
% 	=
% 	\int_{0}^{\|x\|} r(s)\, \mathrm{d}s.
%     \]
%     Let $P\in\PPP(\R^d)$ and assume that $P$ is not supported on a half-line of $\R^d$. Then, for any $u\in \S^{d-1}$, the function 
%     \[
%     x\mapsto M_{1,u}^{\rr,P}(x)
%     =
%     \int_{\R^d} \big(R(x-z)-R(z)\big)\, \mathrm{d}P(z)
%     -
%     \ps{u}{x}
%     \]
%     does not have a minimum in $\R^d$.}
% \end{Lem}
% \vspace{1mm}
% %
% \begin{Proof}{Lemma \ref{LemExtremesQuali2}} \textcolor{blue}{Note that 
% %
% \begin{align*}
% M_{1,u}^{\rr,P}(n u) &= \int_{\R^d} \big\{R(n u-z)-(R(z)+n)\big\}\, \mathrm{d}P(z) \\
%     &= \int_{\R^d} \frac{(R(n u-z)-R(z)-n)(R(n u-z)+R(z)+n)}{R(n u-z)+R(z)+n}\, \mathrm{d}P(z).
% \end{align*}
% }
% %
% \end{Proof}
% \vspace{3mm}

\begin{Proof}{Proposition \ref{PropExtremesQuali}}
Our assumptions imply that $\rr$-geometric quantiles exist and are unique by virtue of Theorem \ref{TheorUniqueness}. Now assume, \emph{ad absurdum}, that $\|\Qpr(\alpha_k u_k)\|$ does not diverge, so that up to considering a subsequence we may assume that $q_k\equiv (\Qpr(\alpha_k u_k))$ is bounded in $\R^d$. Along a further subsequence, we may assume that $q_k\to \mu_0$ for some $\mu_0\in\R^d$ and $\alpha_k u_k\to u_0$ for some $u_0\in\S^{d-1}$. By definition, we have 
$$
M_{\alpha_k, u_k}^{\rr, P}(\Qpr(\alpha_k u_k)) 
\le 
M_{\alpha_k, u_k}^{\rr, P}(\mu)
,\spa 
\forall\ \mu\in\R^d 
.
$$
Joint continuity of $(\alpha, u, \mu)\to M_{\alpha,u}^{\rr,P}(\mu)$ (Lemma \ref{LemExtremesQuali1}) entails that 
$$
M_{1,u_0}^{\rr,P}(\mu_0) 
\le 
M_{1,u_0}^{\rr, P}(\mu) 
,\spa 
\forall\ \mu\in\R^d 
.
$$
It follows that $\mu_0$ is a global minimizer of $M_{1,u_0}^{\rr,P}$, that is $\mu_0$ is an extreme $\rr$-quantile of $P$ in direction $u_0$, which furnishes a contradiction by virtue of Proposition \ref{PropQuantileOrderOne}.
% Lemma \ref{LemExtremesQuali2} implies that this contradicts the fact that $P$ is not supported on a halfline. 
We deduce that $\|q_k\|\to \infty$. We will show that an arbitrary subsequence of $q_{k_j}$ has a further subsequence $q_\ell\equiv q_{k_{j_\ell}}$ such that $q_\ell/\|q_\ell\|\to u$, which will prove the desired result. Let $(q_{k_j})$ denote a subsequence of $(q_k)$ and $q_\ell$ be a further subsequence of $q_{k_j}$ such that $(q_\ell/\|q_\ell\|)$ converges to some $q_\infty\in \S^{d-1}$. By~(\ref{eq:DirCond}), we have $\|\Fpr(q_\ell)-\alpha_\ell u_\ell\|\le \rr(0)P[\{q_\ell\}]$ for all~$j$.
% a converging subsequence of $(\Qpr(\alpha_k u_k)/\|\Qpr(\alpha_k u_k)\|)$ and $q_{\infty}\in \S^{d-1}$ its limit.
Since $\|q_\ell\|\to \infty$, the countable additivity of probability measures combined with a subsequence argument entails that $P[\{q_\ell\}]\to 0$. 
% $(q_{k_j})$ contains a subsequence $q_\ell=q_{k_{j_\ell}}$ consisting of distinct elements. In particular, the countable additivity of probability measures implies that $P[\{q_\ell\}]\to 0$ as $\ell\to \infty$. 
It follows that 
$$
\alpha_\ell u_\ell 
=
\E\Big[ \rr(\|q_\ell -X\|) \frac{q_\ell-X}{\|q_\ell-X\|} \Ind{X\neq q_\ell} \Big]
+
o(1)
.
$$
Since $\|q_\ell\|\to\infty$ with $q_\ell/\|q_\ell\|\to q_\infty$, the dominated convergence theorem entails that the right-hand side of the preceding display converges to $q_\infty$. Since the left-hand side converges to $u$, we deduce that $q_\infty=u$, so that $q_\ell/\|q_\ell\|\to u$, which concludes the proof.
\end{Proof}
\vspace{3mm}

The proof of Theorem~\ref{TheoExtremesDirNorm} requires Lemmas~\ref{LemExtremesLinAlg}--\ref{LemExtremesNorm} below. We start with Lemma~\ref{LemExtremesLinAlg}, which is a simple result from linear algebra in Euclidean spaces.

\begin{Lem}\label{LemExtremesLinAlg}
    Fix a sequence $(u_k)\subset\S^{d-1}$ such that $u_k\to u\in \S^{d-1}$. Then for any $k$, one can construct an orthonormal basis $(v_{1,k},\ldots,v_{d-1,k})$ of the orthogonal complement of $\R u_k$ that converges to an orthonormal basis $(v_1,\ldots,v_{d-1})$ of the orthogonal complement of $\R u$.
\end{Lem}
\vspace{1mm}

\begin{Proof}{Lemma \ref{LemExtremesLinAlg}} 
Let $(e_1,e_2,\ldots,e_d)$ be the canonical basis of $\mathbb{R}^d$. Up to reshuffling the vectors in this basis, we may assume that the first coordinate of $u$ in $\mathbb{R}^d$ is nonzero. Then, for $k$ sufficiently large, $(u_k,e_2,\ldots,e_d)$ is a basis of $\mathbb{R}^d$. Apply the Gram-Schmidt procedure to this basis to obtain vectors $w_{j,k}$ defined iteratively as 
\[
w_{j,k}=e_{j+1}- \ps{e_{j+1}}{u_k} u_k - \sum_{i=1}^{j-1} \frac{\ps{e_{j+1}}{w_{i,k}} w_{i,k}}{\| w_{i,k} \|^2}, \mbox{ for } 1\leq j\leq d-1. 
\]
Letting $k\to \infty$ in the expressions of $w_{1,k},w_{2,k},\ldots,w_{d-1,k}$, we find that each sequence $(w_{j,k})$ converges to a nonzero vector $w_j$ satisfying 
\[
w_j=e_{j+1}- \ps{e_{j+1}}{u} u - \sum_{i=1}^{j-1} \frac{\ps{e_{j+1}}{w_i} w_i}{\| w_i \|^2}, \mbox{ for } 1\leq j\leq d-1. 
\]
In other words, $(w_1,\ldots,w_{d-1})$ is the basis of the orthogonal complement of $\R u$ obtained by the Gram-Schmidt process applied to $(u,e_2,\ldots,e_d)$. Set $v_{j,k}=w_{j,k}/\| w_{j,k} \|$ and $v_j=w_j/\| w_j \|$ to complete the proof.
\end{Proof}
\vspace{3mm}

Lemma~\ref{LemExtremesDir1} below is essentially a first-order asymptotic expansion of the $\rr$-quantile distribution function $\Fpr$ evaluated at an extreme $\rr$-quantile.

\begin{Lem}\label{LemExtremesDir1}
    Let $P\in\PPP(\R^d)$ and $\rr\in\RRR$. Assume that $P$ is not supported on a line or that $\rr$ is increasing. Fix $\beta>0$ and assume that $\ell_{\rr}(\beta)$ defined in (\ref{defEll}) is finite. If $P$ has a finite moment of order $\beta$, then for any sequences $(\alpha_k)\subset [0,1)$ and $(u_k)\subset\S^{d-1}$ such that $\alpha_k \to 1$ and $u_k\to u\in \S^{d-1}$, we have 
    \[
    \| \Qpr(\alpha_k u_k) \|^{\beta}\ 
    \E\bigg[ \Big(1-\rr\big(\|\Qpr(\alpha_k u_k) -X\|\big)\Big)\frac{\Qpr(\alpha_k u_k)-X}{\|\Qpr(\alpha_k u_k)-X\|} \Ind{X\neq \Qpr(\alpha_k u_k)} \bigg] 
    \to \ell_{\rr}(\beta) u.
    \]
\end{Lem}
\vspace{1mm}

\begin{Proof}{Lemma \ref{LemExtremesDir1}} 
Write $q_k$ as shorthand for $\Qpr(\alpha_k u_k)$. We have $\| q_k \|\to\infty$ and $q_k/\| q_k \|\to u$ by virtue of By Proposition~\ref{PropExtremesQuali}. Further recalling the definition of $\ell_\rr(\beta)$ from (\ref{defEll}) yields
\begin{align*}
    &\| q_k \|^{\beta} (1-\rr(\|q_k -X\|))\frac{q_k-X}{\|q_k-X\|} \Ind{X\neq q_k} \\
    &= \frac{\| q_k \|^{\beta}}{\| q_k-X \|^{\beta}} \| q_k-X \|^{\beta} (1-\rr(\|q_k -X\|))\frac{q_k-X}{\|q_k-X\|} \Ind{X\neq q_k} \to \ell_{\rr}(\beta) u
\end{align*}
with $P$-probability 1. Now, we have
\begin{align*}
    &\left\| \| q_k \|^{\beta} (1-\rr(\|q_k -X\|))\frac{q_k-X}{\|q_k-X\|} \Ind{X\neq q_k} \right\| 
    \\
    &\leq \| q_k \|^{\beta} (1-\rr(\|q_k -X\|)) 
    \Big\{ 
    \I\big[\| q_k-X \|\leq \| q_k \|/2\big] + \I\big[\| q_k-X \|>\| q_k \|/2\big]
    \Big\} 
    \\
    &\leq \| q_k \|^{\beta} \I\big[\| X \|\geq \| q_k \|/2\big] 
    + 
    \left( \frac{\| q_k \|}{\| q_k - X \|} \right)^{\beta} 
    \| q_k - X \|^{\beta} (1-\rr(\|q_k -X\|)) \I\big[\| q_k-X \|>\| q_k \|/2\big]
    \\
    &\leq 
    2^{\beta} ( \| X \|^{\beta} + C )  
\end{align*}
where $C=\sup_{s>0} s^{\beta}(1-\rr(s))<\infty$. The upper bound is an integrable random variable, so that the dominated convergence theorem applies and concludes the proof.
\end{Proof}
\vspace{3mm}

Lemma \ref{LemTaylorNorm} is Lemma 5.3 in \cite{Romon2022}. 

\begin{Lem}\label{LemTaylorNorm}
    Let $v\in\R^d\sm\{0\}$. Then, for all $h\in\R^d\sm\{v\}$ we have
    $$
    \bigg\| \frac{v-h}{\|v-h\|}-\frac{v}{\|v\|}+\frac{1}{\|v\|}\Big({\rm I}_d - \frac{vv^T}{\|v\|^2}\Big)h\bigg\| 
    \leq 
    2 \min\bigg\{ \frac{\|h\|}{\|v\|}, \frac{\|h\|^2}{\|v\|^2}\bigg\}
    .
    $$
\end{Lem}
\vspace{1mm}

Lemmas~\ref{LemExtremesDir2} and~\ref{LemExtremesNorm} below are important technical results to analyze the asymptotic behavior of the direction and of the norm of extreme $\rr$-quantile, respectively. %They are both based on the following elementary identity: if $(z_k)\subset \mathbb{R}^d$ is such that $\| z_k \| \to\infty$, then for any $x\in\R^d$ fixed,
% %
%     \begin{align}
%     \nonumber
%     \frac{z_k-x}{\|z_k-x\|}-\frac{z_k}{\|z_k\|}
%     &=
%     \frac{
%     \big(2\ps{x}{z_k}-\|x\|^2\big)z_k
%     }{
%     \|z_k\|\|z_k-x\|\big(\|z_k\|+\|z_k-x\|\big)
%     } 
%     - \frac{x}{\|z_k-x\|}
%     \\[2mm]
%     \label{LemExtremesBase}
%     &= 
%     \frac{
%     2\ps{x}{z_k}z_k-\|z_k\|\big(\|z_k\|+\|z_k-x\|\big)x
%     }{
%     \|z_k\|\|z_k-x\|\big(\|z_k\|+\|z_k-x\|\big)
%     } 
%     - \frac12 \frac{\|x\|^2}{\|z_k\|^2} \frac{z_k}{\|z_k\|} + o\Big(\frac{1}{\|z_k\|^2}\Big)
%     .
%     \end{align}
% %
% as $k\to\infty$.} 

\begin{Lem}\label{LemExtremesDir2}
    Let $P\in\PPP(\R^d)$ and assume that $P$ has a finite first moment. %Let $X$ be a random $d$-vector with law $P$. 
    Then, for any sequence
    % for any sequences $(\alpha_k)\subset [0,1)$ and $(u_k)\susbet\S^{d-1}$ such that $\alpha_k\to 1$ and $u_k\to u\in\S^{d-1}$ as $k\to\infty$, we have 
    $(z_k)\subset \mathbb{R}^d$ such that $\| z_k \| \to\infty$ and $z_k/\| z_k \|\to u\in \S^{d-1}$ for some $u\in \S^{d-1}$, we have 
    \[
    \| z_k \| \left( \E\bigg[ \frac{z_k-X}{\|z_k-X\|} \Ind{X\neq z_k} \bigg] - \frac{z_k}{\| z_k \|} \right) 
    \to 
    -\E\big[X-\ps{u}{X} u\big].
    \]
\end{Lem}
\vspace{1mm}

\begin{Proof}{Lemma \ref{LemExtremesDir2}} 
The vector-valued map $v\mapsto v/\| v\|$ is of class $C^1$ on $\R^d\sm\{0\}$ with Jacobian 
\[
\frac{1}{\|v\|}\Big({\rm I}_d - \frac{vv^T}{\|v\|^2}\Big).
\]
Consequently, for any $x\in \R^d\sm\{z_k\}$, we have
\[
\| z_k \| \left( \frac{z_k-x}{\|z_k-x\|}-\frac{z_k}{\|z_k\|}+\frac{1}{\|z_k\|}\Big({\rm I}_d - \frac{z_k z_k^T}{\|z_k\|^2}\Big)x \right) = \int_0^1 \left( \frac{(z_k-tx) (z_k-tx)^T}{\|z_k-tx\|^2} - \frac{z_k z_k^T}{\|z_k\|^2} \right) x \, \mathrm{d}t.
\]
Since $\|z_k\|\to 0$, the dominated convergence theorem entails that the last display converges to $0$ as $k\to \infty$ for all $x\in\R^d$ fixed. Then, apply Lemma \ref{LemTaylorNorm} to obtain, by dominated convergence again, 
\[
\| z_k \| \left\| \E\bigg[ \frac{z_k-X}{\|z_k-X\|} \Ind{X\neq z_k} \bigg] - \frac{z_k}{\| z_k \|} +\frac{1}{\|z_k\|}\Big({\rm I}_d - \frac{z_k z_k^T}{\|z_k\|^2}\Big) \E[X] \right\| \to 0.
\]
Finally, observe that 
\[
\Big({\rm I}_d - \frac{z_k z_k^T}{\|z_k\|^2}\Big) \E[X] \to  ({\rm I}_d - u u^T) \E[X] = \E\big[X-\ps{u}{X} u\big]
,
\]
which completes the proof.
\end{Proof}
\vspace{3mm}

\begin{Lem}\label{LemExtremesNorm}
    Let $P\in\PPP(\R^d)$ and assume $P$ has a finite second moment. %Let $X$ be a random $d$-vector with law $P$. 
    Then, for any sequence $(z_k)\subset \mathbb{R}^d$ such that $\| z_k \| \to\infty$ and $z_k/\| z_k \|\to u\in \S^{d-1}$, we have 
    % for any sequence
    % for any sequences $(\alpha_k)\subset [0,1)$ and $(u_k)\susbet\S^{d-1}$ such that $\alpha_k\to 1$ and $u_k\to u\in\S^{d-1}$ as $k\to\infty$, we have 
    \[
    \|z_k\|^2
    \ps{ \E\bigg[ \frac{z_k-X}{\|z_k-X\|} \Ind{X\neq z_k} \bigg] - \frac{z_k}{\|z_k\|}}{\frac{z_k}{\|z_k\|}}
    \to 
    -\frac12 \E\Big[\|X\|^2 - \ps{u}{X}^2\Big].
    \]
    % $(z_k)\subset \mathbb{R}^d$ such that $\| z_k \| \to\infty$ and $z_k/\| z_k \|\to u\in \S^{d-1}$ as $k\to\infty$, we have 
    % $$
    % \|z_k\|^2\ 
    % \E\bigg[
    % \ps{\alpha u - \frac{\Qpr(\alpha_k u_k)}{\|\Qpr(\alpha_k u_k)\|}}{\frac{\Qpr(\alpha_k u_k)}{\|\Qpr(\alpha_k u_k)\|}}
    % \bigg] 
    % \to 
    % \frac{1}{2}\Big(
    % $$
\end{Lem}
\vspace{1mm}

\begin{Proof}{Lemma \ref{LemExtremesNorm}} For any $x\in\R^d$ fixed and all $k$ so large that $z_k\neq x$, a direct computation provides
    \[
    \frac{z_k-x}{\|z_k-x\|}-\frac{z_k}{\|z_k\|}
    =
    \frac{
    \big(2\ps{x}{z_k}-\|x\|^2\big)z_k
    }{
    \|z_k\|\|z_k-x\|\big(\|z_k\|+\|z_k-x\|\big)
    } 
    - \frac{x}{\|z_k-x\|}.
    \]
    It follows that
    \[
    \frac{z_k-x}{\|z_k-x\|}-\frac{z_k}{\|z_k\|}
    = 
    \frac{
    2\ps{x}{z_k}z_k-\|z_k\|\big(\|z_k\|+\|z_k-x\|\big)x
    }{
    \|z_k\|\|z_k-x\|\big(\|z_k\|+\|z_k-x\|\big)
    } 
    - \frac12 \frac{\|x\|^2}{\|z_k\|^2} \frac{z_k}{\|z_k\|} + o\Big(\frac{1}{\|z_k\|^2}\Big).
    \]
    We deduce that 
    \begin{align*}
        \ps{\frac{z_k-x}{\|z_k-x\|}-\frac{z_k}{\|z_k\|}}{\frac{z_k}{\|z_k\|}}
        &=
        \frac{\big(\|z_k\|-\|z_k-x\|\big)\ps{x}{z_k}}{\|z_k\|\|z_k-x\|\big(\|z_k\|+\|z_k-x\|\big)}
        -\frac12 \frac{\|x\|^2}{\|z_k\|^2}
        + o\Big(\frac{1}{\|z_k\|^2}\Big)
        .
    \end{align*}
    Since 
    \begin{align*}
        \frac{\big(\|z_k\|-\|z_k-x\|\big)\ps{x}{z_k}}{\|z_k\|\|z_k-x\|\big(\|z_k\|+\|z_k-x\|\big)}
        &=
        \frac{2\ps{x}{z_k}-\|x\|^2}{\big(\|z_k\|+\|z_k-x\|\big)^2} \frac{\ps{x}{z_k}}{\|z_k\|\|z_k-x\|}
        ,
    \end{align*}
    and recalling that $z_k/\|z_k\|\to u$, we have 
    \[
    \|z_k\|^2 \times \frac{\big(\|z_k\|-\|z_k-x\|\big)\ps{x}{z_k}}{\|z_k\|\|z_k-x\|\big(\|z_k\|+\|z_k-x\|\big)}
    \to 
    \frac{1}{2}\ps{u}{x}^2.
    \]
    Consequently, we have
    \[
    \|z_k\|^2 \ps{\frac{z_k-x}{\|z_k-x\|}-\frac{z_k}{\|z_k\|}}{\frac{z_k}{\|z_k\|}}
    \to 
    -\frac12 \Big(\|x\|^2-\ps{u}{x}^2\Big)
    ,
    \]
    for all $x\in\R^d$. Lemma \ref{LemTaylorNorm} further entails that 
    \begin{align*}
    \bigg\lvert\ps{\frac{z_k-x}{\|z_k-x\|}-\frac{z_k}{\|z_k\|}}{\frac{z_k}{\|z_k\|}}\bigg\rvert 
    &=
    \bigg\lvert\ps{\frac{z_k-x}{\|z_k-x\|}-\frac{z_k}{\|z_k\|} + \frac{1}{\|z_k\|}\Big({\rm I}_d-\frac{z_kz_k^T}{\|z_k\|^2}\Big)x }{\frac{z_k}{\|z_k\|}}\bigg\rvert
    % \\[2mm]
    % &
    \leq 
    2\ \frac{\|x\|^2}{\|z_k\|^2}
    .
    \end{align*}
    The dominated convergence theorem then yields 
    \[
    \|z_k\|^2\ \E\bigg[\ps{\frac{z_k-X}{\|z_k-X\|} \Ind{X\neq z_k}-\frac{z_k}{\|z_k\|}}{\frac{z_k}{\|z_k\|}}\bigg]
    \to 
    -\frac12 \E\Big[\|X\|^2 - \ps{u}{X}^2\Big] 
    \]
    which concludes the proof.
\end{Proof}
\vspace{3mm}

\begin{Proof}{Theorem~\ref{TheoExtremesDirNorm}} Write $q_k$ as shorthand for $\Qpr(\alpha_k u_k)$. Lemma \ref{LemExtremesLinAlg} entails that, for each $k$, there is an orthonormal basis $(v_{1,k},\ldots,v_{d-1,k})$ of the orthogonal complement of $\R u_k$ such that each sequence $(v_{j,k})_k$ converges to some $v_j\in \S^{d-1}$ and such that $(v_1,\ldots,v_{d-1})$ forms an orthonormal basis of the orthogonal complement of $\R u$. Since $\|q_k\|\to\infty$, the first-order condition (\ref{eq:DirCond}) combined with our assumptions entails that $F_P^{\rr}(q_k)=\alpha_k u_k$ for all $k$ large enough. This rewrites 
\begin{align}
\nonumber
    &\frac{q_k}{\| q_k \|} - \alpha_k u_k \\
\label{eqn:DirNorm0}
    &= \E\bigg[ (1-\rr(\|q_k -X\|)) \frac{q_k-X}{\|q_k-X\|} \Ind{X\neq q_k} \bigg] - \left( \E\bigg[ \frac{q_k-X}{\|q_k-X\|} \Ind{X\neq q_k} \bigg] - \frac{q_k}{\| q_k \|} \right).
\end{align}
(i) Assume that $P$ has a finite first moment. Lemma~\ref{LemExtremesDir1} and Lemma~\ref{LemExtremesDir2} applied to~\eqref{eqn:DirNorm0} provides
\begin{equation}
\label{eqn:Dir1}
\| q_k \|^{\beta} \left( \frac{q_k}{\| q_k \|} - \alpha_k u_k \right) \to \ell_{\rr}(\beta) u 
\end{equation}
when $\beta<1$, and 
\begin{equation}
\label{eqn:Dir2}
\| q_k \| \left( \frac{q_k}{\| q_k \|} - \alpha_k u_k \right) \to \ell_{\rr}(1) u+\E\big[X-\ps{u}{X} u\big] 
\end{equation}
when $\beta\geq 1$. Now, observe that
\begin{align*}
\frac{q_k}{\| q_k \|} %&= \left\langle \frac{q_k}{\| q_k \|}, u_k \right\rangle u_k + \sum_{j=1}^{d-1} \left\langle \frac{q_k}{\| q_k \|}, v_{j,k} \right\rangle v_{j,k} \\
&= \left\langle \frac{q_k}{\| q_k \|}, u_k \right\rangle u_k + \sum_{j=1}^{d-1} \left\langle \frac{q_k}{\| q_k \|} - \alpha_k u_k, v_{j,k}\right\rangle v_{j,k}
.
\end{align*}
In particular, \eqref{eqn:Dir1} and~\eqref{eqn:Dir2} yield
\begin{equation}
\label{eqn:Dir3}
1 - \left\langle \frac{q_k}{\| q_k \|}, u_k \right\rangle^2 = \sum_{j=1}^{d-1} \left\langle \frac{q_k}{\| q_k \|} - \alpha_k u_k, v_{j,k}\right\rangle^2 = O\left( \frac{1}{\| q_k \|^{2\min(\beta,1)}} \right).
\end{equation}
Then write
\begin{align*}
\frac{q_k}{\| q_k \|} - u_k &= \left( \left\langle \frac{q_k}{\| q_k \|}, u_k \right\rangle - 1 \right) u_k + \sum_{j=1}^{d-1} \left\langle \frac{q_k}{\| q_k \|} - \alpha_k u_k, v_{j,k}\right\rangle v_{j,k} \\
    &= -\frac{1}{2} \left(1 - \left\langle \frac{q_k}{\| q_k \|}, u_k \right\rangle^2\right)(1+o(1)) u_k + \sum_{j=1}^{d-1} \left\langle \frac{q_k}{\| q_k \|} - \alpha_k u_k, v_{j,k}\right\rangle v_{j,k} \\
    &= \sum_{j=1}^{d-1} \left\langle \frac{q_k}{\| q_k \|} - \alpha_k u_k, v_{j,k}\right\rangle v_{j,k} + o\left( \frac{1}{\| q_k \|^{\min(\beta,1)}} \right).
\end{align*}
Using again~\eqref{eqn:Dir1} and~\eqref{eqn:Dir2}, we find that, when $\beta\leq 1$,
\begin{align*}
\| q_k \|^{\beta}(1-\alpha_k)u_k &= \| q_k \|^{\beta} \left( \frac{q_k}{\| q_k \|} - \alpha_k u_k \right) - \| q_k \|^{\beta} \left( \frac{q_k}{\| q_k \|} - u_k \right) \\
    &= \| q_k \|^{\beta} \left\langle \frac{q_k}{\| q_k \|} - \alpha_k u_k, u_k\right\rangle u_k + o\left( \frac{1}{\| q_k \|^{\min(\beta,1)}} \right) \to \ell_{\rr}(\beta) u, 
\end{align*}
namely, $\| q_k \|^{\beta}(1-\alpha_k) \to \ell_{\rr}(\beta)$.\smallskip

(ii) Assume now that $P$ has a finite second moment. Combine Lemma~\ref{LemExtremesDir1}, Lemma~\ref{LemExtremesNorm}, and~\eqref{eqn:DirNorm0} to obtain 
\begin{equation}
\label{eqn:Norm1}
\| q_k \|^{\beta} \left\langle\frac{q_k}{\| q_k \|} - \alpha_k u_k, \frac{q_k}{\| q_k \|} \right\rangle \to \ell_{\rr}(\beta)
\end{equation}
when $\beta<2$, and 
\begin{equation}
\label{eqn:Norm2}
\| q_k \|^2\left\langle\frac{q_k}{\| q_k \|} - \alpha_k u_k, \frac{q_k}{\| q_k \|} \right\rangle \to \ell_{\rr}(2)+\frac12 \E\Big[\|X\|^2 - \ps{u}{X}^2\Big] = \ell_{\rr}(2)+\frac12 \sum_{j=1}^{d-1} \E\Big[ \ps{v_j}{X}^2 \Big] 
\end{equation}
when $\beta=2$. Now, we have
\begin{align*}
\left\langle\frac{q_k}{\| q_k \|} - \alpha_k u_k, \frac{q_k}{\| q_k \|} \right\rangle &= (1-\alpha_k) \left\langle\frac{q_k}{\| q_k \|}, u_k \right\rangle + 1-\left\langle\frac{q_k}{\| q_k \|}, u_k \right\rangle \\
    &= (1-\alpha_k) \left\langle\frac{q_k}{\| q_k \|}, u_k \right\rangle + \frac{1}{2} \left( 1-\left\langle\frac{q_k}{\| q_k \|}, u_k \right\rangle^2 \right) (1+o(1)). 
\end{align*}
The second term on the right-hand side is $O( 1/\| q_k \|^{2\min(\beta,1)} )$, see~\eqref{eqn:Dir3}, so the desired conclusion in (ii) follows from~\eqref{eqn:Norm1} when $\beta<2$. When $\beta=2$, we further have
\[
\left\langle\frac{q_k}{\| q_k \|} - \alpha_k u_k, \frac{q_k}{\| q_k \|} \right\rangle = (1-\alpha_k) (1+o(1)) + \frac{1}{2} \sum_{j=1}^{d-1} \left\langle \frac{q_k}{\| q_k \|} - \alpha_k u_k, v_{j,k}\right\rangle^2 (1+o(1)). 
\]
Then, \eqref{eqn:Dir2} and~\eqref{eqn:Norm2} yield
\begin{align*}
\| q_k \|^2 (1-\alpha_k) &\to \ell_{\rr}(2) + \frac12 \sum_{j=1}^{d-1} \operatorname{Var}(\ps{X}{v_j}) \\
    &= \ell_{\rr}(2) + \frac12 \left( \ps{u}{\Sigma u} + \sum_{j=1}^{d-1} \ps{v_j}{\Sigma v_j} -  \ps{u}{\Sigma u} \right) = \ell_{\rr}(2) + \frac{1}{2} (\operatorname{tr}(\Sigma) - \ps{u}{\Sigma u})
\end{align*}
as required.
\end{Proof}

\subsection{Proofs for Section \ref{sec:Empiric}}\label{sec:ProofsEmpiric}

% \subsubsection{Proofs for Section \ref{sec:EmpiricPoint}}

\begin{Prop}\label{PropCoercivMempiric}
    Let $P\in\PPP(\R^d)$ and $\rr\in\RRR$. Fix $\alpha\in [0,1)$ and $u\in\S^{d-1}$. Let $Z_1, Z_2\ldots, $ be an infinite sample drawn from $P$. For any $n\geq 1$, let $P_n$ be the empirical distribution $n^{-1}\sum_{i=1}^n \delta_{Z_i}$ associated with $Z_1,\ldots, Z_n$, and let $M_n\equiv  M_{\alpha,u}^{\rr,P_n}$ denote the corresponding objective function. Then, for any deterministic sequence $(x_n)\subset\R^d$ such that $\|x_n\|\to\infty$ as $n\to\infty$, we have
    $$
    \liminf_{n \to \infty} \frac{M_n(x_n)}{\|x_n\|}
    \geq 
    1-\alpha
    >
    0
    ,
    $$
    with $P$-probability $1$.
\end{Prop}
\vspace{1mm}

\begin{Proof}{Proposition \ref{PropCoercivMempiric}}
    On the one hand, since $\rr\leq 1$, we have 
    $$
    \rR(x) 
    =
    \int_0^{\|x\|} \rr(s)\, ds 
    \leq 
    \|x\|
    ,
    $$
    for all $x\in\R^d$. On the other hand, the fact that $\rr$ is non-negative and non-decreasing implies that 
    $$
    \rR(x)
    \ge \int_D^{\|x\|} \rr(s)\, ds
    \ge \rr(D)(\|x\|-D)
    ,
    $$
    for all $x\in\R^d$ and $0\leq D\le \|x\|$. Then, $D>0$,  $\eta \in (0,1/2)$, and $n$ large enough so that $\|x_n\|> D$. Recalling that $\rR$ is $1$-Lipschitz since $\nabla \rR(z)=\rr(\|z\|)z/\|z\|$ satisfies $\|\nabla \rR(z)\|\leq 1$ for all $z\in\R^d\sm\{0\}$, we have
    \begin{eqnarray*}
    \lefteqn{
    \hspace{2mm}
     M_n(x_n) + \ps{\alpha u}{x_n}
    }
    \\[2mm]
    &&
    \hspace{-5mm}
    =
    \frac1n \sum_{i=1}^n \big( \rR(Z_i-x_n) - \rR(Z_i)\big) 
    \\[2mm]
    &&
    \hspace{-5mm}
    =
    \frac1n \sum_{i=1}^n \big( \rR(Z_i-x_n) - \rR(Z_i)\big)\I\big[\|Z_i\|\le \eta\|x_n\|\big]
    % \\[2mm]
    % &&
    % \hspace{20mm}
    +
    \frac1n \sum_{i=1}^n \big( \rR(Z_i-x_n) - \rR(Z_i)\big)\I\big[\|Z_i\| > \eta\|x_n\|\big]
    \\[2mm]
    &&
    \ge 
    \frac1n \sum_{i=1}^n \big( \rr(D)(\|Z_i-x_n\| - D) - \|Z_i\|\big)\I\big[\|Z_i\|\le \eta\|x_n\|\big]
    -
    \|x_n\| \frac1n \sum_{i=1}^n \I\big[\|Z_i\| > \eta\|x_n\|\big]
    \\[2mm]
    &&
    \ge 
    \|x_n\| \big(r(D)(1-\eta) - \eta)\frac1n \sum_{i=1}^n \I\big[\|Z_i\|\le \eta\|x_n\|\big]
    -
    \|x_n\| \frac1n \sum_{i=1}^n \I\big[\|Z_i\| > \eta\|x_n\|\big]
    \\[2mm]
    &&
    \ge 
    \|x_n\| \big(r(D)(1-\eta) - \eta)\frac1n \sum_{i=1}^n \I\big[\|Z_i\|\le \eta D\big]
    -
    \|x_n\| \frac1n \sum_{i=1}^n \I\big[\|Z_i\| > \eta D\big]
    .
    \end{eqnarray*}
    Consequently, using that $\liminf (a_k + b_k)\ge \liminf a_k + \liminf b_k$ for any real sequences $(a_k)$ and $(b_k)$, then the strong law of large numbers entails that 
    $$
    \liminf_{n\to\infty} \frac{M_n(x_n) + \ps{\alpha u}{x_n}}{\|x_n\|}
    \ge 
    \big(r(D)(1-\eta) - \eta) \P(\|Z_1\|\le \eta D) - \P(\|Z_1\|>\eta D)
    ,
    $$
    holds with $P$-probability $1$.
    Since $r(s)\to 1$ as $s\to\infty$, then taking $D\to \infty$ yields 
    $$
    \liminf_{n\to\infty} \frac{M_n(x_n) + \ps{\alpha u}{x_n}}{\|x_n\|}
    \ge 
    1-2\eta 
    .
    $$
    Since $\eta\in (0,1/2)$ was arbitrary, then taking $\eta\to 0$ and using that $\liminf (a_k + b_k)\le \liminf a_k + \limsup b_k$ for any real sequences $(a_k)$ and $(b_k)$, we find 
    $$
    \liminf_{n\to\infty} \frac{M_n(x_n)}{\|x_n\|}
    \ge 
    1 - \limsup_{n\to\infty} \frac{\ps{\alpha u}{x_n}}{\|x_n\|}
    \ge 
    1-\alpha 
    ,
    $$
    which concludes the proof.
\end{Proof}
\vspace{3mm}

\begin{Proof}{Theorem \ref{TheorConsistency}}
    Denote by $(\Om,\FF, \P)$ the common probability space underlying the random vectors $Z_1,Z_2,\ldots $. When $\star={\rm a.s.}$, \ie for convergence with $P$-probability $1$, we will prove the result by establishing the existence of a set $\Om^*\in\FF$ with $\P(\Om^*)=1$ such that, for all $\omega\in\Om^*$ fixed and corresponding $\hat q_n = \hat q_n(\omega)$, the (deterministic) sequence $d(\hat q_n, \MM_{\alpha,u}^P)$ converges to $0$ as $n\to\infty$. To do so, we will prove that any subsequence of $d(\hat q_n, \MM_{\alpha,u}^P)$ admits a further subsequence converging to $0$, which will establish the result. When $\star={\rm P}$, \ie for convergence in $P$-probability, we will consider subsequences for which the $o_{\rm P}(1)$ is in fact $o_{\rm a.s.}(1)$, which will reduce to the previous case.
    \medskip

    {\it Step 1. Almost sure uniform convergence over compact sets.}
    \smallskip
    
    \noindent
    Let $M_n\equiv  M_{\alpha,u}^{\rr,P_n}$ be the objective function associated with the empirical distribution $P_n$, \ie for all $x\in\R^d$
    $$
    M_n(x) 
    =
    \frac1n \sum_{i=1}^n \big( \rR(Z_i-x) - \rR(Z_i)\big) - \ps{\alpha u}{x}
    .
    $$
    For all $x\in\R^d$, the strong law of large numbers ensures that there exists $\Om_x\in\FF$ with $\P(\Om_x)=1$ such that $M_n(x)\to M_{\alpha,u}^{\rr,P}(x)$ as $n\to\infty$ on the event $\Om_x$. We will now construct an event $\Om^*\in\FF$ of full probability on which the convergence $M_n(x)\to M_{\alpha,u}^{\rr,P}(x)$ holds \textit{simultaneously} for all $x\in\R^d$ and is furthermore uniform over all compact sets. For this purpose, let $(x_k)_{k\geq 1}$ be a countable and dense subset of $\R^d$, and define $\Om^*\equiv  \cap_{k\geq 1} \Om_{x_k}$. In particular, we have $\P(\Om^*)=1$. Now fix $\omega\in\Om^*$ and consider $M_n$ on the fixed event $\omega$. Then, for all $k\geq 1$, we have $M_n(x_k)\to M_{\alpha,u}^{\rr,P}(x_k)$ as $n\to\infty$. Since $M_n$ is convex for all $n$, then Theorem 10.8 in \cite{Rock70} entails that $M_n(x)\to M_{\alpha,u}^{\rr,P}(x)$ for all $x\in\R^d$ as $n\to\infty$ and, in addition, that for any compact subset $K\subset\R^d$, we have
    \begin{equation}\label{eq:UnifComp}
    \sup_{x\in K} |M_n(x) - M_{\alpha,u}^{\rr,P}(x)| 
    \to 
    0
    ,
    \end{equation}
    as $n\to\infty$. Since this holds for any $\omega\in\Om^*$, we deduce that, with $P$-probability $1$, the convergence in (\ref{eq:UnifComp}) holds for all compact sets $K\subset\R^d$.
    \medskip
    
    {\it Step 2. The case $\star={\rm a.s.}$.}
    \smallskip 

    \noindent
    Up to further intersecting $\Om^*$ with the subset $E\in\FF$ of probability $1$ on which $o_\star(1) \to 0$ in the statement, we may assume that $o_\star(1)\to 0$ for all $\omega\in \Om^*$. Now fix $\omega\in\Om^*$ with corresponding (deterministic) sequence $\hat q_n = \hat q_n(\omega)$. Let $\hat q_{n_k}$ be an arbitrary subsequence of $\hat q_n$. We then have
    $$
    M_{n_k}(\hat q_{n_k}) 
    \leq 
    M_{n_k}(x)
    +
    o(1)
    ,
    $$
    for all $x\in\R^d$, as $k\to \infty$, where $o(1)$ is a deterministic (since $\omega\in\Om^*$ is fixed) sequence converging to $0$ when $k\to\infty$. We established above that $M_{n_k}(x)\to M_{\alpha,u}^{\rr,P}(x)$ as $k\to\infty$. In particular, the sequence $M_{n_k}(\hat q_{n_k})$ is upper-bounded. Therefore, Proposition \ref{PropCoercivMempiric} entails that $\hat q_{n_k}$ must be bounded. In particular, it admits a subsequence, which we denote $(\hat q_\ell)$, such that $\hat q_\ell\to q^*$, for some $q^*\in\R^d$, as $\ell\to\infty$; we denote by $(M_\ell)$ the corresponding subsequence of $(M_{n_k})$. Since $\hat q_\ell$ is bounded, then the uniform convergence in (\ref{eq:UnifComp}) entails that $M_\ell(\hat q_\ell)\to M_{\alpha,u}^{\rr,P}(q^*)$. Recalling that $M_\ell(\hat q_\ell)\leq M_\ell(x) + o(1)$ for all $x\in\R^d$ as $\ell\to\infty$, then taking the limit as $\ell\to\infty$ yields 
    $$
    M_{\alpha,u}^{\rr,P}(q^*)
    \leq 
    M_{\alpha,u}^{\rr,P}(x)
    ,
    $$
    for all $x\in\R^d$. Consequently, we must have $q^*\in\MM_{\alpha,u}^P$. In particular, we have $d(\hat q_\ell,\MM_{\alpha,u}^P)\to 0$ as $\ell\to\infty$. We thus proved that any subsequence of $d(\hat q_n, \MM_{\alpha,u}^P)$ admits a further subsequence converging to $0$. Since this holds for all $\omega\in\Om^*$, and $\P(\Om^*)=1$, the result follows in this case.
    \medskip 

    {\it Step 3. The case $\star={\rm P}.$}
    \smallskip

    \noindent
    We will show that $d(\hat q_n, \MM_{\alpha,u}^P) = o_{\rm P}(1)$ by, equivalently, establishing that any of its subsequences admits a further subsequence converging to $0$ with $P$-almost probability $1$. For this purpose, fix a subsequence $(\hat q_{n_k})$ of $(\hat q_n)$. We then have 
    $$
    M_{n_k}(\hat q_{n_k}) 
    \le 
    \inf_{x\in\R^d} M_{n_k}(x) + o_{\rm P}(1)
    ,
    $$
    as $k\to\infty$. For $o_{\rm P}(1)$ in the last display, we may extract a subsequence such that the corresponding $o_{\rm P}(1)$ in fact converges to $0$ with $P$-probability $1$ as $\ell\to\infty$. Denoting by $(\hat q_\ell)$ and $(M_\ell)$ the corresponding subsequences of $(\hat q_{n_k})$ and $(M_{n_k})$, respectively, we then have
    $$
    M_{\ell}(\hat q_\ell) 
    \le 
    \inf_{x\in\R^d} M_{\ell}(x) + o_{\rm a.s.}(1)
    ,
    $$
    as $\ell\to\infty$. Consequently, Step 2 implies that $d(\hat q_\ell, \MM_{\alpha,u}^P)\to 0$ wiht $P$-probability $1$. This concludes the proof.
    % \begin{equation}\label{eq:TightQhat}
    % \P(\hat q_n\in K) \to 1
    % ,
    % \end{equation}
    % as $n\to 1$. For this purpose, let $M_n\equiv M_{\alpha,u}^{\rr,P_n}$ denote the objective function associated with the empirical distribution $P_n$, \ie 
    % $$
    % M_n(x) 
    % =
    % \frac1n \sum_{i=1}^n \big( \rR(Z_i-x) - \rR(Z_i)\big) - \ps{\alpha u}{x}
    % ,
    % $$
    % for all $x\in\R^d$. Then, by definition, $q_n$ minimizes $M_n$ over $\R^d$. In addition, Lemma \ref{PropCoercivMempiric} entails that, for $\P$-almost all $\omega$ in the common probability space $(\Om, \P)$ underlying $Z_1,Z_2, \ldots$, there exists $D=D(\omega)>0$ and $N=N(\omega)\geq 1$ such that
    % $$
    % M_n(x)
    % \ge 
    % \frac{1-\alpha}{2}
    % \|x\|
    % ,
    % $$
    % holds for all $n\geq N$ and $x\in\R^d$ with $\|x\|>D$. Denote by $\Om_1$ the collection of such $\omega$'s, with $\P(\Om_1)=1$. Now assume, \textit{ad absurdum}, that for any compact set $K\subset\R^d$, we have $\liminf_{n\to\infty} \P(q_n\in K) < 1$, \ie $\limsup_{n\to\infty} \P(q_n\in K^c) > 0$. In particular, for any compact $K$, then $q_n$ admits a subsequence, $(q_\ell)$ say, and there exists $\Om_0\subset \Om$ with $\P(\Om_0)>0$, such that $q_\ell \in K^c$ for all $\ell$ and all $\omega\in \Om_0$. 
    % for $K=K_\rho$ a centered ball in $\R^d$ of radius 
    % $$
    % \rho 
    % =
    % \frac{2}{1-\alpha} \inf_{x\in\R^d} M_{\alpha,u}^{\rr,P}(x) 
    % <
    % \infty 
    % ,
    % $$
    % then there exists $\Om_0\subset \Om$ with $\P(\Om_0)>0$ such that, along a subsequence, $q_n$ belongs to $K_\rho^c$ for all $n$ and all $\omega\in \Om_1$. Note that $\P(\Om_0\cap \Om_1)>0$ since $\P(\Om_1)=1$. Then, fix $\omega\in \Om_0\cap \Om_1$. 
\end{Proof}
\vspace{3mm}

\begin{Lem}\label{LemPLine}
    Let $P\in\PPP(\R^d)$ and assume that $P$ is not supported on a single line of $\R^d$. Then 
    $$
    \sup_{\LL} P[\LL] 
    <
    1,
    $$
    where the supremum ranges over all lines of $\R^d$.
\end{Lem}
\vspace{1mm}

\begin{Proof}{Lemma \ref{LemPLine}}
    Assume, {\it ad absurdum}, that $\sup_{\LL} P[\LL]=1$. Then, there exists a sequence of lines $(\LL_n)$ of $\R^d$ such that $P[\LL_n]\to 1$ as $n\to\infty$. In particular, the sequence $(\LL_n)$ contains infinitely many distinct elements. Indeed, if it only contained finitely many distinct lines, then one would be able to extract a constant subsequence the probability of which would converge (hence, eventually, be equal) to $1$, which would contradict the fact that $P$ is not supported on a line. Consequently, up to extraction of a subsequence, we may assume that $(\LL_n)$ is an injective sequence in the sense that $\LL_n \neq \LL_{n'}$ for all $n\neq n'$. Now observe that we have, by monotonicity,
    $$
    P\Big[ \bigcap_{k} \bigcup_{n\geq k} \LL_n \Big]
    =
    \lim_{k\to \infty} P\Big[ \bigcup_{n\geq k} \LL_n \Big]
    % =
    % \limsup_{k\to \infty} P\Big[ \bigcup_{n\geq k} \LL_n \Big]
    \ge 
    \lim_{k\to \infty} P[ \LL_k ]
    =
    1
    .
    $$
    Therefore, $\cap_{k} \cup_{n\ge k} \LL_n$ is non-empty and there exists $x\in \cap_{k} \cup_{n\ge k} \LL_n$, \ie there exists a sequence $(n_k)$ such that $n_k\to \infty$ and $x\in \LL_{n_k}$ for all $k$. Now recall that $\LL_{n_k}$ are all distinct. Since the intersection of two distinct lines consist of at most one element, we deduce that $\LL_{n_k}\cap \LL_{n_{k'}}=\{x\}$, hence that $\LL_{n_k}\sm\{x\}$ is disjoint from $\LL_{n_{k'}}\sm\{x\}$, for all $k\neq k'$. Consequently, we have 
    $$
    \sum_{k} P[\LL_{n_k}\sm\{x\})]
    =
    P\Big(\bigcup_{k} (\LL_{n_k}\sm\{x\})\Big)
    \le 
    1
    .
    $$
    Since the series in the last display converges, we must have $P[\LL_{n_k}\sm\{x\}]\to 0$ as $k\to\infty$. Recalling that $P[\LL_{n_k}]\to 1$ as $k\to\infty$ thus implies that $P[\{x\}]=1$. In particular, $P$ is concentrated on a single line of $\R^d$, a contradiction. We deduce that $\sup_{\LL} P[\LL]<1$, which concludes the proof.
\end{Proof}

The following result is a standard inequality for the difference of two unit vectors; see, e.g., Lemma~A.1 in \cite{KonPai2025a} for a proof.

\begin{Lem} 
\label{LemKP22}
	For any~$v,w\in\R^d\setminus\{0\}$,
	$$
	\bigg\|\frac{v}{\|v\|}-\frac{w}{\|w\|}\bigg\|
	\leq 
	\frac{2\|v-w\|}{\|w\|}
	\cdot
	$$
\end{Lem}
\vspace{1mm}

\begin{Lem}\label{LemRApprox}
    Let $\rr\in\RRR$ be $L$-Lipschitz over $[0,\infty)$ for some $L>0$. Then, for all $x\in\R^d\sm\{0\}$ and $h\in\R^d$, we have
    $$
    \lvert\rR(x+h) - \rR(x) - \ps{\nabla \rR(x)}{h}\rvert
    \le
    \Big(3L + \frac{2 \rr(0)}{\|x\|}\Big)\|h\|^2
    ,
    $$
    and 
    $$
    \big\lvert\rR(x+h) - \rR(x) - \ps{\nabla \rR(x)}{h} - \frac12 \ps{\nabla^2 \rR(x)h}{h}\hspace{-1mm}\big\rvert
    \le
    \Big(5L + \frac{3 \rr(0)}{\|x\|}\Big)\|h\|^2
    $$
\end{Lem}
\vspace{1mm}

% \begin{Lem}\label{LemRApprox}
%     Let $\rr\in\RRR$, and assume that one of the following holds: (i) $\rr$ is $L$-Lipschitz over $(0,\infty)$ for some $L>0$, or (ii) there exists  $\delta>0$ such that $t\mapsto \rr(t)/t$ is non-increasing over $(0,\delta)$, and $\rr$ is $L$-Lipschitz over $[\delta, \infty)$ for some $L>0$. Then, for all $x\in\R^d\sm\{0\}$ and $h\in\R^d$, we have, in case (i),
%     $$
%     \|\rR(x+h) - \rR(x) - \ps{\nabla \rR(x)}{h}\| 
%     \le
%     3\max\Big\{L, \frac{\rr(\|x\|)}{\|x\|}\Big\}\|h\|^2
%     ,
%     $$
%     whereas, in case (ii), we have
%     $$
%     \|\rR(x+h) - \rR(x) - \ps{\nabla \rR(x)}{h}\| 
%     \leq 
%     3\max\Big\{L, \frac{\rr(\delta)}{\delta}, \frac{\rr(\|x\|)}{\|x\|}\Big\} \|h\|^2
%     .
%     $$
% \end{Lem}
\vspace{1mm}

\begin{Proof}{Lemma \ref{LemRApprox}}
    Fix $h\in\R^d$ and $x\in\R^d\sm D_h$, where $D_h = \{ -th : t\in [0,1] \}$, so that $x+th\neq 0$ for all $t\in [0,1]$. Let $f(t)\equiv  \rR(x+th)$ for all $t\in [0,1]$. Then, $f$ is continuous over $[0,1]$ and differentiable over $(0,1)$ so that the mean value theorem entails that 
    $$
    \rR(x+h) - \rR(x)
    =
    f(1)-f(0) 
    =
    f'(t)
    =
    \ps{\nabla \rR(x+th)}{h}
    ,
    $$
    for some $t\in (0,1)$. It follows that 
    \begin{eqnarray*}
    \lefteqn{
    \lvert\rR(x+h) - \rR(x) - \ps{\nabla \rR(x)}{h}\rvert
    =
    \lvert \ps{\nabla \rR(x+th) - \nabla \rR(x)}{h}\rvert
    }
    \\[2mm]
    &&=
    \Big\lvert \Big\langle \rr(\|x+th\|) \frac{x+th}{\|x+th\|} - \rr(\|x\|) \frac{x}{\|x\|} , h\Big\rangle \Big\rvert
    \\[2mm]
    &&
    \leq 
    \bigg( 
    \big| \rr(\|x+th\|) - \rr(\|x\|)\big| + \rr(\|x\|) \bigg\| \frac{x+th}{\|x+th\|} - \frac{x}{\|x\|} \bigg\|
    \bigg) \|h\|
    .
    \end{eqnarray*}
    Lemma \ref{LemKP22} entails that 
    $$
    \rr(\|x\|) \bigg\| \frac{x+th}{\|x+th\|} - \frac{x}{\|x\|} \bigg\|
    \leq 
    2\frac{\rr(\|x\|)}{\|x\|} \|h\| 
    .
    $$
    Since $\rr$ is Lipschitz over $[0,\infty)$, we have 
    $$
    \rr(\|x\|) \bigg\| \frac{x+th}{\|x+th\|} - \frac{x}{\|x\|} \bigg\|
    \le 
    2\Big(L + \frac{\rr(0)}{\|x\|}\Big)\|h\|
    .
    $$
    Let us now turn to $\rr(\|x+h\|)-\rr(\|x\|)$. Since $\rr$ is $L$-Lipschitz over $[0,\infty)$, we have 
    $$
    |\rr(\|x+h\|)-\rr(\|x\|)|\le L \|h\|
    ,
    $$
    which, combined with the previous inequalities, establishes the first inequality in the statement for all $h\in\R^d$ and $x\in\R^d\sm D_h$. This extends to all $x\in\R^d\sm\{0\}$ by continuity of the desired inequality in $x$ over $\R^d\sm\{0\}$.
    \smallskip 

    For the second inequality in the statement, it is enough to upper-bound $\ps{h}{\nabla^2 \rR(x)h}$. Recall from (\ref{eq:HessR}) that, for all $x\in\R^d\sm\{0\}$, we have
    $$
    \nabla^2 \rR(x)
    =
    \rr'(\|x\|) \frac{xx^T}{\|x\|^2} + \frac{\rr(\|x\|)}{\|x\|}\Big({\rm I}_d - \frac{xx^T}{\|x\|^2}\Big) 
    .
    $$
    Since $\rr$ is differentiable and $L$-Lipschitz over $(0,\infty)$, then $\rr'$ is bounded by $L$ over $(0,\infty)$. For all $x\in\R^d\sm\{0\}$ and $h\in\R^d$, we then have
    $$
    \lvert\ps{h}{\nabla^2 \rR(x)h}\rvert
    \le 
    \Big(L + \frac{\rr(\|x\|)}{\|x\|}\Big)\|h\|^2
    \le 
    \Big( 2L + \frac{\rr(0)}{\|x\|} \Big)\|h\|^2
    .
    $$
    Combined with the first part of this proof, this yields the conclusion.
\end{Proof}
\vspace{3mm}

\begin{Proof}{Theorem \ref{TheorNormal}}
Proposition \ref{TheorUniqueness}(i)--(ii) implies that $M_{\alpha,u}^{\rr,P}$ admits a unique minimizer $q_{\alpha,u}^{\rr}$ over $\R^d$. We will prove the result by applying Theorem 5.23 in \cite{van1998}, with objective function 
$$
m_\theta(x)=\rR(x-\theta)-\rR(x) - \ps{\alpha u}{\theta}
,\spa 
\forall\ \theta, x\in\R^d 
,
$$
so that $\E[m_\theta(X)]=M_{\alpha,u}^{\rr,P}(\theta)$ when $X\sim P$, and $\theta_0=q_{\alpha,u}^{\rr}$.
\medskip 

{\it Step 1. Differentiability in probability.}
\smallskip 

\noindent
Let us show that $\theta\mapsto m_\theta(x)$ is differentiable in $P$-probability at $\theta_0$, \ie for $X\sim P$ we have
\begin{equation}\label{eq:DiffProb}
\frac{m_{\theta_0+h}(X)-m_{\theta_0}(X) - \ps{\dot m_{\theta_0}(X)}{h}}{\|h\|}
=
o_{\rm P}(1)
,
\end{equation}
as $\|h\|\to 0$, for some measurable map $\dot m_{\theta_0}$. Take 
$$
\dot m_{\theta_0}(x) 
\equiv 
\nabla \rR(\theta_0-x)\I[x\neq \theta_0] - \alpha u
,\spa 
\forall\ x\in\R^d
,
$$
and let us compute 
\begin{eqnarray*}
\lefteqn{
\|h\|^{-1}\ | m_{\theta_0+h}(X)-m_{\theta_0}(X) - \ps{\dot m_{\theta_0}(X)}{h}\rvert
}
\\[2mm]
&&
=
\|h\|^{-1}\ \big\lvert \rR(X-(q_{\alpha,u}^\rr+h))-\rR(X-q_{\alpha,u}^\rr) - \langle \nabla \rR(q_{\alpha,u}^\rr-X) \I[X\neq q_{\alpha,u}^\rr] , h \rangle \big\rvert
\\[2mm]
&&\equiv
Q_h
.
\end{eqnarray*}
Now observe that, although $Q_h$ does not converge $P$-almost surely to $0$ when $q_{\alpha,u}^\rr$ is an atom of $P$, it still converges to $0$ in $P$-probability as $\|h\|\to 0$. Indeed, on the one hand, it follows from (\ref{eq:rR}), and the fact that $\rR$ is continuous over $[0,\infty)$, that the map $z\mapsto \rR(z)$ is differentiable over $\R^d\sm\{0\}$. In particular, $Q_h \I[X\neq q_{\alpha,u}^\rr] \to 0$ with $P$-probability $1$ as $\|h\|\to 0$. On the other hand, the fundamental theorem of calculus yields 
$$
\E\big[ Q_h \I[X=q_{\alpha,u}^\rr] \big]
=
\|h\|^{-1}\ \rR(h) P[\{q_{\alpha,u}^\rr\}]
\to 
\rr(0) P[\{q_{\alpha,u}^\rr\}]
=
0
,
$$
% \begin{eqnarray*}
% \lefteqn{
% \E\big[ Q_h \I[X=q_{\alpha,u}^\rr] \big]
% =
% \|h\|^{-1}\ \rR(-h) \E\big[ \I[X=q_{\alpha,u}^\rr]  \big] 
% }
% \\[2mm]
% &&
% =
% \|h\|^{-1}\ \rR(h) P[\{q_{\alpha,u}^\rr\}]
% \to 
% \rr(0) P[\{q_{\alpha,u}^\rr\}]
% =
% 0
% ,
% \end{eqnarray*}
as $\|h\|\to 0$. It follows from Markov's inequality that $Q_h\to 0$ in $P$-probability as $\|h\|\to 0$, which establishes (\ref{eq:DiffProb}).
\medskip 

{\it Step 2. Quadratic Taylor expansion.}
\smallskip 

\noindent 
Let us show that there exists a non-negative definite matrix $V_{\theta_0}$ such that
\begin{equation}\label{eq:Taylor2}
\Delta_h 
\equiv 
\E[m_{\theta_0+h}(X)] - \E[m_{\theta_0}(X)] - \frac12 h^T V_{\theta_0} h 
=
o(\|h\|^2)
,\spa 
h\to 
0
.
\end{equation}
Because $M_{\alpha,u}^{\rr,P}$ is convex on $\R^d$, then $q_{\alpha,u}^{\rr}$ is characterized by the set of equations 
$$
\frac{\partial M_{\alpha,u}^{\rr,P}}{\partial v}(q_{\alpha,u}^\rr)
\ge 
0
,\spa 
\forall\ v\in\R^d 
.
$$
Since $\rr(0)P[\{q_{\alpha,u}^\rr\}]=0$, by virtue of Assumption (ii) in the statement, then Proposition \ref{Propdirectionalderiv} entails that the last display is equivalent to 
$$
\nabla M_{\alpha,u}^{\rr,P}(q_{\alpha,u}^\rr) 
=
0
;
$$
see Definition \ref{DefinFpr} and Theorem \ref{TheorDiffM}. 
We then have
$$
\Delta_h 
=
M_{\alpha,u}^{\rr,P}(q_{\alpha,u}^\rr + h)
-
M_{\alpha,u}^{\rr,P}(q_{\alpha,u}^\rr)
-
\ps{\nabla M_{\alpha,u}^{\rr,P}(q_{\alpha,u}^\rr)}{h}
-
\frac12 h^T V_{\theta_0} h
.
$$
Letting
$$
V_{\theta_0}
\equiv 
\E\big[\nabla^2 \rR(X-q_{\alpha,u}^\rr) \I[X\neq q_{\alpha,u}^\rr] \big]
,
$$
and
\begin{align*}
J_h(x) 
&\equiv
m_{\theta_0+h}(x) - m_{\theta_0}(x) - \ps{\dot m_{\theta_0}(x)}{h}
\\[2mm]
&= \rR(x-(q_{\alpha,u}^\rr+h))-\rR(x-q_{\alpha,u}^\rr) - \langle \nabla \rR(q_{\alpha,u}^\rr-x) \I[x\neq q_{\alpha,u}^\rr] , h \rangle
,
\end{align*}
we then have
\begin{align*}
\frac{\Delta_h}{\|h\|^2} 
&=
\frac{1}{\|h\|^2} \E\Big[J_h(X) - \frac12 \ps{h}{\nabla^2 \rR(X-q_{\alpha,u}^\rr)h} \I[X\ne q_{\alpha,u}^\rr]\Big]
\\[2mm]
&=
\frac{1}{\|h\|^2} \E\Big[\big( J_h(X) - \frac12 \ps{h}{\nabla^2 \rR(X-q_{\alpha,u}^\rr)h} \big) \I[X\ne q_{\alpha,u}^\rr]\Big]
+
\frac{1}{\|h\|^2} \rR(-h) P[\{q_{\alpha,u}^\rr\}]
.
\end{align*}
Since $\rR$ is twice-differentiable over $\R^d\sm\{0\}$, we have 
$$
T_h 
\equiv 
\frac{1}{\|h\|^2}\big( J_h(X) - \ps{h}{\nabla^2 \rR(X-q_{\alpha,u}^\rr}{h} \big) \I[X\ne q_{\alpha,u}^\rr]
\to 
0
,\spa 
h\to 
0
,
$$
with $P$-probability $1$. In addition, Lemma \ref{LemRApprox} entails that 
$$
|T_h| 
\le 
5L + \frac{3\rr(0)}{\|X-q_{\alpha,u}^\rr\|}
,
$$
holds uniformly in $\|h\|\ne 0$. Notice that, under our assumptions, the r.h.s. in the last display is integrable (we have $\rr(0)=0$ when the moment condition is not satisfied). The dominated convergence theorem then entails that 
$$
\frac{1}{\|h\|^2} \E\Big[\big( J_h(X) - \frac12 \ps{h}{\nabla^2 \rR(X-q_{\alpha,u}^\rr)h} \big) \I[X\ne q_{\alpha,u}^\rr]\Big]
=
\E[ T_h ]
\to 
0.
$$
Since $(\rr(0)+\rr'(0))P[\{q_{\alpha,u}^\rr\}]=0$, we further have 
\begin{align*}
\frac{1}{\|h\|^2}\rR(-h) P[\{q_{\alpha,u}^\rr\}]
&=
\frac{1}{\|h\|^2} \int_0^{\|h\|} \rr(s)\, ds
\\[2mm]
&=
\frac{1}{\|h\|^2} \Big\{ \rr(0)\|h\| + \frac12 \rr'(0)\|h\|^2 + o(\|h\|^2) \Big\} P[\{q_{\alpha,u}^\rr\}]
\\[2mm]
&=
\frac{o(\|h\|^2)}{\|h\|^2} P[\{q_{\alpha,u}^\rr\}]
\to 
0
,
\end{align*}
as $h\to 0$. We deduce that $\|h\|^{-2} \Delta_h \to 0$ as $h\to 0$, which establishes (\ref{eq:Taylor2}).
\medskip 

{\it Step 3. Invertibility of the Hessian.} 
\smallskip 

\noindent 
Since $\rr(s)>0$, we see from the expression for $\nabla^2 \rR$ from (\ref{eq:HessR}) that if there exists $v\in\R^d\sm\{0\}$ such that $\ps{v}{V_{\theta_0}v}=0$, then 
$$
1-\frac{\ps{v}{z-q_{\alpha,u}^\rr}}{\|z-q_{\alpha,u}^\rr\|}
$$
for $P$-almost all $z\in\R^d$, \ie $P$ is supported on a single line. This contradicts our assumption if $P$ is not supported a line. Assume then, that $P$ is supported on a line and, hence, that $\rr'(s)>0$ for all $s>0$ and $P$ is not a Dirac mass. It follows again from the expression of $\nabla^2 \rR$ and the previous display that $\rr'(\|z-q_{\alpha,u}^\rr\|)=0$ for $P$-almost all $z$ which, given that $P$ is not a Dirac, further implies that $\rr'(s)=0$ for some $s>0$, a contradiction. We deduce that $\ps{v}{V_{\theta_0}v}>0$ for all $v\in\R^d\sm\{0\}$ which, given that $V_{\theta_0}$ is symmetric, entails that $V_{\theta_0}$ is invertible.
\medskip 

{\it Step 4. Conclusion.} 
\smallskip 

\noindent 
Because $0\le \rr\leq 1$, then $\rR$ is Lipschitz over $\R^d$, so that the map $\theta\mapsto m_\theta(x)$ is Lipschitz over $\R^d$ for every $x\in\R^d$ fixed with Lipschitz constant uniformly bounded over $x\in\R^d$. Since 
$$
M_{\alpha,u}^{\rr,P_n}(\hat q_n)
\le 
\inf_{y\in\R^d} M_{\alpha,u}^{\rr,P_n}(y) 
,
$$
and $\hat q_n = q_{\alpha,u}^\rr + o_{\rm P}(1)$ by Theorem \ref{TheorConsistency}, the conclusion follows from Theorem 5.23 in \cite{van1998}.
% For this purpose, by virtue of Markov's inequality, it is enough to identify a map $\dot m_{\theta_0}$ such that 
% $$
% \|h\|^{-1}\ \E\big[| m_{\theta_0+h}(X)-m_{\theta_0}(X) - \ps{\dot m_{\theta_0}(X)}{h}\rvert \big]
% =
% o(1)
% ,
% $$
% as $\|h\|\to 0$. Take 
% $$
% \dot m_{\theta_0}(x) 
% \equiv 
% \nabla \rR(\theta_0-x)\I[x\neq \theta_0] - \alpha u
% ,\spa 
% \forall\ x\in\R^d
% ,
% $$
% and let us compute 
% \begin{eqnarray*}
% \lefteqn{
% \|h\|^{-1}\ \E\big[| m_{\theta_0+h}(X)-m_{\theta_0}(X) - \ps{\dot m_{\theta_0}(X)}{h}\rvert \big]
% }
% \\[2mm]
% &&
% =
% \|h\|^{-1}\ \E\Big[\big\lvert \rR(X-(\theta_0+h))-\rR(X-\theta_0) - \langle \nabla \rR(\theta_0-X) , h \rangle \big\rvert \Big]
% \end{eqnarray*}
% It follows from (\ref{eq:rR}) and the fact that $\rR$ is continuous over $[0,\infty)$ that the map $z\mapsto \rR(z)$ is differentiable over $\R^d\sm\{0\}$, so that $\theta\mapsto m_\theta(x)$ is differentiable over 
\end{Proof}

\subsection{Proof of Theorem \ref{TheorSummary}(vi)}\label{sec:ProofsSummary}

% \gs{This would deserve a theorem in the main text (it is the only result in Theorem \ref{TheorSummary} which is not taken from the main text I think).}

Before proceeding with the proof, let us introduce a collection $\FF\subset \PPP(\R^d)$ of Borel probability measures such that for any compact set $C\subset\R^d$ there exists $a>0$ and a bounded Borel subset $E\subset\R^d$ such that 
\begin{equation}\label{eq:InfCond}
\inf_{P\in \FF} \inf_{x\in C}
% \big(
P[E \sm B(x,a)]
% -
% \sup_{x\in C} P[B(x,a)] \big)
>
0
,
\end{equation}
where $B(x,a)$ denotes the open ball centred at $x$ with radius $a$. We exhibit two important examples of such classes $\FF$. 
\begin{itemize}
    \item[a)] Consider fixed distinct locations $x_1,\ldots, x_n\in\R^d$ with $n\ge 3$. Then, take $\FF$ as the collection of discrete probability measures $P=n^{-1}\sum_{i=1}^n \delta_{y_i}$ such that at most $\ell$ of the $y_i$'s differ from the $x_i$'s, for some $\ell \le n-2$. Then fix a compact set $C\subset\R^d$, and take $a=\min_{i,j} \|x_i-x_j\|>0$ and $E=\cup_{x\in C} B(x,a) \cup \{x_1,\ldots, x_n\}$. Fix $P\in \FF$ and write $P=n^{-1}(\sum_{i=1}^\ell \delta_{y_i} + \sum_{i=\ell+1}^n \delta_{x_i})$. Then, $P[E\sm B(x,a)]$ is minimal when $B(x,a)$ contains the maximum number of points possible from $\{y_1,\ldots, y_\ell, x_{\ell+1}, \ldots, x_n\}$. Assume that a number $k\le \ell$ of the $y_i$'s lies in $E$. Since any $B(x,a)$ contains at most one of the $x_i$'s, we deduce that 
    $$
    P[E\sm B(x,a)]
    \ge 
    \frac{(n-\ell+k)-(k+1)}{n}
    =
    \frac{n-\ell-1}{n}
    ,
    $$
    which is positive by virtue of $\ell\le n-2$. The lower bound is also uniform in $x$, $C$ and $P$.

    \item[b)] Consider a fixed $P\in\PPP(\R^d)$ and $\FF=\{P\}$. Fix a compact set $C\subset \R^d$. Unless $P$ is a single Dirac mass, we let $\ve=1-\sup_{x\in\R^d} P[\{x\}]>0$. Then, take a compact set $E\subset\R^d$ such that $P[E]\ge 1-\ve/2$ and assume, \emph{ad absurdum}, that $\inf_{x\in C} P[E\sm B(x,a)]=0$ for all $a>0$. Then take a sequence $a_n\to 0$ and $(x_n)\subset C$ such that $P[E\sm B(x_n, a_n)]\to 0$. Up to a subsequence, we may assume that $x_n\to x$ for some $x\in\R^d$ since $C$ is compact. For all $r>0$ and $n$ large enough such that $a_n + |x_n-x| < r$, we have 
    $$
    P[E\sm B(x_n, a_n)]
    \ge 
    P[E]-P[B(x_n,a_n)] 
    \ge 
    P[E] - P[B(x,r)]
    .
    $$
    We deduce that 
    $$
    \liminf_{n\to \infty} P[E\sm B(x_n, a_n)]
    \ge 
    P[E] - P[B(x,r)]
    ,\spa 
    \forall\ r>0 
    .
    $$
    We deduce that the r.h.s. in the last display can be replaced by $P[E]-P[\{x\}]$, which yields 
    $$
    \liminf_{n\to \infty} P[E\sm B(x_n, a_n)] 
    \ge 
    1-\frac{\ve}{2} - P[\{x\}]
    \ge 
    \ve - \frac{\ve}{2}
    >
    0
    ,
    $$
    a contradiction.
\end{itemize}\smallskip
% $ for all $x\in C$, we have
% \begin{eqnarray*}
% \lefteqn{
% P[E\sm B(x,a)]
% \ge 
% P\big[(E\sm B(x,a))\sm \{y_1,\ldots, y_\ell\}]
% =
% }
% \\[2mm]
% &&
% \frac{1}{n} \sum_{i=\ell+1}^n \I\big[x_i\in (E\sm B(x,a))\sm \{y_1,\ldots, y_\ell\}\big]
% =
% \frac{1}{n} \sum_{i=\ell+1}^n \Big( 1 - \I\big[x_i\in B(x,a)\cup \{y_1,\ldots, y_\ell\}\big] \Big)
% .
% \end{eqnarray*}
% Since there are at most $\ell+1$ points $x_i$'s satisfying $x_i\in B(x,a)\cup\{y_1,\ldots, y_\ell\}$, we find 
% $$
% P[E\sm B(x,a)]
% \ge 
% \frac{n-\ell-(\ell+1)}{n}
% =
% \frac{n-1-2\ell}{n}
% >
% 0
% ,
% $$
% since $\ell < (n-1)/2$.
% $$
% \\[2mm]
% &&
% \frac{1}{n} \sum_{i=\ell+1}^n \Big( 1 - \I\big[x_i\in B(x,a)\cup \{y_1,\ldots, y_\ell\}\big] \Big)
% \end{eqnarray*}
% $$
% \inf_{P\in\FF} P [E]
% =
% \frac{n-\ell}{n} 
% >
% \frac{\ell+1}{n}
% \ge
% \sup_{x\in C} \sup_{P\in \FF} P[B(x,a)]
% .
% $$

We now proceed with the proof of Theorem~\ref{TheorSummary}(vi). Let $\FF\subset\PP(\R^d)$ be a collection satisfying~(\ref{eq:InfCond}) and fix $P\in\FF$. To ease notation, we write $q=\Qpr(\alpha u)$. It follows from Step 2 in the proof of Theorem~\ref{TheorNormal} that $M_{\alpha,u}^{\rr,P}$ is twice-differentiable over $\R^d$ with 
$$
\nabla^2 M_{\alpha,u}^{\rr,P}(x) 
=
\int_{\R^d} (\nabla^2 \rR)(x-z) \I[z\ne x]\, dP(z)
,\spa 
\forall\ x\in\R^d 
.
$$
Moreover, since $\rr'$ and $\rr$ are bounded, and $\rr(0)=0$, a simple application of the dominated convergence theorem entails that the Hessian $\nabla^2 M_{\alpha,u}^{\rr,P}$ is continuous. Straightforward calculations provide
$$
M_{\alpha,u}^{\rr,P}(x) - M_{\alpha,u}^{\tilde\rr, P}(x) 
=
\int_0^\infty P[\R^d \sm B(x,s)]\big( \rr(s)-\tilde\rr(s)\big)\, ds 
,\spa 
\forall\ x\in\R^d 
. 
$$
It follows that
$$
\sup_{x\in\R^d} \big| M_{\alpha,u}^{\rr,P}(x) - M_{\alpha,u}^{\tilde\rr, P}(x) \big| 
\le 
\int_0^\infty |\rr(s)-\tilde\rr(s)|\, ds 
\equiv
\delta 
.
$$
We can assume, without loss of generality, that $\delta<\infty$ since the statement is otherwise trivial. Now, on the one hand, since $\tilde q$ is a global minimizer of $M_{\alpha,u}^{\tilde r,P}$, we have
$$
M_{\alpha,u}^{\rr,P}(\tilde q)
\le 
M_{\alpha,u}^{\tilde \rr,P}(\tilde q) + \delta 
\le 
M_{\alpha,u}^{\tilde \rr,P}(q) + \delta 
\le 
M_{\alpha,u}^{\rr,P}(q) + 2\delta 
.
$$
On the other hand, since $q$ is a global minimum of $M_{\alpha,u}^{\rr,P}\in C^2$, then $\nabla M_{\alpha,u}^{\rr,P}(q)=0$ and Taylor's theorem provides 
$$
M_{\alpha,u}^{\rr,P}(\tilde q)
- 
M_{\alpha,u}^{\rr,P}(q)
=
\int_0^1 (1-t) 
\ps{
\tilde q-q
}{
\nabla^2 M_{\alpha,u}^{\rr,P}\big(q + t(\tilde q-q)\big) (\tilde q-q)
}\, dt
.
$$
Now, let $K\subset \R^d$ be a compact set containing the line segment joining $q$ and $\tilde q$. In fact, since $\Qpr$ is continuous over $\B^d$ by virtue of Theorem \ref{TheorContQuant} since $\rr$ is increasing, we can take $K$ large enough so that it contains the compact set $\Qpr(I\times\S^{d-1})$, hence only depends on $I$ and $\|\tilde q\|$. This provides
$$
\frac12 \|\tilde q-q\|^2 \inf_{x\in K} \inf_{\|v\|=1} \ps{v}{\nabla^2 M_{\alpha,u}^{\rr,P}(x)v}
\le 
M_{\alpha,u}^{\rr,P}(\tilde q)
- 
M_{\alpha,u}^{\rr,P}(q)
\le 
2\delta 
.
$$
Given the expression of $\nabla^2 \rR(x)$ in (\ref{eq:rR}), and using that $\rr$ is strictly increasing and $\rr'>0$ over $(0,\infty)$, we find for all $0<a<b$ and $\|v\|=1$
\begin{eqnarray*}
\lefteqn{
\ps{v}{\nabla^2 M_{\alpha,u}^{\rr,P}(x)v}
\ge 
% \frac{\rr(A)}{B} 
\min_{s\in [a,b]} \frac{\rr(s)}{s}
\int_{a\le \|x-z\|\le b}
 \Big\{ 1 - \Big( \frac{\ps{v}{x-z}}{\|x-z\|}\Big)^2\Big\} \, dP(z)
}
\\[2mm]
&&\hspace{30mm}
+
\min_{s\in [a,b]} \rr'(s) \int_{a\le \|x-z\|\le b} \Big( \frac{\ps{v}{x-z}}{\|x-z\|}\Big)^2 \, dP(z)
\\[2mm]
&&\hspace{22mm}
\ge 
\min\Big(\min_{s\in [a,b]} \frac{\rr(s)}{s}, \min_{s\in [a,b]} \rr'(s)\Big) P\big[B(x,b)\sm B(x,a)\big]
.
\end{eqnarray*}
Let $a>0$ and $E\subset\R^d$ be a bounded measurable set such that (\ref{eq:InfCond}) holds with compact set $C=K$. Then, take $b$ large enough such that $B(x,b)\supset E$ for all $x\in K$. This provides 
$$
\inf_{x\in K} \inf_{\|v\|=1} \ps{v}{\nabla^2 M_{\alpha,u}^{\rr,P}(x)v}
\ge 
c P[E\sm B(x,a)]
\ge 
c \inf_{x\in K} \inf_{P\in \FF} P[E\sm B(x,a)]
\equiv 
C
>0
,
$$
where $c=c(a,b,\rr)>0$ since $\rr>0$ and $\rr'>0$ over $(0,\infty)$. We deduce that
$$
\|\tilde q-q\|^2 
\le 
4C \delta 
=
4C \int_0^\infty |\rr(s)-\tilde\rr(s)|\, ds 
.
$$
The constant $C$ only depends on the regularizer $\rr$, the class of distributions $\FF$, and the compact set $K$. We deduce that $C$ depends only on $\rr$, $\FF$, $I$, and $\|\tilde q\|$, which concludes the proof.

% \subsubsection{Proofs for Section \ref{sec:EmpiricDonsker}}

\end{document}